\documentclass{gtart_a}
\pdfoutput=1
\usepackage[all]{xy}
\usepackage{calc,dratex,stmaryrd}

\title[Blanchfield and Seifert algebra]{Blanchfield and Seifert 
algebra in high-dimensional\\boundary link theory I: Algebraic $K$--theory}

\author{Andrew Ranicki}
\givenname{Andrew}
\surname{Ranicki}
\address{School of Mathematics\\University of Edinburgh\\\newline
Edinburgh EH9 3JZ\\Scotland, UK}
\email{a.ranicki@ed.ac.uk}
\urladdr{}

\author{Desmond Sheiham}
\givenname{Desmond}
\surname{Sheiham}

\volumenumber{10}
\issuenumber{}
\publicationyear{2006}
\papernumber{43}
\lognumber{0664}
\startpage{1761}
\endpage{1853}

\doi{}
\MR{}
\Zbl{}


\received{6 October 2005}
\revised{14 July 2006}
\accepted{2 September 2006}
\proposed{Wolfgang L\"uck}
\seconded{Peter Teichner, Steve Ferry}
\published{2 November 2006}
\publishedonline{2 November 2006}
\corresponding{}
\editor{}
\version{}

\arxivreference{math.AT/0508405}

\let\xysavmatrix\xymatrix
\def\xymatrix{\disablesubscriptcorrection\xysavmatrix}
\AtBeginDocument{\let\bar\wbar}

\newcommand{\llangle}{\langle\!\langle}
\newcommand{\rrangle}{\rangle\!\rangle}

\makeop{Hom}
\makeop{Tor}
\makeop{coker}
\makeop{adj}
\makeop{min}
\makeop{max}

\makeatletter
\def\cnewtheorem#1[#2]#3{\newtheorem{#1}{#3}[section]
\expandafter\let\csname c@#1\endcsname\c@thm}

\newtheorem{thm}{Theorem}[section]
\newtheorem{athm}{Theorem}

\cnewtheorem{lemma}[thm]{Lemma}
\cnewtheorem{cor}[thm]{Corollary}
\cnewtheorem{proposition}[thm]{Proposition}
\theoremstyle{remark}
\cnewtheorem{remark}[thm]{Remark}
\cnewtheorem{example}[thm]{Example}
\cnewtheorem{definition}[thm]{Definition}
\cnewtheorem{notation}[thm]{Notation}
\newtheorem*{rem}{Remark}

\makeatother  %  move after \newtheorem block

\makeautorefname{athm}{Theorem}

\newcommand{\A}{\mathcal{A}}
\newcommand{\B}{\mathcal{B}}
\newcommand{\Bla}{{\mathcal{B}\mathit{la}}}
\newcommand{\BLA}{{\rm Bla}}
\renewcommand{\C}{\mathcal{C}}
\renewcommand{\D}{\mathcal{D}}
\newcommand{\E}{\mathcal{E}}
\newcommand{\Endd}{{\mathcal{E}\mathit{nd}}}
\newcommand{\ENDD}{{\rm End}}
\newcommand{\F}{\mathcal{F}}
\newcommand{\Flk}{{\mathcal{F}\mathit{lk}}}
\newcommand{\FLK}{{\rm Flk}}
\newcommand{\G}{\mathcal{G}}
\newcommand{\Mod}{{\mathcal{M}\mathit{od}}}
\newcommand{\Nil}{{\mathcal{N}\mathit{il}}}
\newcommand{\NIL}{{\rm Nil}}
\newcommand{\PP}{{\mathcal{P}}}
\newcommand{\Prim}{{\mathcal{P}\mathit{rim}}}
\newcommand{\PRIM}{{\rm Prim}}
\newcommand{\Proj}{{\mathcal{P}\mathit{roj}}}

\newcommand{\RR}{\mathbb{R}}
\newcommand{\Sei}{{\mathcal{S}\mathit{ei}}}
\newcommand{\SEI}{{\rm Sei}}

\renewcommand{\Z}{\mathbb{Z}}
\newcommand{\di}{\displaystyle}
\DeclareMathOperator{\End}{End}

\begin{document}

\begin{asciiabstract}
The classification of high-dimensional mu-component boundary links
motivates decomposition theorems for the algebraic K-groups of the
group ring A[F_mu] and the noncommutative Cohn localization
Sigma^{-1}A[F_mu], for any mu>0 and an arbitrary ring A, with F_mu the
free group on mu generators and Sigma the set of matrices over A[F_mu]
which become invertible over A under the augmentation A[F_mu] to A.
Blanchfield A[F_mu]-modules and Seifert A-modules are abstract
algebraic analogues of the exteriors and Seifert surfaces of boundary
links.  Algebraic transversality for A[F_mu]-module chain complexes is
used to establish a long exact sequence relating the algebraic
K-groups of the Blanchfield and Seifert modules, and to obtain the
decompositions of K_*(A[F_mu]) and K_*(Sigma^{-1}A[F_mu]) subject to a
stable flatness condition on Sigma^{-1}A[F_mu] for the higher
K-groups.
\end{asciiabstract}

\begin{htmlabstract}
The classification of high-dimensional &mu;&ndash;component boundary
links motivates decomposition theorems for the algebraic K&ndash;groups
of the group ring A[F<sub>&mu;</sub>] and the noncommutative Cohn
localization &Sigma;<sup>-1</sup>A[F<sub>&mu;</sub>], for any &mu;&ge;1
and an arbitrary ring A, with F<sub>&mu;</sub> the free
group on &mu; generators and &Sigma; the set of matrices over
A[F<sub>&mu;</sub>] which become invertible over A under the augmentation
A[F<sub>&mu;</sub>]&rarr;A.  Blanchfield A[F<sub>&mu;</sub>]&ndash;modules
and Seifert A&ndash;modules are abstract algebraic analogues of
the exteriors and Seifert surfaces of boundary links.  Algebraic
transversality for A[F<sub>&mu;</sub>]&ndash;module chain complexes
is used to establish a long exact sequence relating the algebraic
K&ndash;groups of the Blanchfield and Seifert modules, and to
obtain the decompositions of K<sub>*</sub>(A[F<sub>&mu;</sub>]) and
K<sub>*</sub>(&Sigma;<sup>-1</sup>A[F<sub>&mu;</sub>]) subject to a
stable flatness condition on &Sigma;<sup>-1</sup>A[F<sub>&mu;</sub>]
for the higher K&ndash;groups.
\end{htmlabstract}

\begin{abstract}
The classification of high-dimensional $\mu$--component boundary
links motivates decomposition theorems for the algebraic $K$--groups
of the group ring $A[F_{\mu}]$ and the noncommutative Cohn localization
$\Sigma^{-1}A[F_{\mu}]$, for any $\mu \geqslant 1$ and an arbitrary ring
$A$, with $F_{\mu}$ the free group on $\mu$ generators and $\Sigma$ the
set of matrices over $A[F_{\mu}]$ which become invertible over $A$ under
the augmentation $A[F_{\mu}] \to A$.  Blanchfield $A[F_{\mu}]$--modules
and Seifert $A$--modules are abstract algebraic analogues of the exteriors
and Seifert surfaces of boundary links.  Algebraic transversality for
$A[F_{\mu}]$--module chain complexes is used to establish a long exact
sequence relating the algebraic $K$--groups of the Blanchfield and
Seifert modules, and to obtain the decompositions of $K_*(A[F_{\mu}])$
and $K_*(\Sigma^{-1}A[F_{\mu}])$ subject to a stable flatness condition
on $\Sigma^{-1}A[F_{\mu}]$ for the higher $K$--groups.
\end{abstract}

\maketitle
\cl{\small\it Desmond Sheiham died 25 March 2005.}
\cl{\small This paper is dedicated to the memory of Paul Cohn and Jerry Levine.}

\section*{Introduction}

For any integer $\mu \geqslant 1$ let $F_\mu$ be the free group on
$\mu$ generators $z_1,z_2,\ldots,z_\mu$.  The classification theory of
high-dimensional $\mu$--component boundary links involves `Seifert
$\Z$--modules' and `Blanchfield $\Z[F_\mu]$--modules', corresponding to
the algebraic invariants obtained from $\mu$--component Seifert surfaces
and the boundary link exterior.  This paper concerns the algebraic
relationship between f.g.~projective Seifert $A$--modules and h.d.~1
Blanchfield $A[F_\mu]$--modules for any ring $A$, extending the work of
Sheiham \cite{Sh2}.  Part I deals with the algebraic $K$--theory of the
Seifert and Blanchfield modules.  Part II will deal with the algebraic
$L$--theory of the Seifert and Blanchfield forms, such as arises in the
computation of the cobordism groups of boundary links.
The algebraic $K$-- and $L$--theory in the knot case $\mu = 1$
have already been done by Ranicki \cite{RBS}.

\subsection*{Combinatorial transversality}
\fullref{combinatorial transversality} develops a combinatorial
construction of fundamental domains for $F_\mu$--covers of $CW$ complexes
which will serve as a role model for the algebraic transversality
of $A[F_\mu]$--module chain complexes in the subsequent sections.
The $F_\mu$--covers $p\co \wwtilde{W}{\to}W$ of a space $W$ are
classified by the homotopy classes of maps
$$c\co W \longrightarrow BF_\mu = {\bigvee_\mu} S^1$$
with $\wwtilde{W} = c^*EF_\mu$ the pullback to $W$ of the universal cover
$EF_\mu$ of $BF_\mu$.
Let $0 \in BF_\mu$ be the point at which the circles $S^1$ are joined,
and choose points $1,2,\ldots,\mu \in BF_\mu\backslash \{0\}$, one in
each circle.  If $W$ is a compact manifold then $c$ is homotopic to a map which
is transverse regular at $\{1,2,\ldots,\mu\} \subset BF_\mu$, so that
$$V  =  c^{-1}(\{1,2,\ldots,\mu\})  =  V_1 \sqcup V_2 \sqcup \ldots \sqcup V_\mu
 \subset W$$
is a disjoint union of $\mu$ codimension--1 submanifolds
$V_i = c^{-1}(\{i\})\subset W$ (which may be empty) and cutting $W$ at
$V$ there is obtained a fundamental domain $U \subset
\wwtilde{W}$, a compact manifold with boundary
$$\partial U = \bigsqcup\limits_{i = 1}^{\mu}(V_i \sqcup z_iV_i).$$
If $W$ is connected and $c_*\co \pi_1(W) \to F_\mu$ is surjective then $U$
is connected and $V_1,V_2,\ldots,V_\mu$ are non-empty, and may be chosen
to be connected.  In the combinatorial version of transversality it is
only required that $W$ be a finite $CW$ complex, and $W$ may be
replaced by a simple homotopy equivalent finite $CW$ complex also
denoted by $W$, with disjoint subcomplexes $V_1,V_2,\ldots,V_\mu \subset
W$ and a fundamental domain $U \subset \wwtilde{W}$ which is a finite
subcomplex with a subcomplex
$$\partial U = \bigsqcup\limits_{i = 1}^{\mu}(V_i \sqcup z_iV_i) \subset U,
\quad V_i = U \cap z^{-1}_iU$$
such that
$$\bigcup\limits_{g \in F_\mu}gU = \wwtilde{W},~
gU \cap hU = \emptyset \text{ unless}~g^{-1}h \in \{1,z_1,z_1^{-1},
\ldots,z_{\mu},z_{\mu}^{-1}\}.$$
Ranicki \cite{RAC} developed  combinatorial transversality at $Y \subset X$
for maps of finite $CW$ complexes
$$W\to X = X_1\cup_YX_2$$
with $X,X_1,X_2,Y$ connected and $\pi_1(Y) \to \pi_1(X_1)$, $\pi_1(Y) \to
\pi_1(X_2)$ injective. The essential difference from \cite{RAC} is that we
are here using the Cayley tree $EF_\mu = G_\mu$ of $F_{\mu}$ rather than
the Bass--Serre tree of the amalgamated free product given by the
Seifert--van Kampen Theorem
$$\pi_1(X) = \pi_1(X_1)*_{\pi_1(Y)}\pi_1(X_2)$$
for bookkeeping purposes.  We show that $W$ can be replaced by a simple
homotopy equivalent finite $CW$ complex $W$ with disjoint subcomplexes
$V_1,V_2,\ldots,V_\mu \subset W$, such that the $F_\mu$--cover
$\wwtilde{W}$ can be constructed from a fundamental domain finite
subcomplex $U \subset \wwtilde{W}$ obtained by cutting $W$ at
$V = V_1\sqcup V_2 \sqcup \ldots \sqcup V_\mu \subset W$.

\subsection*{Algebraic transversality}

Let $A$ be an associative ring with 1. All $A$--modules will be understood
to be left $A$--modules, unless a right $A$--module structure is specified.

\fullref{algebraic transversality} develops an `algebraic
transversality' technique for cutting $A[F_\mu]$--modules along
$A$--modules, which mimics the geometric transversality method of
\fullref{combinatorial transversality}.  In \fullref{algebraic
transversality} we shall prove:

\begin{athm}
\label{thm1}
Every $A[F_\mu]$--module
chain complex $E$ admits a `Mayer--Vietoris presentation'
$$\xymatrix{0 \ar[r] & \bigoplus\limits^{\mu}_{i = 1} C^{(i)}[F_\mu]
\ar[r]^-{\di{f}} & D[F_\mu] \ar[r]& E \ar[r] &0}$$
with $C^{(i)},D$ $A$--module chain complexes, and
$f = (f^+_1z_1-f_1^- \ldots  f^+_\mu z_\mu -f_\mu^-)$
defined using $A$--module chain maps $f^+_i,f^-_i\co C^{(i)} \to D$.
If $E$ is a f.g.~free $A[F_\mu]$--module chain complex
then $C^{(i)}$, $D$ can be chosen to be f.g.~free $A$--module chain complexes,
with $D \subset E$ and $f^+_i,f^-_i\co C^{(i)} = D \cap z^{-1}_iD \to D$
given by $f^+_i(x) = x$, $f^-_i(x) = z_ix$.
\end{athm}

\begin{rem}
For $\mu = 1$ \fullref{thm1} was first proved by
Waldhausen \cite{Wald1}, being the chain complex version of the
Higman linearization trick for matrices with entries in the
Laurent polynomial extension $A[F_1] = A[z,z^{-1}]$. The algebraic
transversality theory of \cite{Wald1} applies to chain complexes
over the group rings $A[G_1*_HG_2]$ of injective amalgamated free
products $G_1*_HG_2$, using the Bass--Serre theory of groups acting
on trees. In principle, \fullref{thm1} for $\mu \geqslant 2$ could be proved by
applying \cite{Wald1} to the successive free products in
$$F_\mu = F_1*F_{\mu-1} = F_1*(F_1*F_{\mu-2}) = \cdots = 
F_1*(F_1*(F_1*\cdots*(F_1)))$$
but this would be quite awkward in practice. In view of both the
geometric motivation and the algebraic applications it is better
to prove \fullref{thm1} (as will be done in \fullref{algebraic transversality}) using the Cayley tree of
$F_\mu$ with respect to the generator set $\{z_1,z_2,\ldots,z_{\mu}\}$.
\end{rem}

\subsection*{Boundary links}
A $\mu$--component link is a (locally flat, oriented) embedding
$$\ell\co \bigsqcup\limits_\mu S^n \subset S^{n+2}.$$
Every link admits a Seifert surface $V^{n+1} \subset S^{n+2}$, a
codimension--1 submanifold with boundary
$$\partial V = \ell\Bigl(\bigsqcup\limits_\mu S^n\Bigr) \subset
    S^{n+2}.$$
By definition, $\ell$ is a \emph{$\mu$--component boundary link}
if there exists a $\mu$--component Seifert surface
$$V^{n+1} = V_1 \sqcup V_2 \sqcup \ldots \sqcup V_\mu \subset S^{n+2}.$$
The exterior of a link $\ell$ is the $(n{+}2)$--dimensional manifold
with boundary
$$(W^{n+2},\partial W) = 
\Bigl({\rm cl}\Bigl(S^{n+2}-\Bigl(\ell\Bigl(\bigsqcup\limits_\mu
S^n\Bigr)\times D^2\Bigr)\Bigr), \ell\Bigl(\bigsqcup \limits_\mu S^n\Bigr)
\times S^1\Bigr).$$
In particular, a knot $S^n \subset S^{n+2}$ is a 1--component boundary link.

The trivial $\mu$--component boundary link
$$\ell_0\co \bigsqcup\limits_\mu S^n \subset S^{n+2}$$
is defined by the connected sum of $\mu$ copies of the trivial knot
$$S^n \subset (S^n \times D^2) \cup (D^{n+1} \times S^1) = S^{n+2},$$
so that
$$\ell_0\co \bigsqcup\limits_\mu S^n \subset \mathop{\#}\limits_\mu S^{n+2} = 
    S^{n+2} = \Bigl(\bigsqcup\limits_\mu S^n \times D^2\Bigr)\cup W_0$$
has Seifert surface and exterior
$$V_0 = \bigsqcup\limits_\mu D^{n+1} ,~
W_0 = \mathop{\#}\limits_\mu (D^{n+1} \times S^1)\subset
    S^{n+2}.$$
The exterior $W_0$ has the homotopy type of
$\bigvee_\mu S^1\vee \bigvee_{\mu-1}S^{n+1}$,
with $\pi_1(W_0) = F_\mu$.

We shall make much use of the fact that the universal cover of
$BF_\mu = \bigvee_\mu S^1$ is the contractible
space with free $F_\mu$--action defined by the Cayley tree
$EF_\mu = G_\mu$ of $F_\mu$, with vertices $g \in F_\mu$ and edges
$(g,gz_i)$ $(g \in F_\mu,~1 \leqslant i \leqslant \mu)$. The
cellular chain complex $C(EF_\mu) = C(G_\mu)$ is the
standard 1--dimensional f.g.~free $\Z[F_\mu]$--module resolution of $\Z$
$$\xymatrix{
0 \ar[r] &} C_1(G_\mu) = \bigoplus\limits_{i = 1}^{\mu} \Z[F_\mu]
\xymatrix{\ar[r]^-{\displaystyle{d}}& C_0(G_\mu) = \Z[F_\mu]
\ar[r] & \Z \ar[r] &0},$$
the Mayer--Vietoris presentation with
$d = (z_1-1~z_2-1~\ldots~z_\mu-1)$.

The exterior $W$ of an $n$--dimensional link
$\ell\co \bigsqcup_\mu S^n \subset S^{n+2}$ is homotopy
equivalent to the complement $S^{n+2}\backslash
\ell\bigl(\bigsqcup_\mu S^n\bigr)$, so that
$$\begin{array}{ll}
H_*(W)& = H_*\Bigl(S^{n+2}\backslash \ell\Bigl(\bigsqcup \limits_\mu
S^n\Bigr)\Bigr)\\[2ex]
& = H^{n+2-*}\Bigl(S^{n+2},\ell\Bigl(\bigsqcup\limits_\mu S^n\Bigr)\Bigr)
 =  H^{n+1-*}\Bigl(\bigsqcup\limits_\mu S^n\Bigr)\quad(* \neq 0,n+2)
\end{array}$$
by Alexander duality. The homology groups $H_*(W),H_*(W_0)$ are thus
the same:
$$H_r(W) = H_r(W_0) = \begin{cases} \Z&\hbox{if $r = 0$} \\
\bigoplus\limits_\mu \Z&\hbox{if $r = 1$}\\
\bigoplus\limits_{\mu-1} \Z&\hbox{if $r = n+1$}\\
0&\hbox{otherwise}.
\end{cases}$$
The homotopy groups $\pi_*(W),\pi_*(W_0)$ are in general not the
same, on account of linking.  By Smythe \cite{Sm} and Gutierrez
\cite{Gu} $\ell$ is a boundary link if and only if there exists a
surjection $\pi_1(W) \to \pi_1(W_0) = F_\mu$ sending the meridians
$m_1,m_2,\ldots,m_{\mu}\co S^1 \subset W$ around the $\mu$ components
$\ell_1,\ell_2,\ldots,\ell_\mu\co S^n \subset S^{n+2}$ of $\ell$ to
the generators $z_1,z_2,\ldots,z_\mu \in F_\mu$.  We shall only be
considering boundary links $\ell$ with a particular choice of such
a surjection $\pi_1(W)\to F_\mu$, the \emph{$F_\mu$--links} of
Cappell and Shaneson \cite{CS2}.  For any such $\ell$ there exists
a map $c\co W \to W_0$ which induces a surjection $c_*\co \pi_1(W) \to
\pi_1(W_0)$ and isomorphisms $c_*\co H_*(W)\cong H_*(W_0)$.  Let
$\wwtilde{W} = c^*\wwtilde{W}_0$ be the pullback $F_\mu$--cover
of $W$, with a f.g.~free $\Z[F_\mu]$--module cellular chain
complex $C(\wwtilde{W})$.  An $F_\mu$--equivariant lift
$\wtilde{c}\co \wwtilde{W} \to \wwtilde{W}_0$ of $c$ induces a
$\Z[F_\mu]$--module chain map $\wtilde{c}\co C(\wwtilde{W}) \to
C(\wwtilde{W}_0)$ and a $\Z$--module chain equivalence $c\co C(W)
\to C(W_0)$.  A $\mu$--component Seifert surface $V = V_1\sqcup V_2
\sqcup \ldots \sqcup V_\mu \subset S^{n+2}$ for $\ell$ has a
neighbourhood $V \times [-1,1] \subset S^{n+2}$, with $V = V \times
\{0\}$. The $F_\mu$--cover $\wwtilde{W}$ can be constructed from
$F_\mu$ copies of $S^{n+2}\backslash V$, glued together using the
inclusions $f^+_i,f^-_i\co V_i \to S^{n+2}\backslash V$ defined by
$$f^{\pm}_i(v_i) = (v_i,\pm 1) \in V \times [-1,1] \subset S^{n+2}.$$
It follows that $C(\wwtilde{W})$ has a f.g.~free  $\Z[F_\mu]$--module
Mayer--Vietoris presentation
$$\xymatrix@C-5pt{0 \ar[r] &}
\bigoplus\limits^\mu_{i = 1}C(V_i)[F_\mu]
\xymatrix@C-5pt{\ar[r]^-{\di{f}}&}
C(S^{n+2}\backslash V)[F_\mu] \xymatrix@C-5pt{\ar[r] &} C(\wwtilde{W})
\xymatrix@C-5pt{\ar[r]&0}$$
with $f = f^+z-f^- = (f^+_1z_1-f_1^-~\ldots~f^+_\mu z_\mu -f_\mu^-)$.

\subsection*{Seifert and Blanchfield modules}
There are four fundamental notions in our abstract version for any
ring $A$ of the Seifert and Blanchfield modules of $\mu$--component
boundary links:
\begin{itemize}
\item[(i)] A \emph{Seifert $A$--module} is a triple
$$(P,e,\{\pi_i\}) = (~\hbox{$A$--module},~\hbox{endomorphism},~
\{\pi_i\})$$
where $\{\pi_i\co P \to P\}$ is a system of idempotents expressing $P$ as
a $\mu$--fold direct sum, with
\begin{align*}
\pi_i\co P = P_1 \oplus P_2 \oplus \cdots \oplus P_{\mu} &\to P;\\[1ex]
(x_1,x_2,\ldots,x_\mu) &\mapsto (0,\ldots,0,x_i,0,\ldots,0).
\end{align*}
Let $\Sei_{\infty}(A)$ be the category of Seifert $A$--modules. A
Seifert $A$--module $(P,e,\{\pi_i\})$ is \emph{f.g.~projective} if
$P$ is a f.g.~projective $A$--module. Let $\Sei(A) \subset
\Sei_{\infty}(A)$ be the full subcategory of the f.g.~projective
Seifert $A$--modules.
\item[(ii)] A \emph{Blanchfield $A[F_\mu]$--module} $M$ is an
$A[F_\mu]$--module such that
$$\Tor^{A[F_\mu]}_*(A,M) = 0,$$
regarding $A$ as a right $A[F_\mu]$--module via the augmentation map
$$\epsilon\co A[F_\mu] \to A;~z_i \mapsto 1.$$
Let $\Bla_{\infty}(A)$ be the category of Blanchfield $A[F_\mu]$--modules.
In \fullref{bla} Blanchfield $A[F_\mu]$--modules will be identified
with the \emph{$F_\mu$--link modules} in the sense of Sheiham \cite{Sh2},
that is $A[F_\mu]$--modules $M$ which admit an $A[F_\mu]$--module presentation
$$\xymatrix{0 \ar[r]& P[F_\mu] \ar[r]^-{\di{d}} & Q[F_\mu] \ar[r]&  M \ar[r] &    0}$$
for $A$--modules $P,Q$ with the augmentation $\epsilon(d)\co P \to Q$ an
$A$--module isomorphism.  Thus $\Bla_{\infty}(A)$ is just the
$F_\mu$--link module category $\Flk_{\infty}(A)$ of \cite{Sh2}.  A
Blanchfield $A[F_\mu]$--module $M$ has \emph{homological dimension $1$}
(or \emph{h.d.~1} for short) if it has a 1--dimensional f.g.~projective
$A[F_\mu]$--module resolution
$$\xymatrix{0 \ar[r]& K \ar[r]^-{\di{d}} & L \ar[r]&  M \ar[r] & 0}$$
with (necessarily)
$\epsilon(d) = 1\otimes d\co A\otimes_{A[F_{\mu}]}K \to A\otimes_{A[F_{\mu}]}L$
an $A$--module isomorphism.  Let $\Bla(A)
\subset \Bla_{\infty}(A)$ be the full subcategory of the h.d.~1
Blanchfield $A[F_\mu]$--modules. Let $\Flk(A) \subset \Bla(A)$ be the full
subcategory of the h.d.~1 Blanchfield modules $M$ which admit a 1--dimensional
induced f.g.~projective $A[F_\mu]$--module resolution
$$\xymatrix{0 \ar[r]& P[F_\mu] \ar[r]^-{\di{d}} & Q[F_\mu] \ar[r]&  M \ar[r] & 0}$$
with $P,Q$ f.g.~projective $A$--modules. As in \cite{Sh2}
the objects of $\Flk(A)$ will be called \emph{h.d.~1 $F_\mu$--link modules}.
\item[(iii)] The \emph{covering} of a Seifert $A$--module $(P,e,\{\pi_i\})$
is the Blanchfield $A[F_\mu]$--module
$$B(P,e,\{\pi_i\}) =  \coker(1-e+ez\co P[F_\mu] \to
P[F_\mu])$$
with $z = \sum\limits^{\mu}_{i = 1}\pi_iz_i\co P[F_{\mu}] \to
P[F_{\mu}]$, defining functors
$$B_{\infty}\co \Sei_{\infty}(A) \to
\Bla_{\infty}(A),~ B\co \Sei(A) \to \Flk(A).$$
\item[(iv)] A Seifert
$A$--module $(P,e,\{\pi_i\})$ is \emph{primitive} if
$B(P,e,\{\pi_i\}) = 0$.  Let
$$\Prim_{\infty}(A) = {\rm
ker}(B_{\infty}\co \Sei_{\infty}(A)\to \Bla_{\infty}(A))$$
be the full subcategory of $\Sei_{\infty}(A)$ with objects the primitive Seifert
$A$--modules, and let
$$\Prim(A) = \ker(B\co \Sei(A)\to \Flk(A)) \subset \Sei(A)$$
be the full subcategory of $\Sei(A)$ with objects the
primitive f.g.~projective Seifert $A$--modules.
\end{itemize}

\subsection*{Simple boundary links}

The motivational examples of f.g.~projective Seifert $\Z$--modules
and h.d.~1 $F_\mu$--link $\Z[F_\mu]$--modules come from the
$(2q{-}1)$--dimensional $\mu$--component boundary links
$\ell\co \bigsqcup_{\mu} S^{2q-1} \subset S^{2q+1}$ which are
\emph{simple}, meaning that the exterior $W$ has homotopy groups
$$\pi_r(W) = \begin{cases} F_\mu&{\rm if}~r = 1\\
0&{\rm if}~2 \leqslant r \leqslant q-1,
\end{cases}$$
so that the universal cover $\wwtilde{W}$ is $(q{-}1)$--connected.
These conditions are equivalent to the existence of a $\mu$--component Seifert
surface $V = V_1 \sqcup V_2 \sqcup \cdots \sqcup V_{\mu}$
with each component $V_i$ $(q{-}1)$--connected:
$$\pi_r(V_i) = 0\quad
  (1 \leqslant i \leqslant \mu,~1 \leqslant r \leqslant q-1).$$
The homology of the Seifert surface defines a f.g.~projective
(actually f.g.~free) Seifert $\Z$--module $(P,e,\{\pi_i\})$, with
$$\pi_i = 0 \oplus \cdots \oplus 0 \oplus 1 \oplus 0 \oplus \cdots \oplus 0\co 
P = \smash{\bigoplus\limits^\mu_{i = 1}}H_q(V_i)\to
P = \smash{\bigoplus\limits^\mu_{i = 1}}H_q(V_i)$$
and
\begin{multline*}
e = (f^+_1~f^+_2\ldots f^+_\mu)\co 
  P = H_q(V) = \bigoplus\limits^\mu_{i = 1}H_q(V_i)\longrightarrow \\
    H_q(S^{2q+1}\backslash V) = H^q(V) = H_q(V) = P
\end{multline*}
the endomorphism induced by the inclusions $f_i^+\co V_i \to
S^{2q+1}\backslash V$, identifying
$$H_q(S^{2q+1}\backslash V) = H^q(V)$$
by Alexander duality and $H^q(V) = H_q(V)$ by Poincar\'e
duality. The covering of $(P,e,\{\pi_i\})$ is the h.d.~1 $F_\mu$--link
$\Z[F_\mu]$--module
$$B(P,e,\{\pi_i\}) = H_q(\wwtilde{W})$$
defined by the homology of the $F_\mu$--cover $\wwtilde{W}$ of the
exterior $W$.  The f.g.~projective Seifert $\Z$--module
$(P,e,\{\pi_i\})$ is primitive if and only if
$H_q(\wwtilde{W}) = 0$; for $q \geqslant 2$ this is the case if
and only if $\ell$ is unlinked (Gutierrez \cite{Gu}).

\subsection*{Blanchfield  =  Seifert/primitive}

\fullref{modules} uses algebraic transversality to prove that
every h.d.~1 $F_\mu$--link module $M$ is isomorphic to the covering
$B(P,e,\{\pi_i\})$ of a f.g.~projective Seifert $A$--module
$(P,e,\{\pi_i\})$, and that morphisms of h.d.~1 $F_\mu$--link
modules can be expressed as fractions of morphisms of f.g.~projective
Seifert $A$--modules.

The algebraic relation between Seifert $A$--modules and Blanchfield
$A[F_\mu]$--modules for $A = \Z$ was first investigated
systematically in the knot case $\mu = 1$, by Levine \cite{L1,L2} and
Trotter \cite{T}, and for the link case $\mu \geqslant 1$ by
Farber \cite{Fa1,Fa3} and Sheiham \cite{Sh2}. In particular,
\cite{Sh2} expressed the Blanchfield module
category $\Bla_\infty(A) = \Flk_\infty(A)$ as the quotient of the Seifert
$A$--module category $\Sei_{\infty}(A)$ by the primitive Seifert
$A$--module subcategory $\Prim_{\infty}(A)$, as we now recall.

Let $\A$ be an abelian category.  By definition, a \emph{Serre
subcategory} $\C \subset \A$ is a non-empty full
subcategory such that for every exact sequence  in $\A$
$$0 \to M' \to M \to M'' \to 0$$
$M$ is an object in $\C$ if and
only if $M',M''$ are objects in $\C$.  Gabriel \cite{Ga}
defined the quotient abelian category $\A/\C$ with
the same objects as $\A$ but different morphisms: if
$M,N$ are objects in $\A$ then
$$\Hom_{\A/\C}(M,N) = 
\varinjlim \Hom_\A(M',N'')$$
with the direct limit taken over all the exact sequences in $\A$
$$0 \to M' \to M \to M'' \to 0,~0 \to N' \to N \to N'' \to 0$$
with $M'',N'$ objects in $\C$.  The canonical functor
$F\co \A \to \A/\C;A \mapsto A$ sends each object $C$ in $\C$ to $F(C) = 0$,
and has the universal property that for any exact functor
$G\co \A \to \B$ such that $G(C) = 0$ for all objects in $\C$
there exists a unique functor $\wwbar{G}\co \A/\C \to \B$ such that
$\wwbar{G}F = G$. In particular, if $\B$ is an exact category and
$G\co \A \to \B$ is an exact functor then the full subcategory $\C \subset \A$
with objects $C$ such that $G(C) = 0$ is a Serre subcategory, and there
is induced a functor $\wwbar{G}\co \A/\C \to \B;A \mapsto G(A)$
such that $G = \wwbar{G}F$.

By definition, a category is \emph{small} if the class of morphisms is a set.
In order to avoid set-theoretic difficulties we shall only be
dealing with categories which are \emph{essentially small}, ie  equivalent
to a small category.

Let $\A$ be an essentially small category, and let $\Sigma$ be a set
of morphisms in $\mathcal A$, e.g. the morphisms of a subcategory.
A \emph{category of fractions} $\Sigma^{-1}\A$
is a category with a universally $\Sigma$--inverting functor
$F\co \A \to \Sigma^{-1}\A$, meaning that:
\begin{itemize}
\item[(i)] $F$ sends each $f \in \Sigma$ to an isomorphism $F(f)$
in $\Sigma^{-1}\A$,
\item[(ii)] for any functor $G\co \A \to \B$
which sends each $f \in \Sigma$ to an isomorphism $G(f)$
there exists a unique functor $\wwbar{G}\co \Sigma^{-1}\A \to
\B$ such that $\wwbar{G}F = G$.
\end{itemize}
An essentially small category of fractions $\Sigma^{-1}\A$ exists,
with the same objects as $\A$, and such a category
is unique up to isomorphism (Gabriel and Zisman \cite{GZ},
Borceux \cite[5.2.2]{Bo}).
For example, if $\mathcal A$ is an abelian category and $\C
\subset \A$ is a Serre subcategory, then
$$\A/\C = \Sigma^{-1}\A$$
is a category of fractions inverting the set $\Sigma$ of morphisms
$f$ in $\A$ with $\ker(f)$ and $\coker(f)$ in $\C$.

An $A[F_\mu]$--module $M$ is Blanchfield if and only if the
$A$--module morphism
$$\gamma\co \bigoplus\limits_\mu M \to M;~
(m_1,m_2,\ldots,m_\mu) \mapsto \sum\limits^\mu_{i = 1}(z_i-1)m_i$$
is an isomorphism, called the \emph{Sato isomorphism} (after \cite{Sa},
the case $A = \Z$).
As in Sheiham \cite{Sh2}, for any Blanchfield
$A[F_\mu]$--module $M$ use the $A$--module morphisms
\begin{align*}
p_i\co \bigoplus\limits_\mu M &\to M;
  \quad(m_1,m_2,\ldots,m_\mu) \mapsto m_i,\\
\omega\co \bigoplus\limits_\mu M &\to M;
  \quad(m_1,m_2,\ldots,m_\mu) \mapsto \smash{\sum\limits^\mu_{i = 1}}m_i,\\
\pi_i = \gamma p_i \gamma^{-1}\co M &\to M,\\
 e = \omega \gamma^{-1}\co M &\to M
\end{align*}
to define a Seifert $A$--module $U(M) = (M,e,\{\pi_i\})$.

The categories $\Prim_{\infty}(A)$, $\Sei_{\infty}(A)$ are abelian,
while $\Bla_{\infty}(A)$ is in general only exact.  The covering
functor
$B_{\infty}\co \Sei_{\infty}(A) \to \Bla_{\infty}(A)$ was shown in
\cite[5.2]{Sh2} to be exact, so that $\Prim_{\infty}(A) \subset
\Sei_{\infty}(A)$ is a Serre subcategory.  Thus if $\Xi_{\infty}$ is
the set of morphisms $f$ in $\Sei_{\infty}(A)$ such that $B(f)$ is an
isomorphism in $\Bla_{\infty}(A)$, or equivalently $\ker(f)$ and
$\coker(f)$ are in $\Prim_{\infty}(A)$, then
$$\Sei_{\infty}(A)/\Prim_{\infty}(A) = \Xi_{\infty}^{-1}\Sei_{\infty}(A).$$
The induced exact functor
$\wbar{B}_{\infty}\co \Sei_{\infty}(A)/\Prim_{\infty}(A) \to \Bla_{\infty}(A)$
is such that
$$B_{\infty}\co \Sei_{\infty}(A) \to \Sei_{\infty}(A)/\Prim_{\infty}(A)
\xymatrix{\ar[r]^-{\di{\wbar{B}_\infty}}&} \Bla_{\infty}(A)$$
and has the universal property of inverting $\Xi_{\infty}$. The
functor $\wbar{B}_{\infty}$ was shown to be an equivalence in
\cite[5.15]{Sh2} using the fact that the functor
$$U_{\infty}\co \Bla_{\infty}(A) \to \Sei_{\infty}(A);~M \mapsto U(M)$$
is right adjoint to $B$:
for any Seifert $A$--module $V$ there is a natural isomorphism
$$\Hom_{\Bla_{\infty}(A)}(B(V),M)\cong~\Hom_{\Sei_{\infty}(A)}(V,U(M)).$$
The functor $U_{\infty}$ is fully faithful, allowing $\Bla_{\infty}(A)$
to be regarded as a full subcategory of $\Sei_{\infty}(A)$.  By
\cite[5.15]{Sh2} $U_{\infty}$ induces a functor
$$\wwbar{U}_{\infty}\co \Bla_{\infty}(A)\to \Sei_{\infty}(A)/\Prim_{\infty}(A)$$
which is an equivalence inverse to $\wbar{B}_{\infty}$. Thus up
to equivalence
$$\Sei_{\infty}(A)/\Prim_{\infty}(A) = \Xi_{\infty}^{-1}\Sei_{\infty}(A) = 
\Bla_{\infty}(A).$$
\indent The categories $\Prim(A),\Sei(A),\Flk(A),\Bla(A)$ are exact but
not in general abelian.  As in \cite{Sh2} let
$\Sei(A)/\Prim_{\infty}(A) \subset \Sei_{\infty}(A)/\Prim_{\infty}(A)$
be the full subcategory with objects in $\Sei(A)$.  The equivalence
$$\wbar{B}_{\infty}\co \Sei_{\infty}(A)/\Prim_{\infty}(A)
\xymatrix{\ar[r]^-{\approx}&}\Bla_{\infty}(A)$$
was shown in \cite[5.17]{Sh2} to restrict to an
equivalence of exact sequences
$$\wbar{B}\co \Sei(A)/\Prim_{\infty}(A) \xymatrix{\ar[r]^-{\approx}&}
\Flk(A)$$
with
$$B\co \Sei(A) \to  \Sei(A)/\Prim_{\infty}(A)
\xymatrix{\ar[r]^{\di{\wbar{B}}}_{\approx}&} \Flk(A).$$
From the construction of $\Sei_{\infty}(A)/\Prim_{\infty}(A)$
a morphism in $\Sei(A)/\Prim_{\infty}(A)$ may involve objects
in $\Sei_{\infty}(A)$ which are not in $\Sei(A)$, so that the
equivalence $\wbar{B}$
cannot be used to relate the algebraic $K$--theories of
$\Sei(A)$ and $\Flk(A)$.

A category of fractions $\Sigma^{-1}\A$ has a
\emph{left calculus of fractions} if:
\begin{itemize}
\item[(i)] $(1\co A \to A) \in \Sigma$ for every object $A$ in $\mathcal A$,
\item[(ii)] if $(s\co A \to B),(t\co B \to C) \in \Sigma$ then $(ts\co A \to C)
\in \Sigma$,
\item[(iii)] for any
$f\co A \to B$ in $\A$ and $s\co A \to D$ in $\Sigma$
there exist
$g\co D \to C$ in $\A$ and $t\co B \to C$ in $\Sigma$
such that $tf = gs\co A \to C$,
\item[(iv)] for any  $f,g\co A \to B$ in $\A$ and
$s\co D \to A$ in $\Sigma$ with $fs = gs\co D \to B$ there
exists $(t\co B \to C) \in \Sigma$ with $tf = tg\co A \to C$.
\end{itemize}
It then follows that a morphism $A \to B$ in $\Sigma^{-1}\A$
can be regarded as an equivalence class $s^{-1}f$ of pairs
$(f\co A \to C,s\co B \to C)$ of morphisms in $\mathcal A$
with $s \in \Sigma$, where
$$\begin{array}{ll}
(f,s) \sim (f',s')&\hbox{\rm if there exist morphisms}~
g\co C \to D,~g'\co C' \to D~{\rm in}~\A\\[1ex]
&{\rm with}~(gs = g's'\co B \to D) \in \Sigma~{\rm and}~gf = g'f'\co A\to D
\end{array}$$
so that
$$s^{-1}f = (gs)^{-1}(gf) = (g's')^{-1}(g'f') = {s'}^{-1}f'\co 
A \to B~{\rm in}~\Sigma^{-1}\A.$$
\indent
Let $\Xi$ be the set of morphisms $f$ in $\Sei(A)$ such that $B(f)$ is
an isomorphism in $\Flk(A)$, or equivalently such that $\ker(f)$
and $\coker(f)$ are in $\Prim_{\infty}(A)$.
In \fullref{modules} we shall prove:

\begin{athm}\label{thm2}
{\rm (i)}\qua The category of fractions
$\Xi^{-1}\Sei(A)$ has a left calculus of fractions, and the
covering functor $B\co \Sei(A) \to \Flk(A)$ induces an equivalence of exact
categories
$$\wbar{B}\co \Xi^{-1}\Sei(A) \xymatrix{\ar[r]^-{\approx}&} \Flk(A).$$
{\rm (ii)}\qua The h.d.~1 Blanchfield $A[F_\mu]$--module category
$\Bla(A)$ is the idempotent completion of the h.d.~1 $F_\mu$--link
module category $\Flk(A)$.
\end{athm}

The key step in the proof of \fullref{thm2} (i) is the use of the algebraic
transversality \fullref{thm1} to verify that for any h.d.~1 $F_\mu$--link
module $M$ the Seifert $A$--module $U(M)$ is a direct limit of morphisms
in $\Xi$.

\subsection*{Primitive  =  near-projection}

\fullref{kernel} gives an intrinsic characterization of the
primitive f.g.~projective Seifert $A$--modules $(P,e,\{\pi_i\})$ as
generalized near-projections.

An endomorphism $e\co P \to P$ of an $A$--module $P$ is {\it
nilpotent} if $e^N = 0$ for some $N \geqslant 0$.

An endomorphism $e\co P \to P$ is a \emph{near-projection}
if $e(1-e)\co P \to P$ is nilpotent (L\"uck and Ranicki \cite{LR}).

In \fullref{kernel} we shall prove:

\begin{athm}\label{thm3}
A f.g.~projective Seifert $A$--module
$(P,e,\{\pi_i\})$ is primitive if and only if it can be expressed
as
$$(P,e,\{\pi_i\}) = \biggl(P^+\oplus P^-,\begin{pmatrix}
e^{++} & e^{+-} \\
e^{-+} & e^{--}
\end{pmatrix},\{\pi_i^+\}\oplus \{\pi_i^-\}\biggr)$$
and the $2\mu$--component Seifert $A$--module
$$(P',e',\pi') = \left( P^+\oplus P^-\ ,\
\begin{pmatrix}
e^{++} & -e^{+-} \\
e^{-+} & 1-e^{--}
\end{pmatrix}\ ,\ \{\pi_i^+\}\oplus\{\pi_i^-\}\right)$$
is such that $e'z'\co P'[F_{2\mu}] \to P'[F_{2\mu}]$ is
nilpotent, with $F_{2\mu}$ the free group on $2\mu$ generators
$z'_1,\ldots,z'_{2\mu}$.
\end{athm}

For $\mu = 1$ the condition for a f.g.~projective Seifert $A$--module
$(P,e,\{\pi_i\})$ to be primitive is just that $e$ be a near-projection.
For $\mu = 1$ \fullref{thm3} is just the result of Bass, Heller and Swan \cite{BHS}
that $1-e+ez\co P[z,z^{-1}] \to P[z,z^{-1}]$ is an $A[z,z^{-1}]$--module
isomorphism if and only if $e$ is a near-projection, if and only if
$(P,e) = (P^+,e^{++})\oplus (P^-,e^{--})$ with
$e^{++}\co P^+ \to P^+$ and $1-e^{--}\co P^- \to P^-$ nilpotent.

\subsection*{Algebraic $K$--theory}

\fullref{ktheory} obtains results on the algebraic $K$--theory of
$A[F_\mu]$, $\Prim(A)$, $\Sei(A)$, $\Flk(A)$ and $\Bla(A)$, using
the algebraic $K$--theory noncommutative localization exact sequences
of Schofield \cite{Sc} and Neeman--Ranicki \cite{NR1,NR2}.

The class group $K_0(\E)$ of an exact category $\E$ is the Grothendieck
group with one generator $[M]$ for each object $M$ in $\E$, and one relation
$[K]-[L]+[M] = 0$ for each exact sequence in $\E$
$$0 \to K \to L \to M \to 0.$$
The algebraic $K$--groups $K_n(\E)$ are defined by Quillen \cite{Q} for
$n \geqslant 1$ and by Schlichting \cite{Schl} for $n \leqslant -1$.
Write
\begin{align*}
\PRIM_*(A) & =  K_*(\Prim(A)), &\SEI_*(A) & =  K_*(\Sei(A)),\\
\BLA_*(A) & =  K_*(\Bla(A)), &\FLK_*(A) & =  K_*(\Flk(A)),
\end{align*}
noting that $\BLA_n(A) = \FLK_n(A)$ for $n \neq 0$.

\begin{athm}\label{thm4}
{\rm (i)}\qua The algebraic $K$--groups of $A[F_\mu]$ split as
$$K_*(A[F_\mu]) = K_*(A) \oplus \bigoplus\limits_{\mu} K_{*-1}(A)
\oplus \widetilde{\PRIM}_{*-1}(A).$$
{\rm (ii)}\qua The sequence of functors
$$\xymatrix{\Prim(A) \ar[r]& \Sei(A) \ar[r]^-{\di{B}} & \Bla(A)}$$
induces a long exact sequence of algebraic $K$--groups
$$\cdots \to \PRIM_n(A) \to \SEI_n(A) \xymatrix{\ar[r]^-{\displaystyle{B}}&}
\BLA_n(A) \to \PRIM_{n-1}(A) \to \cdots$$
with
$${\rm im}(B\co \SEI_0(A) \to \BLA_0(A)) = \FLK_0(A) \subseteq \BLA_0(A).$$
{\rm (iii)}\qua The exact sequence in {\rm (ii)} splits as a direct sum of
exact sequences
$$\begin{array}{l}
\cdots \to \bigoplus \limits_{2\mu} K_n(A) \to
\bigoplus \limits_\mu K_n(A) \xymatrix@C-10pt{\ar[r]^-{\di{0}}&}
\bigoplus \limits_\mu K_{n-1}(A) \to
\bigoplus \limits_{2\mu} K_{n-1}(A) \to \cdots,\\[2ex]
\cdots \to \widetilde{\PRIM}_n(A) \to \widetilde{\SEI}_n(A) \to
\widetilde{\BLA}_n(A) \to \widetilde{\PRIM}_{n-1}(A) \to \cdots.
\end{array}$$
\end{athm}

For $\mu = 1$ $\Prim(A)$ is the exact category of f.g.~projective
$A$--modules $P$ with a near-projection $e\co P \to P$, which is equivalent
to the product $\Nil(A) \times \Nil(A)$ of two copies of
the exact category $\Nil(A)$ of f.g.~projective $A$--modules $P$
with a nilpotent endomorphism $e\co P \to P$, and
\begin{eqnarray*}
\PRIM_*(A) & = & K_*(\Prim(A)) = \NIL_*(A) \oplus \NIL_*(A),\\[1ex]
\NIL_*(A) & = & K_*(\Nil(A)) = K_*(A) \oplus \widetilde{\NIL}_*(A),\\[1ex]
\widetilde{\PRIM}_*(A) & = & \widetilde{\NIL}_*(A) \oplus
\widetilde{\NIL}_*(A).
\end{eqnarray*}
Thus for $\mu = 1$ \fullref{thm4} (i) is just the splitting theorem
of Bass, Heller and Swan \cite{BHS}, \cite{B2} for $K_1(A[z,z^{-1}])$
and its generalization to the higher $K$--groups
$$K_*(A[z,z^{-1}]) = K_*(A) \oplus K_{*-1}(A) \oplus
\widetilde{\NIL}_{*-1}(A) \oplus \widetilde{\NIL}_{*-1}(A).$$
\fullref{thm4} (ii)--(iii) is new even in the case $\mu = 1$.

Let $\Sigma^{-1}A[F_\mu]$ be the noncommutative Cohn (ie  universal)
localization of $A[F_\mu]$ inverting the set $\Sigma$ of the morphisms of
f.g.~projective $A[F_\mu]$--modules which induce isomorphisms of
f.g.~projective $A$--modules under the augmentation $\epsilon:
A[F_\mu] \to A$. The exact category $H(A[F_\mu],\Sigma)$
of h.d.~1 $\Sigma$--torsion $A[F_\mu]$--modules is such that
$$H(A[F_\mu],\Sigma) = \Bla(A),~K_*(H(A[F_\mu],\Sigma)) = \BLA_*(A).$$

\begin{athm}\label{thm5}
{\rm (i)}\qua The localization exact sequence
$$K_1(A[F_\mu]) \to K_1(\Sigma^{-1}A[F_\mu])
\to K_0(H(A[F_\mu],\Sigma)) \to K_0(A[F_\mu]) \to \cdots $$
splits as a direct sum of the exact sequences
$$\begin{array}{l}
K_1(A)\oplus\bigoplus \limits_\mu K_0(A) \to K_1(A)
\xymatrix@C-10pt{\ar[r]^-{\di{0}}&} \bigoplus \limits_\mu K_{-1}(A) \to
K_0(A) \oplus\bigoplus \limits_\mu K_{-1}(A) \to \cdots,\\[1ex]
\widetilde{\PRIM}_0(A) \to \widetilde{\SEI}_0(A)
\to \widetilde{\BLA}_0(A) \to \widetilde{\PRIM}_{-1}(A) \to \cdots.
\end{array}$$
{\rm (ii)}\qua If $\Sigma^{-1}A[F_\mu]$ is stably flat (ie if
$\Tor^{A[F_\mu]}_*(\Sigma^{-1}A[F_\mu],\Sigma^{-1}A[F_\mu]) = 0$
for $* \geqslant 1$)
the exact sequences and the splitting in {\rm (i)} extend to the left,
involving the algebraic $K$--groups $K_n$ for $n \geqslant 2$, with
$$K_*(\Sigma^{-1}A[F_\mu]) = K_*(A) \oplus \widetilde{\SEI}_{*-1}(A).$$
\end{athm}

For $\mu = 1$ $\Sei(A)$ is the exact category $\Endd(A)$ of f.g.~projective
$A$--modules $P$ with an endomorphism $e\co P \to P$, and
$$\SEI_*(A) = K_*(\Endd(A)) = \ENDD_*(A) = K_*(A) \oplus
\widetilde{\ENDD}_*(A).$$
The special case of \fullref{thm5} (i)
$$K_1(\Sigma^{-1}A[z,z^{-1}]) = K_1(A) \oplus \widetilde{\ENDD}_0(A)$$
is the splitting theorem of Ranicki \cite[10.21]{RHK}.

We are grateful to Pere Ara, Warren Dicks, Marco Schlichting and the
referee for helpful comments on the preprint version of the paper,
which have led to various improvements.  In particular, it was Pere Ara
who pointed out that the Blanchfield $A[F_\mu]$--module category
$\Bla_{\infty}(A)$ is the same as the $F_\mu$--link module category
$\Flk_{\infty}(A)$ of \cite{Sh2}.

\section{Combinatorial transversality for $F_\mu$--covers}
\label{combinatorial transversality}

For $\mu \geqslant 1$ let
$F_{\mu} = \langle z_1,z_2,\ldots,z_{\mu} \rangle$
be the free group with generators $z_1,z_2,\ldots,z_{\mu}$.

\subsection{$F_\mu$--covers}

\begin{definition} {\rm An \emph{$F_{\mu}$--cover} of a space $W$
is a regular covering $p\co \wwtilde{W}\to W$ with group of
covering translations $F_{\mu}$.}
\end{definition}

\indent A classifying space $BF_\mu$ for $F_\mu$--covers is a
connected space such that
$$\pi_j(BF_{\mu}) = \begin{cases}
F_{\mu}&{\rm if}~j = 1\\
0&{\rm if}~j \geqslant 2.
\end{cases}$$
The universal cover of $BF_\mu$ is an $F_\mu$--cover
$$p_{\mu}\co EF_{\mu} = \widetilde{BF}_{\mu}\to BF_{\mu}$$
with $EF_{\mu}$ a contractible space with a free $F_{\mu}$--action.

\begin{proposition}
{\rm (i)}\qua Given an $F_{\mu}$--cover $p\co \wwtilde{W} \to W$ and a
map $f\co V \to W$ there is defined a pullback square
$$\xymatrix@C+10pt@R+10pt{\wtilde{V} \ar[r]^-{\di\wtilde{f}}
\ar[d]_-{\di{f^*p}} &
\wwtilde{W} \ar[d]^-{\di{p}} \\
V \ar[r]^-{\di{f}} & W}$$ with
$$\begin{array}{l}
\wtilde{V} = f^*\wwtilde{W} = \bigl\{(x,y) \in V\times
\wwtilde{W}\,\vert\,
f(x) = p(y) \in W\bigr\},\\[1ex]
f^*p\co \wtilde{V} \to V;~(x,y) \mapsto x,\quad
\wtilde{f}\co \wtilde{V} \to \wwtilde{W};~(x,y) \mapsto y
\end{array}$$
such that $f^*p\co \wtilde{V} \to V$ is the pullback $F_{\mu}$--cover.

{\rm (ii)}\qua The $F_{\mu}$--covers $p\co \wwtilde{W} \to W$ of a space
$W$ are classified by the homotopy classes of maps $c\co W \to BF_{\mu}$
with
$$\begin{array}{l}
\wwtilde{W} = c^*EF_{\mu} = \bigl\{(x,y) \in W \times EF_{\mu}\,\vert\,
    c(x) = [y] \in BF_{\mu}\bigr\},\\[1ex]
    p(x,y) = c^*p_{\mu}(x,y) = x.
\end{array}$$
For a connected space $W$ the homotopy classes of maps $c\co W \to
BF_{\mu}$ are in one-one correspondence with the morphisms
$c_*\co \pi_1(W) \to F_{\mu}$; the connected $F_{\mu}$--covers
$\wwtilde{W}$ correspond to surjections $c_*\co \pi_1(W) \to F_{\mu}$.
\end{proposition}
\begin{proof} Standard.
\end{proof}

\subsection{The Cayley tree $G_\mu$}

We shall be working with the following explicit constructions of
$BF_{\mu}$ and $EF_{\mu}$, as well as the Cayley tree of $F_\mu$:

\begin{definition} \label{cayley}
{\rm The \emph{Cayley tree} $G_{\mu}$ is the tree with vertex set
$$G_{\mu}^{(0)} = F_{\mu}$$
and edge set
$$G_{\mu}^{(1)} = \{(g,gz_i)\,\vert\,g \in
F_{\mu},\,1 \leqslant i \leqslant \mu\} \subset G_{\mu}^{(0)} \times
G_{\mu}^{(0)}.$$
$$\Draw
\MoveTo(0,0) \Text(--$\bullet$--) \MoveTo(0,60)
\Text(--$\bullet$--) \MoveTo(0,-60) \Text(--$\bullet$--)
\MoveTo(60,0) \Text(--$\bullet$--) \MoveTo(-60,0)
\Text(--$\bullet$--) \LineAt(-90,0,90,0) \LineAt(0,-90,0,90)
\LineAt(-60,40,-60,-40) \LineAt(-40,-60,40,-60)
\LineAt(60,-40,60,40) \LineAt(-40,60,40,60) \MoveTo(-70,10)
\Text(--$z_i^{-1}$--) \MoveTo(-25,10) \Text(--$(z_i^{-1},1)$--)
\MoveTo(7,10) \Text(--$1$--) \MoveTo(35,10) \Text(--$(1,z_i)$--)
\MoveTo(70,7) \Text(--$z_i$--) \MoveTo(10,-70) \Text(--$z_j$--)
\MoveTo(16,-30) \Text(--$(1,z_j)$--) \MoveTo(20,30)
\Text(--$(z_j^{-1},1)$--) \MoveTo(12,70) \Text(--$z_j^{-1}$--)
\EndDraw$$
Define a transitive $F_{\mu}$--action on $G_\mu$
$$F_\mu \times G_\mu \to G_\mu ;~(g,x) \mapsto gx$$
with quotient the one-point union of $\mu$ circles
$$G_\mu/F_\mu = BF_\mu = \smash{\bigvee_{\mu}} S^1.$$}
\end{definition}

Let
$$I_{\mu} = 
\bigcup^{\mu}_{i = 1} [e^-_i,e^+_i] \subset \RR^{\mu}$$
with
\begin{align*}
e^+_i & =  (0,\ldots,0,1,0,\ldots,0),~e^-_i = (0,\ldots,0,-1,0,\ldots,0)
\in \RR^{\mu},\\
[e^-_i,e^+_i] & =  \{(0,\ldots,0,t,0,\ldots,0)\,\vert\, -1 \leqslant t \leqslant 1\}
\subset \RR^{\mu}.
\end{align*}
Thus $I_\mu$ is the one-point union of $\mu$ copies of the interval
$[-1,1] \subset \RR$, identifying the $\mu$ copies of $0 \in [-1,1]$.

We regard $BF_{\mu}$ as the quotient space of $I_{\mu}$
$$BF_{\mu} = I_{\mu}/\{e^+_i \sim e^-_i\,\vert\,1 \leqslant i \leqslant \mu\} = 
    \bigvee_{\mu}S^1,$$
the one-point union of $\mu$ copies of the circle $S^1 = [-1,1]/(-1
\sim 1)$ in which the $\mu$ copies of $[0] \in S^1$ are
identified, with
$$e_i = [e^+_i] = [e^-_i]\neq [0] \in BF_{\mu}$$
a point in the $i^{th}$ circle.  The universal cover $EF_{\mu}$ of
$BF_{\mu}$ is
$$EF_{\mu} = (F_{\mu} \times I_{\mu})/\{(g,e^+_i) \sim (gz_i,e^-_i)\,\vert\,g
    \in F_{\mu},1 \leqslant i \leqslant \mu\},$$
a contractible space with a free $F_{\mu}$--action
$$F_{\mu} \times EF_{\mu} \to EF_{\mu};~(g,(h,x)) \mapsto (gh,x)$$
and covering projection
$$p_{\mu}\co EF_{\mu} \to BF_{\mu};~[g,x] \mapsto [x].$$
Define an $F_{\mu}$--equivariant homeomorphism
    $\xymatrix{G_{\mu} \ar[r]^-{\cong} & EF_{\mu}}$
by sending the vertex $g \in \smash{G^{(0)}_{\mu}} = F_{\mu}$ to
the point $(g,0) \in EF_{\mu}$, and the edge $(g,gz_i) \in
\smash{G^{(1)}_{\mu}}$ to the line segment
$$\{(g,te^+_i)\,\vert\, 0 \leqslant t \leqslant 1\} \cup
\{(gz_i,te^-_i)\,\vert\, 0 \leqslant t \leqslant 1\} \subset
    EF_{\mu}$$
with endpoints $(g,0)$, $(gz_i,0) \in EF_{\mu}$. The projection
$G_\mu \to G_\mu/F_\mu$ can thus be identified with the universal
cover $p_\mu\co EF_\mu \to BF_\mu$.

\subsection{Fundamental domains}

\begin{definition}{\rm
A \emph{fundamental domain} of an $F_{\mu}$--cover
$p\co \wwtilde{W}\to W$ is a closed subspace $U \subset
\wwtilde{W}$ such that
\begin{itemize}
\item[(a)] $F_{\mu}U = \wwtilde{W}$, or equivalently $p(U) = W$,
\item[(b)] for any $g,h \in F_{\mu}$
$$gU \cap hU = \begin{cases}
gV_i&{\rm if}~g = hz_i\\
hV_i&{\rm if}~g = hz^{-1}_i\\
gU&{\rm if}~g = h\\
\emptyset&{\rm otherwise}
\end{cases}$$
with $V_i = U \cap z_i^{-1}U$.
\end{itemize}}
\end{definition}

{\tiny
%%%Start of picture.
\tracingstats = 2
\newcounter{rot}
\newcounter{rot2}
\newcounter{posx}
\newcounter{posy}
%
%The counter psz controls the overall size of the picture.
\newcounter{psz}
%The counters pcentx and pcenty give `origin' coordinates relative to which
%a picture is to be drawn.
\newcounter{pcentx}
\newcounter{pcenty}
\newcounter{subpcentx}
\newcounter{subpcenty}
\newcounter{subsubpcentx}
\newcounter{subsubpcenty}
%

%
%
%The counters ovalthk and ovalthn determine the dimensions of
%ovals in the `thick' and `thin' directions.
\newcounter{ovalthk}
\newcounter{ovalthn}
\newcounter{subovalthk}
\newcounter{subovalthn}
\newcounter{subsubovalthk}
\newcounter{subsubovalthn}
\setcounter{psz}{50}
\setcounter{subsubpcentx}{\value{psz}*19/6}
\setcounter{subsubpcenty}{0}
\setcounter{subpcentx}{\value{psz}*21/10}
\setcounter{subpcenty}{0}
\setcounter{pcentx}{0}
\setcounter{pcenty}{0}
\setcounter{subsubovalthk}{\value{psz}/8}
\setcounter{subsubovalthn}{\value{psz}/10}
\setcounter{subovalthk}{\value{psz}/6}
\setcounter{subovalthn}{\value{psz}/8}
\setcounter{ovalthk}{\value{psz}/4}
\setcounter{ovalthn}{\value{psz}/5}
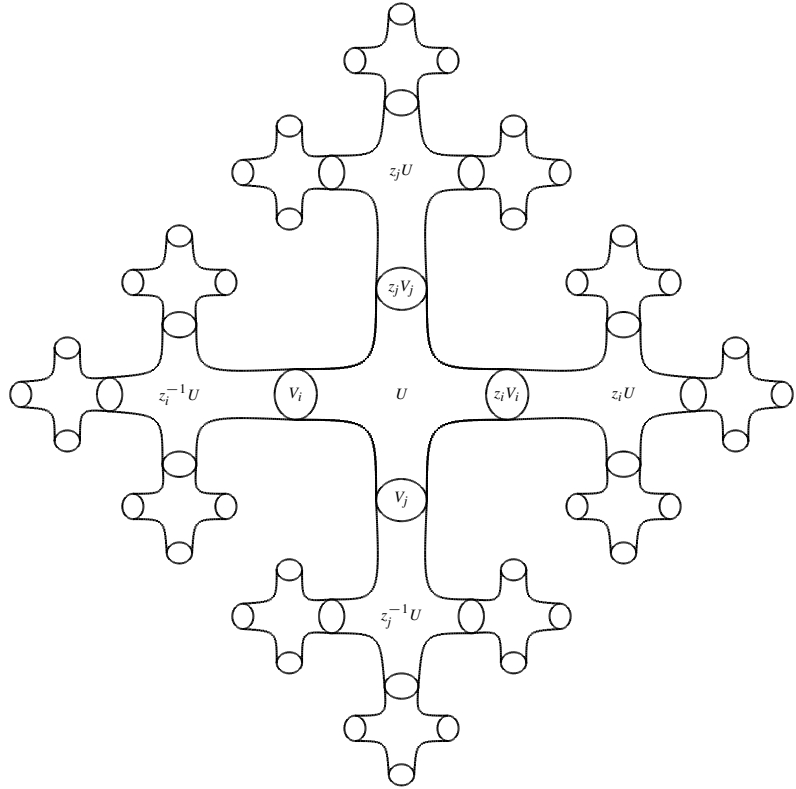
\begin{equation*}\label{picture}
\Draw(0.8pt,0.8pt)
%\Scale(0.5,0.5)
%
%
%\MoveTo(\value{pcentx},\value{pcenty})
%\Text(--$y_1^2P$--)
%
%\setcounter{posx}{\value{pcentx}-\value{psz}*2/3}
%\setcounter{posy}{\value{pcenty}+\value{subsubovalthk}}
%\MoveTo(\value{posx},\value{posy})
%\MarkLoc(b-1+)
%
%
%
%\TranslateAxes
%
%%%%%%%%%%%%%%%
%Now for some really small domains
%%%%%%%%%%%%%%%
\def\smallestpic{%
\setcounter{posx}{\value{subsubpcentx}-\value{subsubovalthk}}
\setcounter{posy}{\value{subsubpcenty}+\value{psz}*9/20}
\MoveTo(\value{posx},\value{posy})
\MarkLoc(c+2-)
\setcounter{posx}{\value{subsubpcentx}-\value{subsubovalthk}}
\setcounter{posy}{\value{subsubpcenty}+\value{subsubovalthk}*5/4}
\MoveTo(\value{posx},\value{posy})
\MarkLoc(C-+)
\setcounter{posx}{\value{subpcentx}+\value{psz}*2/3}
\setcounter{posy}{\value{subpcenty}+\value{subovalthk}}
\MoveTo(\value{posx},\value{posy})
\MarkLoc(b+1+)
\Curve(b+1+,C-+,C-+,c+2-)
\setcounter{posx}{\value{subsubpcentx}-\value{subsubovalthk}}
\setcounter{posy}{\value{subsubpcenty}-\value{psz}*9/20}
\MoveTo(\value{posx},\value{posy})
\MarkLoc(c-2-)
\setcounter{posx}{\value{subsubpcentx}-\value{subsubovalthk}}
\setcounter{posy}{\value{subsubpcenty}-\value{subsubovalthk}*5/4}
\MoveTo(\value{posx},\value{posy})
\MarkLoc(C--)
\setcounter{posx}{\value{subpcentx}+\value{psz}*2/3}
\setcounter{posy}{\value{subpcenty}-\value{subovalthk}}
\MoveTo(\value{posx},\value{posy})
\MarkLoc(b+1-)
\Curve(b+1-,C--,C--,c-2-)
\setcounter{posx}{\value{subsubpcentx}+\value{subsubovalthk}}
\setcounter{posy}{\value{subsubpcenty}+\value{psz}*9/20}
\MoveTo(\value{posx},\value{posy})
\MarkLoc(c+2+)
\setcounter{posx}{\value{subsubpcentx}+\value{psz}*9/20}
\setcounter{posy}{\value{subsubpcenty}+\value{subsubovalthk}}
\MoveTo(\value{posx},\value{posy})
\MarkLoc(c+1+)
\setcounter{posx}{\value{subsubpcentx}+\value{subsubovalthk}}
\setcounter{posy}{\value{subsubpcenty}+\value{subsubovalthk}*5/4}
\MoveTo(\value{posx},\value{posy})
\MarkLoc(C++)
\Curve(c+2+,C++,C++,c+1+)
\setcounter{posx}{\value{subsubpcentx}+\value{subsubovalthk}}
\setcounter{posy}{\value{subsubpcenty}-\value{psz}*9/20}
\MoveTo(\value{posx},\value{posy})
\MarkLoc(c-2+)
\setcounter{posx}{\value{subsubpcentx}+\value{psz}*9/20}
\setcounter{posy}{\value{subsubpcenty}-\value{subsubovalthk}}
\MoveTo(\value{posx},\value{posy})
\MarkLoc(c+1-)
\setcounter{posx}{\value{subsubpcentx}+\value{subsubovalthk}}
\setcounter{posy}{\value{subsubpcenty}-\value{subsubovalthk}*5/4}
\MoveTo(\value{posx},\value{posy})
\MarkLoc(C+-)
\Curve(c-2+,C+-,C+-,c+1-)
%
%\setcounter{posx}{\value{subsubpcentx}-\value{psz}*9/20}
%\setcounter{posy}{\value{subsubpcenty}}
%\MoveTo(\value{posx},\value{posy})
%\DrawOval(\value{ovalthn},\value{subsubovalthk})
%
\setcounter{posx}{\value{subsubpcentx}+\value{psz}*9/20}
\setcounter{posy}{\value{subsubpcenty}}
\MoveTo(\value{posx},\value{posy})
\DrawOval(\value{subsubovalthn},\value{subsubovalthk})
\setcounter{posx}{\value{subsubpcentx}}
\setcounter{posy}{\value{subsubpcenty}-\value{psz}*9/20}
\MoveTo(\value{posx},\value{posy})
\DrawOval(\value{subsubovalthk},\value{subsubovalthn})
\setcounter{posx}{\value{subsubpcentx}}
\setcounter{posy}{\value{subsubpcenty}+\value{psz}*9/20}
\MoveTo(\value{posx},\value{posy})
\DrawOval(\value{subsubovalthk},\value{subsubovalthn})
}
\def\quartpic{
%\MoveTo(\value{pcentx},\value{pcenty})
%\Text(--$P$--)
%
\setcounter{posx}{\value{pcentx}-\value{psz}}
\setcounter{posy}{\value{pcenty}}
\MoveTo(\value{posx},\value{posy})
\DrawOval(\value{ovalthn},\value{ovalthk})
\setcounter{posx}{\value{pcentx}+\value{psz}}
\setcounter{posy}{\value{pcenty}+\value{ovalthk}}
\MoveTo(\value{posx},\value{posy})
\MarkLoc(a+1+)
\setcounter{posx}{\value{pcentx}+\value{ovalthk}}
\setcounter{posy}{\value{pcenty}+\value{psz}}
\MoveTo(\value{posx},\value{posy})
\MarkLoc(a+2+)
\setcounter{posx}{\value{pcentx}+\value{ovalthk}}
\setcounter{posy}{\value{pcenty}+\value{ovalthk}}
\MoveTo(\value{posx},\value{posy})
\MarkLoc(A++)
\Curve(a+1+,A++,A++,a+2+)
\setcounter{posx}{\value{pcentx}+\value{psz}}
\setcounter{posy}{\value{pcenty}-\value{ovalthk}}
\MoveTo(\value{posx},\value{posy})
\MarkLoc(a+1-)
\setcounter{posx}{\value{pcentx}+\value{ovalthk}}
\setcounter{posy}{\value{pcenty}-\value{psz}}
\MoveTo(\value{posx},\value{posy})
\MarkLoc(a-2+)
\setcounter{posx}{\value{pcentx}+\value{ovalthk}}
\setcounter{posy}{\value{pcenty}-\value{ovalthk}}
\MoveTo(\value{posx},\value{posy})
\MarkLoc(A+-)
\Curve(a+1-,A+-,A+-,a-2+)
%
%
%
%\MoveTo(\value{pcentx},\value{pcenty})
%\Text(--$y_1^2P$--)
%
%\setcounter{posx}{\value{pcentx}-\value{psz}*2/3}
%\setcounter{posy}{\value{pcenty}+\value{subovalthk}}
%\MoveTo(\value{posx},\value{posy})
%\MarkLoc(b-1+)
\setcounter{posx}{\value{subpcentx}-\value{subovalthk}}
\setcounter{posy}{\value{subpcenty}+\value{psz}*2/3}
\MoveTo(\value{posx},\value{posy})
\MarkLoc(b+2-)
\setcounter{posx}{\value{subpcentx}-\value{subovalthk}}
\setcounter{posy}{\value{subpcenty}+\value{subovalthk}*5/4}
\MoveTo(\value{posx},\value{posy})
\MarkLoc(B-+)
\Curve(a+1+,B-+,B-+,b+2-)
\setcounter{posx}{\value{subpcentx}-\value{subovalthk}}
\setcounter{posy}{\value{subpcenty}-\value{psz}*2/3}
\MoveTo(\value{posx},\value{posy})
\MarkLoc(b-2-)
\setcounter{posx}{\value{subpcentx}-\value{subovalthk}}
\setcounter{posy}{\value{subpcenty}-\value{subovalthk}*5/4}
\MoveTo(\value{posx},\value{posy})
\MarkLoc(B--)
\Curve(a+1-,B--,B--,b-2-)
\setcounter{posx}{\value{subpcentx}+\value{subovalthk}}
\setcounter{posy}{\value{subpcenty}+\value{psz}*2/3}
\MoveTo(\value{posx},\value{posy})
\MarkLoc(b+2+)
\setcounter{posx}{\value{subpcentx}+\value{psz}*2/3}
\setcounter{posy}{\value{subpcenty}+\value{subovalthk}}
\MoveTo(\value{posx},\value{posy})
\MarkLoc(b+1+)
\setcounter{posx}{\value{subpcentx}+\value{subovalthk}}
\setcounter{posy}{\value{subpcenty}+\value{subovalthk}*5/4}
\MoveTo(\value{posx},\value{posy})
\MarkLoc(B++)
\Curve(b+2+,B++,B++,b+1+)
\setcounter{posx}{\value{subpcentx}+\value{subovalthk}}
\setcounter{posy}{\value{subpcenty}-\value{psz}*2/3}
\MoveTo(\value{posx},\value{posy})
\MarkLoc(b-2+)
\setcounter{posx}{\value{subpcentx}+\value{psz}*2/3}
\setcounter{posy}{\value{subpcenty}-\value{subovalthk}}
\MoveTo(\value{posx},\value{posy})
\MarkLoc(b+1-)
\setcounter{posx}{\value{subpcentx}+\value{subovalthk}}
\setcounter{posy}{\value{subpcenty}-\value{subovalthk}*5/4}
\MoveTo(\value{posx},\value{posy})
\MarkLoc(B+-)
\Curve(b-2+,B+-,B+-,b+1-)

%
%\setcounter{posx}{\value{subpcentx}-\value{psz}*2/3}
%\setcounter{posy}{\value{subpcenty}}
%\MoveTo(\value{posx},\value{posy})
%\DrawOval(\value{ovalthn},\value{subovalthk})
%
\setcounter{posx}{\value{subpcentx}+\value{psz}*2/3}
\setcounter{posy}{\value{subpcenty}}
\MoveTo(\value{posx},\value{posy})
\DrawOval(\value{subovalthn},\value{subovalthk})
\setcounter{posx}{\value{subpcentx}}
\setcounter{posy}{\value{subpcenty}-\value{psz}*2/3}
\MoveTo(\value{posx},\value{posy})
\DrawOval(\value{subovalthk},\value{subovalthn})
\setcounter{posx}{\value{subpcentx}}
\setcounter{posy}{\value{subpcenty}+\value{psz}*2/3}
\MoveTo(\value{posx},\value{posy})
\DrawOval(\value{subovalthk},\value{subovalthn})
}
\setcounter{rot}{0} \Do(1,4){
%\addtocounter{pcentx}{-50}
\MoveTo(\value{pcentx},\value{pcenty})
%\addtocounter{pcentx}{50}
\setcounter{rot2}{\value{rot}+90}
\RotatedAxes(\value{rot},\value{rot2}) \quartpic \smallestpic
\EndRotatedAxes \addtocounter{rot}{90} }
%%%%%
\MoveTo(\value{subpcentx},-\value{subpcentx}) \RotatedAxes(90,180)
\smallestpic \EndRotatedAxes
\MoveTo(\value{subpcentx},\value{subpcentx}) \RotatedAxes(180,270)
\smallestpic \EndRotatedAxes
\MoveTo(-\value{subpcentx},\value{subpcentx}) \RotatedAxes(270,0)
\smallestpic \EndRotatedAxes
\MoveTo(-\value{subpcentx},-\value{subpcentx}) \RotatedAxes(0,90)
\smallestpic \EndRotatedAxes
%
%%%%%
\MoveTo(-\value{subpcentx},\value{subpcentx}) \RotatedAxes(0,90)
\smallestpic \EndRotatedAxes
\MoveTo(\value{subpcentx},\value{subpcentx}) \RotatedAxes(270,0)
\smallestpic \EndRotatedAxes
\MoveTo(\value{subpcentx},-\value{subpcentx})
\RotatedAxes(180,270) \smallestpic \EndRotatedAxes
\MoveTo(-\value{subpcentx},-\value{subpcentx})
\RotatedAxes(90,180) \smallestpic \EndRotatedAxes
%
%%%%%%%%%%%%%%%%%
\MoveTo(0,0) \Text(--\mbox{$U$}--)
\MoveTo(-50,0) \Text(--\mbox{$V_i$}--)
\MoveTo(50,0) \Text(--\mbox{$z_iV_i$}--)
\MoveTo(0,-50) \Text(--\mbox{$V_j$}--)
\MoveTo(0,50) \Text(--\mbox{$z_jV_j$}--)
\MoveTo(105,0) \Text(--\mbox{$z_iU$}--)
\MoveTo(-105,0) \Text(--\mbox{$z^{-1}_iU$}--)
\MoveTo(0,105) \Text(--\mbox{$z_jU$}--)
\MoveTo(0,-105) \Text(--\mbox{$z^{-1}_jU$}--)
\EndDraw
\end{equation*}
%%%End of picture.
}

Thus $U \subset \wwtilde{W}$ is sufficiently large for the translates
$gU \subset \wwtilde{W}$ ($g \in F_\mu$) to cover $\wwtilde{W}$, but
sufficiently small for the overlaps $gU \cap hU$ to be non-empty
only if $g^{-1}h = 1$ or $z_i$ or $z^{-1}_i$.

\begin{example} \label{funex}
(i)\qua  The subspace $(1,I_\mu) \subset EF_\mu$ is a fundamental
domain of the universal cover $p_\mu\co EF_\mu \to BF_\mu$.

(ii)\qua Let $G'_\mu$ be the barycentric subdivision of the Cayley
tree $G_\mu$, the tree with
$$\begin{array}{l}
(G'_{\mu})^{(0)} = G_\mu^{(0)} \cup G_\mu^{(1)},\\[1ex]
(G'_{\mu})^{(1)} = \{(h,(g,gz_i))\,\vert\, h  = g~{\rm or}~gz_i\}
\subset (G'_\mu)^{(0)} \times (G'_{\mu})^{(0)}.
\end{array}$$
The $F_\mu$--equivariant homeomorphism $G_\mu = G'_\mu\cong EF_\mu$
sends the vertex $(g,gz_i) \in (G'_\mu)^{(0)}$ to $(g,e^+_i) \in
EF_{\mu}$. The subgraph $U_{\mu} \subset G'_{\mu}$ defined by
$$\begin{array}{l}
U_{\mu}^{(0)} =  \{1\} \cup \{(1,z_i)\} \cup
\{(z_i^{-1},1)\}\\[1ex]
U_\mu^{(1)} = \{(1,(1,z_i))\} \cup \{(z_i^{-1},(z_i^{-1},1))\}
\end{array}$$
is the fundamental domain of the  cover $G_\mu \to G_\mu/F_\mu$
corresponding to $(1,I_\mu) \subset EF_\mu$ under the
$G_\mu$--equivariant homeomorphism $G_\mu\cong EF_\mu$.
{\small
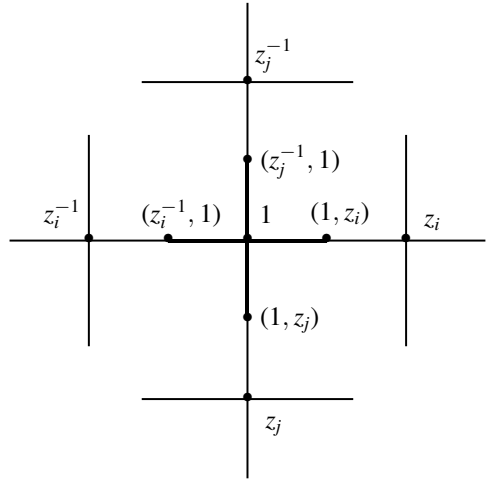
    $$\Draw
    \MoveTo(0,0) \Text(--$\bullet$--)
    \MoveTo(0,60) \Text(--$\bullet$--)
    \MoveTo(0,-60) \Text(--$\bullet$--)
    \MoveTo(60,0) \Text(--$\bullet$--)
    \MoveTo(-60,0)\Text(--$\bullet$--)
    \MoveTo(0,30) \Text(--$\bullet$--)
    \MoveTo(0,-30) \Text(--$\bullet$--)
    \MoveTo(30,0) \Text(--$\bullet$--)
    \MoveTo(-30,0)\Text(--$\bullet$--)
    \LineAt(0,-90,0,-30)
    \LineAt(0,30,0,90)
    \LineAt(-90,0,-30,0)
    \LineAt(30,0,90,0)
    {\PenSize(1.5pt)
    \LineAt(-30,0,30,0)
    \LineAt(0,-30,0,30)}
    \LineAt(-60,40,-60,-40) \LineAt(-40,-60,40,-60)
\LineAt(60,-40,60,40) \LineAt(-40,60,40,60) \MoveTo(-70,10)
\Text(--$z_i^{-1}$--) \MoveTo(-25,10) \Text(--$(z_i^{-1},1)$--)
\MoveTo(7,10) \Text(--$1$--) \MoveTo(35,10) \Text(--$(1,z_i)$--)
\MoveTo(70,7) \Text(--$z_i$--) \MoveTo(10,-70) \Text(--$z_j$--)
\MoveTo(16,-30) \Text(--$(1,z_j)$--) \MoveTo(20,30)
\Text(--$(z_j^{-1},1)$--) \MoveTo(10,70) \Text(--$z_j^{-1}$--)
    \EndDraw$$
}
\end{example}

\begin{proposition} \label{fund}
{\rm (i)}\qua Given an $F_\mu$--cover $p\co \wwtilde{W} \to W$ and a map
$f\co V \to W$ let $f^*p\co \wtilde{V} = f^*\wwtilde{W} \to V$ be the
pullback $F_\mu$--cover. If $U \subset \wwtilde{W}$ is a
fundamental domain of $p$ then
$$\wtilde{f}^{-1}(U) = \{(x,y)\,\vert\, x \in V,y \in U,f(x) = p(y) \in W\}
    \subset \wtilde{V}$$
is a fundamental domain of $f^*p$.

{\rm (ii)}\qua Every $F_{\mu}$--cover $p\co \wwtilde{W} \to W$ has
fundamental domains.
\end{proposition}
\begin{proof} (i)\qua By construction.

(ii)\qua Apply (i), using the fundamental domain $U_\mu \subset
G'_\mu = EF_\mu$ for the cover
$$p_\mu\co EF_{\mu} \to EF_\mu/F_\mu = BF_\mu$$
given by \fullref{funex}, noting that
$$p = c^*p_{\mu}\co \wwtilde{W} = c^*EF_{\mu} \to W$$
is the pullback of the universal $F_{\mu}$--cover
$p_{\mu}\co EF_{\mu} \to BF_{\mu}$ along a classifying map $c\co W \to BF_{\mu}$
$$\xymatrix@C+10pt@R+10pt{\wwtilde{W} \ar[r]^-{\di\wtilde{c}}
\ar[d]_-{\di{p}} &
EF_{\mu} \ar[d]^-{\di{p_{\mu}}} \\
    W \ar[r]^-{\di{c}} & BF_{\mu}}$$
The inverse image of  $U_{\mu} \subset EF_{\mu}$
$$U = \wtilde{c}^{-1}(U_{\mu}) \subset \wwtilde{W}$$
is a fundamental domain of $c\co \wwtilde{W} \to W$.
\end{proof}

\subsection{Combinatorial transversality}

If $p\co \wwtilde{W} \to W$ is an $F_{\mu}$--cover of a space $W$
with an additional structure such as a manifold or finite $CW$
complex, we should like to have fundamental domains $U \subset
\wwtilde{W}$ with the additional structure.  For manifolds this
is achieved by choosing a classifying map $c\co W \to BF_{\mu}$
transverse at $\{e_1,e_2,\ldots,e_{\mu}\} \subset BF_{\mu}$ -- see
\fullref{manifold transversality} below for a more detailed
discussion. For a finite $CW$ complex $W$ we shall develop a
combinatorial version of transversality,
constructing finite subcomplexes $X \subset X(\infty)$
of the Borel construction
$X(\infty) = \wwtilde{W}\times_{F_{\mu}}G_{\mu}$, such that the
projection $f(\infty)\co W(\infty) \to W$ restricts to
a simple homotopy equivalence $f\co X \to W$ such that the
pullback $F_{\mu}$--cover $\widetilde{X} = f^*\wwtilde{W} \to X$
has a fundamental domain $U \subset \widetilde{X}$ which is a
finite subcomplex.

\begin{proposition} \label{Borel} For any $F_{\mu}$--cover
$p\co \wwtilde{W} \to W$ let $F_{\mu}$ act diagonally on
$\wwtilde{W} \times G_{\mu}$
$$F_{\mu} \times (\wwtilde{W} \times G_{\mu}) \to (\wwtilde{W} \times G_{\mu});~
(g,(x,y)) \mapsto (gx,gy).$$
{\rm (i)}\qua The map
$$\pi\co X = \wwtilde{W}\times_{F_{\mu}}G_{\mu} \to W;~
[x,g] \mapsto p(x)$$
is the projection of a fibration
$$\xymatrix{G_{\mu} \ar[r] & X \ar[r]^-{\di{\pi}} & W}$$
with contractible point inverses;
for each $x \in \wwtilde{W}$ there is defined a homeomorphism
$$G_{\mu} \to \pi^{-1}p(x);~g \mapsto [x,g].$$
In particular, $\pi$ is a homotopy equivalence.

{\rm (ii)}\qua The pullback $F_{\mu}$--cover of $X$
$$\pi^*p\co \widetilde{X} = p^*\wwtilde{W} = \wwtilde{W}
\times G_{\mu} \to X = \wwtilde{W} \times_{F_{\mu}} G_{\mu}$$
has fundamental domain $\wwtilde{W} \times U \subset
\widetilde{X} = \wwtilde{W} \times G_{\mu}$, with $U \subset
G_\mu$ any fundamental domain.
\end{proposition}
\begin{proof} Standard.
\end{proof}

\begin{definition} \label{fsplit} {\rm
(i)\qua An \emph{$F_{\mu}$--splitting} $(X,Y,Z,h)$ of a space $W$ is
a homeomorphism $h\co X \to W$ from a space with a decomposition
$$X = Y \times [-1,1] \cup_{Y \times \{-1,1\}}Z$$
with $Y = Y_1\sqcup Y_2 \sqcup \ldots \sqcup Y_{\mu}$ the disjoint union
of spaces $Y_1,Y_2,\ldots,Y_{\mu}$ and
$Y \times [-1,1]$ attached to $Z$ along maps
$$\alpha^-_i\co Y_i \times \{-1\} \to Z,~
\alpha^+_i\co Y_i \times \{1\} \to Z.$$
(ii)\qua An $F_{\mu}$--splitting $(X,Y,Z,h)$ of a connected space $W$
is \emph{connected} if each of $Y_1,Y_2,\ldots,Y_{\mu},Z$ is
non-empty and connected.}
\end{definition}

\begin{proposition} \label{fcover}
Let $W$ be a space with an $F_{\mu}$--splitting $(X,Y,Z,h)$.

{\rm (i)}\qua The $F_\mu$--splitting determines an $F_{\mu}$--cover
$p\co \wwtilde{W} \to W$ with
$$\begin{array}{l}
\wwtilde{W} = \big(F_{\mu} \times (Y \times [-1,1] \sqcup
Z)\big)/\sim,\\[1ex]
\hspace*{50pt}(g,y_i,1) \sim (z_ig,\alpha^+_i(y_i,1)),\\[1ex]
\hspace*{50pt}(g,y_i,-1) \sim (g,\alpha^-_i(y_i,-1))\quad 
(g \in F_{\mu}, y_i \in Y_i,1 \leqslant i \leqslant \mu),\\[1ex]
 p\co \wwtilde{W} \to W;~(g,x) \mapsto [h(x)].
\end{array}$$
The subspace
$$Z' = (1,Y\times [0,1]) \cup (1,Z) \cup \bigcup\limits^{\mu}_{i = 1}
(z_i,Y_i \times [-1,0]) \subset \wwtilde{W}$$
is a fundamental domain of $p\co \wwtilde{W} \to W$.

{\rm (ii)}\qua If there exists a homeomorphism $\phi\co Z' \to Z$ such that
$$\phi(1,y_i,0) = \alpha^-_i(y_i,-1),~\phi(z_i,y_i,0) =
\alpha^+_i(y_i,1)\quad 
(y_i \in Y_i,1 \leqslant i \leqslant \mu)$$
the identification space
$$\wwtilde{W}' = \big(F_{\mu} \times Z\big)/(g,\alpha^-_i(y_i))
\sim (z_ig,\alpha^+_i(y_i))$$
is such that there is defined a homeomorphism
$$(1,\phi)\co \wwtilde{W} \to \wwtilde{W}';~(g,x) \mapsto (g,\phi(x))$$
so that
$$p' = p(1,\phi)^{-1}\co \wwtilde{W}' \to W;~(g,x) \mapsto p\phi^{-1}(x)$$
is an $F_{\mu}$--cover of $W$ which is isomorphic to $p\co \wwtilde{W} \to W$,
with fundamental domain
$$(1,\phi)(Z') = (1,Z) \subset \wwtilde{W}'.$$
{\rm (iii)}\qua The fundamental group of a connected space $W$ with a
connected $F_{\mu}$--splitting $(X,Y,Z,h)$ is an amalgamated free
product
$$\pi_1(W) = \pi_1(Z)*F_{\mu}/\{\alpha^+_i(g_i)z_i = z_i\alpha^-_i(g_i)\,\vert\,
    g_i \in \pi_1(Y_i),1 \leqslant i \leqslant \mu\}.$$
The surjection $\pi_1(W) \to F_{\mu}$ is induced by a map $c\co W \to
BF_{\mu}$ sending $h(Y_i \times \{0\})\subset W$ to $\{e_i\} \subset
BF_{\mu}$.  The surjection $\pi_1(W) \to F_{\mu}$ classifies the
connected $F_{\mu}$--cover $p\co \wwtilde{W} \to W$ in {\rm (i)}.
\end{proposition}
\begin{proof} (i) and (ii) follow by construction.

(iii) follows from the Seifert--van Kampen theorem and obstruction theory.
\end{proof}

\begin{example} \label{J}
Define an $F_{\mu}$--splitting $(H_\mu,\{1,2,\ldots,\mu\},I_{\mu},f)$ of
$BF_{\mu}$ by
$$\begin{array}{l}
H_{\mu} = \{1,2,\ldots,\mu\} \times [-1,1] \cup_{(i,1)\sim
e^+_i,(i,-1)\sim e^-_i}I_{\mu},\\[1ex]
f\co H_{\mu} \to BF_{\mu};~
\begin{cases}
(i,t) \mapsto [(1-t/2)e^+_i]&{\rm for}~0 \leqslant t \leqslant 1\\[1ex]
(i,t) \mapsto [(1+t/2)e^-_i]&{\rm for}~-1 \leqslant t \leqslant 0\\[1ex]
u \mapsto u/2&{\rm for}~u \in I_{\mu}
\end{cases}
\end{array}$$
with
$$f(i,0) = e_i,~f(i,1) = e^+_i/2,~f(i,-1) = e^-_i/2.$$
The corresponding $F_{\mu}$--cover of $BF_{\mu}$ is the universal
$F_\mu$--cover $\widetilde{BF}_{\mu} = G_{\mu} \to BF_{\mu}$, with
fundamental domain $I_{\mu} = (1,I_{\mu}) \subset G_{\mu}$.
Note that $f(I_{\mu}) = J_{\mu}$, with $J_{\mu} \subset I_{\mu}$ the
homeomorphic copy of $I_{\mu}$ defined by
$$J_{\mu} = \{(0,\ldots,0,t,0,\ldots,0) \in I_{\mu}\,\vert\,
-1/2 \leqslant t \leqslant 1/2\}.$$
\end{example}

A subspace $Y \subset X$ is \emph{collared} if the inclusion $i\co Y \to X$
extends to an embedding $j\co Y\times [0,1] \to X$, with $i(y) = j(y,0)\in X$
for $y \in Y$. In particular, $\partial Z \subset Z$ is collared,
for any manifold with boundary $(Z,\partial Z)$.

\begin{example}  \label{manifold transversality}
Use the $F_{\mu}$--splitting $(H_{\mu},\{1,2,\ldots,\mu\},I_{\mu},f)$
of $BF_{\mu}$ given by \fullref{J} to identify
$$BF_{\mu} = H_{\mu} = \{1,2,\ldots,\mu\} \times
[-1,1]\cup_{\{1,2,\ldots,\mu\}\times \{-1,1\}}I_{\mu}.$$
If $p\co \widetilde{X} \to X$ is an $F_{\mu}$--cover of a manifold
$X$ it is possible to choose a classifying map
$$c\co X \to BF_{\mu} = \{1,2,\ldots,\mu\} \times
    [-1,1]\cup_{\{1,2,\ldots,\mu\}\times \{-1,1\}}I_{\mu}$$
which is transverse regular at $\{e_1,e_2,\ldots,e_{\mu}\} \subset
BF_{\mu}$, with the inverse images of $e_i = (i,0) \in BF_{\mu}$
disjoint framed codimension--1 submanifolds
$$Y_i = c^{-1}(e_i) \subset X\quad(1 \leqslant i \leqslant \mu).$$
Cutting $X$ along
$$Y = c^{-1}\{e_1,e_2,\ldots,e_\mu\} = 
    Y_1 \sqcup Y_2 \sqcup \ldots \sqcup Y_{\mu} \subset X$$
there is obtained an $F_{\mu}$--splitting $(X,Y,Z,{\rm id.})$ of $X$, so that
$$X = Y \times [-1,1] \cup_{Y \times \{-1,1\}}Z$$
with $Y = Y\times \{0\} \subset X$ a framed codimension--1
submanifold, and $Z = c^{-1}(I_{\mu}) \subset X$ a codimension--0
submanifold with
$$\alpha^+_i\co Y_i \times \{1\} \to Z,\quad \alpha^-_i\co Y_i \times \{-1\} \to Z$$
components of the inclusion of the boundary $\partial Z = Y \times \{-1,1\} \subset Z$.
Since $\partial Z \subset Z$ is collared
the fundamental domain of the $F_{\mu}$--cover $\widetilde{X} = c^*G_{\mu}$
$$Z' = (1,Y\times [0,1]) \cup (1,Z) \cup \bigcup\limits^{\mu}_{i = 1}
(z_i,Y_i \times [-1,0]) \subset \widetilde{X}$$
is such that there exists a homeomorphism $\phi\co Z' \to Z$ with
$$\phi(1,y_i,0) = \alpha^-_i(y_i,-1),~\phi(z_i,y_i,0) =
\alpha^+_i(y_i,1)\quad
(y_i \in Y_i,1 \leqslant i \leqslant \mu).$$
Thus by \fullref{fcover} (ii) $p\co \widetilde{X} \to X$ is isomorphic
to the $F_{\mu}$--cover $p'\co \widetilde{X}' \to X$ with
$$\begin{array}{l}
\widetilde{X}' = \big(F_{\mu} \times Z\big)/(g,\alpha^-_i(y_i))
\sim (z_ig,\alpha^+_i(y_i)),\\[1ex]
p' = p(1,\phi)^{-1}\co \widetilde{X}' \to X;~(g,x) \mapsto p\phi^{-1}(x).
\end{array}$$
If $X$ and $\widetilde{X}$ are connected it is possible to choose $c$ such
that each $Y_i = p^{-1}(e_i)$ is connected, with
$$p_* = p(Y,Z)_*\co \pi_1(X) \to F_{\mu}.$$
\end{example}

\begin{definition}
(i)\qua A \emph{homotopy $F_{\mu}$--splitting} $(X,Y,Z,h)$ of a space $W$ is
a homotopy equivalence $h\co X \to W$ from a space with an
$F_\mu$--splitting $(X,Y,Z,1)$,  so that
$$X = Y \times [-1,1] \cup_{Y \times \{-1,1\}}Z,~
Y = Y_1 \sqcup Y_2 \sqcup \ldots \sqcup Y_{\mu}.$$
(ii)\qua A homotopy $F_{\mu}$--splitting $(X,Y,Z,h)$ of a finite $CW$
complex $W$ is \emph{simple} if $X$ is a finite $CW$ complex,
$Y_1,Y_2,\ldots,Y_{\mu},Z \subset X$ are subcomplexes and $h\co W \to X$ is a
simple homotopy equivalence.
\end{definition}

\begin{example}
Any finite $CW$ complex $W$ with an $F_{\mu}$--cover
$\wwtilde{W}\to W$ admits simple homotopy $F_{\mu}$--splittings $(X,Y,Z,h)$\,:
embed $W \subset S^N$ ($N$ large) with closed regular
neighbourhood $(X,\partial X)$ and apply the manifold
transversality of \fullref{manifold transversality} to
the $F_\mu$--cover $\widetilde{X} \simeq \wwtilde{W} \to W \simeq X$.
\end{example}

Working as in Ranicki \cite{RAC} we shall now develop a
combinatorial transversality construction of simple homotopy
$F_{\mu}$--splittings of $W$ using finite subcomplexes of
the Borel construction (\fullref{Borel})
$\wwtilde{W}\times_{F_{\mu}}G_{\mu}$, as follows.

\begin{definition} \label{canonical}
The \emph{canonical homotopy $F_{\mu}$--splitting}
$(X(\infty),Y(\infty),Z(\infty),h(\infty))$ of a space $W$
with an $F_{\mu}$--cover $p\co \wwtilde{W} \to W$ is given by
$$X(\infty) = Y(\infty) \times [-1,1] \cup_{Y(\infty)\times
    \{-1,1\}}Z(\infty)$$
with
$$\begin{array}{l}
\alpha(\infty)^+_i\co Y(\infty)_i = \wwtilde{W} \to
Z(\infty) = \wwtilde{W}\times I_{\mu};~
x \mapsto (z_ix,e^+_i),\\[1ex]
\alpha(\infty)^-_i\co Y(\infty)_i = \wwtilde{W} \to
Z(\infty) = \wwtilde{W}\times I_{\mu};~ x \mapsto
(x,e^-_i),\\[1ex]
h(\infty)\co X(\infty) \to W;~(x,y) \mapsto p(x).
\end{array}$$
The map $h(\infty)$ is a homotopy equivalence since it is the
composite
$$h(\infty) = \pi\circ f\co X(\infty) \xymatrix{\ar[r]^-{\di{f}}&}
\wwtilde{W} \times_{F_{\mu}}G_{\mu}
    \xymatrix{\ar[r]^-{\di{\pi}}&} W$$
    of the homeomorphism
$$f\co X(\infty) \to \wwtilde{W} \times_{F_{\mu}}G_{\mu};~
\begin{cases}
(x,i,t) \mapsto (x,(1-t/2)e^+_i)&{\rm for}~0 \leqslant t \leqslant 1\\[1ex]
(x,i,t) \mapsto (z_ix,(1+t/2)e^-_i)&{\rm for}~-1 \leqslant t \leqslant 0\\[1ex]
(x,u) \mapsto (x,u/2)&{\rm for}~u \in I_{\mu}
\end{cases}$$
and the homotopy equivalence
$$\pi\co X = \wwtilde{W}\times_{F_{\mu}}G_{\mu} \to W$$
given by \fullref{Borel}.  For every $y \in G_{\mu}$ there
is a unique $g \in F_{\mu}$ such that $gy \in I_{\mu} \backslash
\{e^+_1,e^+_2,\ldots,e^+_{\mu}\}$, so that either $gy = te^+_i$ with
$0 \leqslant t <1$, or $gy = te^-_i$ with $0 \leqslant t \leqslant
1$, and
$$\begin{array}{l}
f^{-1}\co \wwtilde{W} \times_{F_{\mu}}G_{\mu} \to X(\infty)~:\\[1ex]
[x,y] \mapsto
\begin{cases}
(gx,i,2(1-t))&{\rm if}~gy = te^+_i~{\rm with}~1/2 \leqslant t < 1\\[1ex]
(z^{-1}_igx,i,2(t-1))&{\rm if}~gy = te^-_i~{\rm with}~1/2 \leqslant t \leqslant 1\\[1ex]
(gx,gy)&{\rm if}~2gy \in I_{\mu}~(\hbox{\rm ie if $-1/2
\leqslant t \leqslant 1/2$}).
\end{cases}
\end{array}$$
\end{definition}

\begin{proposition}\label{xv}
Given a space $W$ with $F_{\mu}$--cover $p\co \wwtilde{W} \to W$ and a
subspace $V \subseteq \wwtilde{W}$ let
$$X(V) = Y(V) \times [-1,1] \cup_{Y(V) \times \{-1,1\}}Z(V)
    \subseteq X(\infty)$$
with
$$\begin{array}{l}
\alpha(V)^+_i\co 
Y(V)_i = V \cap z_i^{-1}V \to Z(V) = V \times I_{\mu};~
x \mapsto (z_ix,e^+_i),\\[1ex]
\alpha(V)^-_i\co Y(V)_i = V \cap z_i^{-1}V \to Z(V) = V \times
    I_{\mu};~ x \mapsto (x,e^-_i),
\end{array}$$
and set
$$h(V) = h(\infty)\vert\co X(V) \to W;~(x,t) \mapsto p(x).$$
{\rm (i)}\qua For any $x \in V$
$$\begin{array}{l}
h(V)^{-1}(p(x)) = \{(x,y) \in \wwtilde{W}\times_{F_{\mu}}G_\mu\,\vert\,
y \in G_{\mu}(V,x)\}\\[1ex]
\hphantom{h(V)^{-1}(p(x))~} = \{x\} \times G_{\mu}(V,x)
\subseteq X(V) \subseteq X(\infty) = \wwtilde{W}\times_{F_{\mu}}G_\mu
\end{array}$$
with $G_{\mu}(V,x) \subseteq G_{\mu}$ the subgraph defined by
$$\begin{array}{l}
G_{\mu}(V,x)^{(0)} = \{g \in F_{\mu}\,\vert\,gx \in V\}
\subseteq G_{\mu}^{(0)} = F_{\mu},\\[1ex]
G_{\mu}(V,x)^{(1)} = \{(i,g) \,\vert\,gx,gz_ix \in V\}
\subseteq G_{\mu}^{(1)} = \{1,2,\ldots,\mu\} \times F_{\mu}.
\end{array}$$
{\rm (ii)}\qua The image of $h(V)$ is
$$h(V)(X(V)) = p(V) \subseteq W,$$
so that $h(V)$ is surjective if and only if $p(V) = W$, if and only if
$\bigcup_{g \in F_{\mu}}gV = \wwtilde{W}$.
\end{proposition}
\begin{proof} By construction.
\end{proof}

In particular, if $V = \wwtilde{W}$ then
$$(X(V),Y(V),Z(V),h(V)) = (X(\infty),Y(\infty),Z(\infty),h(\infty))$$
and $h(V)\co X(V) = X(\infty) \to W$ is a homotopy equivalence (since
it has contractible point inverses).

\begin{thm}[Combinatorial transversality]
\label{combtrans}
Let $W$ be a connected finite $CW$
complex with a connected $F_{\mu}$--cover $p\co \wwtilde{W}\to W$.
The canonical homotopy $F_{\mu}$--splitting
$(X(\infty),Y(\infty),Z(\infty),h(\infty))$ of $W$ is a union
$$(X(\infty),Y(\infty),Z(\infty),h(\infty)) = \bigcup\limits_{\{V\}}
 (X(V),Y(V),Z(V),h(V))$$
of simple homotopy $F_{\mu}$--splittings $(X(V),Y(V),Z(V),h(V))$ of $W$,
with $\{V\}$ a collection of finite subcomplexes $V \subset \wwtilde{W}$
such that
$$p(V) = W,\qquad \bigcup\limits_{\{V\}} V = \wwtilde{W}.$$
In particular, there exist simple homotopy $F_{\mu}$--splittings of $W$.
\end{thm}

\begin{proof} Let
$$W = \bigcup D^0 \cup \bigcup D^1 \cup \ldots \cup \bigcup D^n$$
be the given cell structure of $W$, with skeleta
$$W^{(r)} = \bigcup D^0 \cup \bigcup D^1 \cup \ldots \cup \bigcup D^r.$$
The characteristic maps $D^r \to W$ of the $r$--cells restrict to embeddings
$D^r\backslash S^{r-1}\subset W$ on the interiors, and as a set
$W$ is the disjoint union of the interiors
$$W = \bigsqcup D^0 \sqcup \bigsqcup (D^1\backslash S^0) \sqcup
\ldots \sqcup \bigsqcup (D^n\backslash S^{n-1}).$$
Choose a lift of each $r$--cell $D^r$ in $W$ to an $r$--cell
$\widetilde{D}^r$ in $\wwtilde{W}$, so that
$$\wwtilde{W} = \bigcup\limits_{g \in F_\mu}\bigcup g\widetilde{D}^0  \cup
\bigcup\limits_{g \in F_\mu}\bigcup g\widetilde{D}^1 \cup \ldots \cup
\bigcup\limits_{g \in F_\mu}\bigcup g\widetilde{D}^n.$$
Write $\phi\co S^r \to W^{(r)}$ for the attaching maps of the
$(r{+}1)$--cells in $W$, and let $\widetilde{\phi}\co S^r \to \wwtilde{W}^{(r)}$
be the attaching maps of the chosen lifted $(r{+}1)$--cells in $\wwtilde{W}$.
For any subtree $T_n \subseteq G_{\mu}$ there exists a sequence
of subtrees $T_r \subseteq G_{\mu}$ for $r = n-1,n-2,\ldots,0$
such that
$$\widetilde{\phi}(S^r) \subseteq \wwtilde{W}^{(r-1)} \cup
\bigcup\limits_{g_r \in T^{(0)}_r}g_r\widetilde{D}^r.\eqno{(*)}$$
The sequence $T = (T_n,T_{n-1},\ldots,T_0)$ determines a subcomplex
$$V\langle T \rangle = \bigcup\limits_{g_0 \in T^{(0)}_0}\bigcup g_0\widetilde{D}^0 \cup
\bigcup\limits_{g_1 \in T^{(0)}_1}\bigcup g_1\widetilde{D}^1 \cup \ldots \cup
\bigcup\limits_{g_n \in T^{(0)}_n}\bigcup g_n\widetilde{D}^n \subseteq
\wwtilde{W}$$
such that $p(V\langle T \rangle) = W$.
The map $h(V\langle T \rangle)\co X(V\langle T \rangle) \to W$ constructed in
\fullref{xv} is surjective, with contractible point inverses
$$h(V\langle T \rangle)^{-1}(p(x)) = G_{\mu}(V,x) = T_r\hskip10pt
(p(x) \in D^r\backslash S^{r-1} \subset W),$$
so that it is a homotopy equivalence and
$(X(V\langle T \rangle),Y(V\langle T \rangle),Z(V\langle T \rangle),
h(V\langle T \rangle))$ is a homotopy $F_\mu$--splitting of $W$.
For the maximal sequence $T = (G_\mu,G_\mu,\ldots,G_\mu)$
$V\langle T\rangle = \wwtilde{W}$ and we have the canonical homotopy
$F_\mu$--splitting $(X(\infty),Y(\infty),
Z(\infty),\allowbreak h(\infty))$ of $W$.  Any
finite subtree $T_n \subset G_\mu$ can be used to start a sequence
$T = (T_n,T_{n-1},\ldots,T_0)$ of finite subtrees $T_r \subset G_\mu$
satisfying $(*)$, since for each $r = n,n-1,\ldots,1$ the $r$--cells
$\widetilde{D}^r \to \wwtilde{W}$ are attached to a finite subcomplex
of the $(r{-}1)$--skeleton $\wwtilde{W}^{(r-1)}$.  For a sequence $T$ of
finite subtrees $(X(V\langle T \rangle),Y(V\langle T
\rangle),\allowbreak Z(V\langle T \rangle),h(V\langle T \rangle))$ is a
simple homotopy $F_{\mu}$--splitting of $W$.  Finally, note that $G_\mu$
is a union of finite subtrees $T_n \subset G_\mu$, so that
$(F_{\mu},F_{\mu},\ldots,F_{\mu})$ is a union of sequences
$T = (T_n,T_{n-1},\ldots,T_0)$ of finite subtrees $T_r \subset G_\mu$
satisfying $(*)$, with corresponding expressions
\begin{align*}
\wwtilde{W} & =  \bigcup\limits_TV\langle T \rangle,\\
(X(\infty),Y(\infty),Z(\infty),h(\infty)) & = 
\bigcup\limits_T (X(V\langle T \rangle),Y(V\langle T \rangle),
 Z(V\langle T \rangle),h(V\langle T \rangle)).
\end{align*}
This completes the proof.
\end{proof}

\section{Algebraic transversality for $A[F_{\mu}]$--module complexes}
\label{algebraic transversality}

Algebraic transversality for $A[F_\mu]$--module chain complexes is
modelled on the combinatorial transversality for $F_\mu$--covers of
\fullref{combinatorial transversality}.  The procedure replaces
matrices with entries in $A[F_\mu]$ by (in general larger) matrices
with entries of the linear type
$$a_1+\sum\limits_{i = 1}^\mu a_{z_i}z_i \in
        A[F_\mu]\quad(a_1,a_{z_1},\ldots,a_{z_\mu} \in A).$$
Algebraic transversality can be traced back to the work of Higman,
Bass--Heller--Swan, Stallings, Casson and Waldhausen on the algebraic
$K$--theory of polynomial extensions and more general amalgamated
free products. See of Ranicki \cite[Chapter~7]{RHK} for a treatment
of algebraic transversality in the case $\mu = 1$ when
$A[F_\mu] = A[z,z^{-1}]$ is the Laurent polynomial extension of $A$.

\begin{definition}
Given an $A$--module $P$ and a set $F$ let
$$P[F] = \bigoplus\limits_{x \in F}xP$$
be the direct sum of copies $xP$ of $P$, consisting of the formal $A$--linear
combinations $\sum\limits_{x \in F}xa_x$ $(a_x \in P)$
with $\{x \in F\,\vert\, a_x \neq 0\}$ finite.
\end{definition}

In particular, if $F$ is a semigroup with 1 then $A[F]$ is a ring.

We shall be particularly concerned with the case of a free group
$F = F_\mu$ or the free semigroup $F^+_{\mu}$ on $\mu$ generators
$z_1,z_2,\ldots,z_{\mu}$. Thus $F^+_{\mu} \subset F_{\mu}$ consists
of all the products $z_{i_1}^{n_1}z_{i_2}^{n_2}\ldots
z_{i_k}^{n_k}$ with $n_1,n_2,\ldots,n_k \geqslant 0$. The rings
$A[F_\mu]$, $A[F^+_\mu]$ are free products
\begin{align*}
A[F_\mu] & =  A[z_1,z_1^{-1}]*_AA[z_2,z_2^{-1}]*_A\ldots *_AA[z_{\mu},z_{\mu}^{-1}],\\
A[F^+_\mu] & =  A[z_1]*_AA[z_2]*_A\cdots *_AA[z_{\mu}].
\end{align*}

For any ring morphism $k\co A \to B$ induction and restriction
define functors
$$\begin{array}{l}
k_!\co \{\hbox{\rm $A$--modules}\} \to \{\hbox{\rm
$B$--modules}\};~L \mapsto k_!L = B\otimes_AL,\\[1ex]
k^!\co \{\hbox{\rm $B$--modules}\} \to \{\hbox{\rm $A$--modules}\};~M
\mapsto k^!M = M
\end{array}$$
such that $k_!$ is left adjoint to $k^!$, with a natural isomorphism
$$\Hom_A(L,k^!M) \to \Hom_B(k_!L,M);~
f \mapsto (b \otimes x  \mapsto bf(x)).$$

\begin{definition}
An $A[F]$--module is \emph{induced} if it is of the form
$$P[F] = k_!P = A[F]\otimes_AP$$
for an $A$--module $P$, with $k\co A \to A[F]$ the inclusion.
\end{definition}

\begin{proposition}\label{induced}
Let $P,Q$ be $A$--modules.

{\rm (i)}\qua  There is defined a natural isomorphism of additive
groups
$$\Hom_A(P,Q[F]) \to \Hom_{A[F]}(P[F],Q[F]);~
    f  \mapsto \big(\sum\limits_{y \in F}yg_y \mapsto  \sum\limits_{y \in F}yf(g_y)\big).$$
{\rm (ii)}\qua There is defined a natural injection of additive groups
$$\Hom_A(P,Q)[F] \to \Hom_A(P,Q[F]);~
\sum\limits_{x \in F}    xf_x  \mapsto \big(y \mapsto
    \sum\limits_{x \in F}xf_x(y)\big).$$
{\rm (iii)}\qua If $P$ is a f.g.~projective $A$--module the injection
in {\rm (ii)} is also a surjection, so that the composite with
the isomorphism in {\rm (i)} is a natural isomorphism allowing the
identification
$$\Hom_A(P,Q)[F] = \Hom_{A[F]}(P[F],Q[F]).$$
\end{proposition}
\begin{proof} (i)\qua This is just the adjointness of $k_!$ and $k^!$,
with $k\co A \to A[F]$ the inclusion.

(ii)\qua Obvious.

(iii)\qua It is sufficient to consider the case $P = A$.
\end{proof}

\begin{definition} Let $P$ be an $A$--module which is given as a
$\mu$--fold direct sum
$$P = P_1 \oplus P_2 \oplus \cdots \oplus P_{\mu}$$
with idempotents $\pi_i\co P \to P_i \to P$.

(i)\qua Define the $A[F]$--module endomorphism
\begin{multline*}
z  =  \sum\limits^{\mu}_{i = 1}\pi_iz_i  = 
\begin{pmatrix} z_1 & 0 & \cdots & 0 \\
0 & z_2 & \cdots & 0 \\
\vdots & \vdots & \ddots & \vdots \\
0 & 0 & \cdots & z_{\mu} \end{pmatrix} \co
P[F]  =  P_1[F] \oplus P_2[F] \oplus \cdots \oplus P_{\mu}[F]\\
\longrightarrow P[F]  =  P_1[F] \oplus P_2[F]\oplus \cdots \oplus P_{\mu}[F].
\end{multline*}
For $F = F_{\mu}$ this is an automorphism, with inverse
{\small\begin{multline*}
z^{-1}  =  \sum\limits^{\mu}_{i = 1}\pi_iz_i^{-1}  =  
\begin{pmatrix} z_1^{-1} & 0 & \cdots & 0 \\
0 & z_2^{-1} & \cdots & 0 \\
\vdots & \vdots & \ddots & \vdots \\
0 & 0 & \cdots & z_{\mu}^{-1} \end{pmatrix}\co
P[F_{\mu}]  =  P_1[F_{\mu}] \oplus P_2[F_{\mu}] \oplus \cdots \oplus P_{\mu}[F_{\mu}]\\
\longrightarrow P[F_{\mu}]  =  P_1[F_{\mu}] \oplus P_2[F_{\mu}] \oplus
\cdots \oplus P_{\mu}[F_{\mu}].
\end{multline*}}%
(ii)\qua Given a collection of $A$--module morphisms
$$e = \{e_i \in \Hom_A(P_i,Q)\,\vert\, 1 \leqslant i \leqslant \mu\}$$
define the $A[F]$--module morphism
{\small$$
ez  =  \!\sum^{\mu}_{i = 1}\!e\pi_iz_i  = 
\begin{pmatrix} e_1 z_1 & e_2 z_2 & \cdots & e_{\mu}z_{\mu}
\end{pmatrix}\co
P[F]  =  P_1[F] \oplus P_2[F] \oplus \cdots \oplus P_{\mu}[F] \to Q[F].$$}%
(iii)\qua An $A[F]$--module morphism $f\co P[F] \to Q[F]$ is \emph{linear} if
$$\begin{array}{l}
f = f^+z - f^- = \begin{pmatrix} f^{+,1}z_1-f^{-,1} & \ldots &
f^{+,\mu}z_\mu-f^{-,\mu} \end{pmatrix}~:\\[1ex]
\hskip100pt
P[F] = P_1[F]\oplus P_2[F] \oplus \cdots \oplus P_{\mu}[F] \to Q[F]
\end{array}$$
for some $A$--module morphisms $f^{+,i},f^{-,i}\co P_i \to Q$.
\end{definition}

\begin{definition}
(i)\qua A \emph{Mayer--Vietoris presentation} of an
$A[F]$--module $E$ is an exact sequence of the type
$$\xymatrix{
0 \ar[r] & {{\bigoplus\limits_{i = 1}^\mu C^{(i)}[F]}}
\ar[r]^-{\di{f}} & D[F] \ar[r] & E \ar[r] & 0 }$$
with $C^{(i)},D$ $A$--modules and $f = f^+z-f^-$ a linear
$A[F]$--module morphism.

(ii)\qua A \emph{Mayer--Vietoris presentation} of an $A[F]$--module morphism
$\phi\co E \to E'$ is a morphism of Mayer--Vietoris presentations
$$\xymatrix{
0 \ar[r] & {{\bigoplus\limits_{i = 1}^\mu C^{(i)}[F]}}\ar[r]^-{\di{f}}
\ar[d]^{\di{\oplus g^{(i)}}} & D[F] \ar[r]\ar[d]^{\di{h}} &
E \ar[r]\ar[d]^{\di{\phi}} & 0  \\
0 \ar[r] & {\bigoplus\limits_{i = 1}^\mu
C^{\prime(i)} [F]}\ar[r]^-{\di{f'}} & D'[F] \ar[r] & E' \ar[r] & 0
}$$
where $g^{(i)}\co C^{(i)} \to {C'}^{(i)}$ and $h\co D \to D'$ $A$--module morphisms.

(iii)\qua A \emph{Mayer--Vietoris presentation} of an $A[F]$--module
chain complex $E$ is an exact sequence as in (i),
with $C^{(i)},D$ $A$--module chain complexes and
$f^{\pm,i}\co C^{(i)} \to D$ $A$--module chain maps.
Similarly for an $A[F]$--module chain map $\phi\co E \to E'$, with a
morphism of exact sequences as in (ii)
in which $g^{(i)}$, $h$ are $A$--module chain maps.

(iv)\qua A Mayer--Vietoris presentation of a finite induced
f.g.~projective $A[F_\mu]$--module
chain complex $E$ is \emph{finite} if $C^{(i)},D$ are finite f.g.~projective
$A$--module chain complexes.
\end{definition}

\begin{example} Let $X$ be the $CW$ complex
$$X = Z/\bigl\{x \sim \beta_i(x)\,\vert\,x \in Y^+_i,1 \leqslant i
    \leqslant \mu\bigr\}$$
which is obtained from a $CW$ complex $Z$ and  disjoint collared
subcomplexes
$$Y^+_1,~Y^+_2,\ldots,~Y^+_{\mu},~Y^-_1,~Y^-_2,\ldots,~Y^-_{\mu}
    \subset Z$$
using cellular homeomorphisms $\beta_i\co Y^+_i \to Y^-_i$ as
identifications. As in \fullref{fsplit} there is an $F_\mu$--splitting
$(X,Y,Z,h)$, where $Y = Y^+_1\sqcup Y^+_2 \sqcup \ldots \sqcup
Y^+_\mu$ and
$$\begin{array}{l}
\alpha^+_i = {\rm inclusion}_{Y^+_i \subset Z}\co Y_i = Y^+_i \to Z,\\[1ex]
\alpha^-_i = ({\rm inclusion}_{Y^-_i \subset
Z})\beta_i\co Y_i = Y^+_i \to Z.
\end{array}$$
The cellular free $\Z[F_\mu]$--module chain complex $C(\widetilde{X})$
of the $F_\mu$--cover $\widetilde{X}$ of $X$ given by \fullref{fcover} (i)
has a Mayer--Vietoris presentation
$$\xymatrix{0 \ar[r]& C(Y)[F_\mu] \ar[r]^-{\di{\alpha}}&
    C(Z)[F_\mu] \ar[r]& C(\widetilde{X})  \ar[r]& 0}$$
with $C(Y)^{(i)} = C(Y_i)$, $C(Z)$ free $\Z$--module chain complexes,
and $\alpha = \alpha^+z-\alpha^-$ a linear $\Z[F_\mu]$--module chain map.
If $Z$ is a finite $CW$ complex the
Mayer--Vietoris presentation is finite.
\end{example}

We shall construct Mayer--Vietoris presentations of free
$A[F_\mu]$--module chain complexes using the Cayley tree $G_{\mu}$
(\fullref{cayley}) and the subtree $G^+_{\mu} \subset G_{\mu}$
corresponding to $F^+_{\mu} \subset F_{\mu}$.

\begin{definition}
(i)\qua Let $G_{\mu}^+ \subset G_{\mu}$ be the subtree with
$$(G_{\mu}^+)^{(0)} = F^+_{\mu},\quad(G_{\mu}^+)^{(1)} = \{(g,gz_i)\,\vert\,g \in F^+_{\mu},
    1 \leqslant i \leqslant \mu\}.$$
(ii)\qua For any subtree $T \subseteq G_{\mu}$ and
$i = 1,2,\ldots,\mu$ let $T^{(i,1)} \subseteq T^{(1)}$ be the set of
edges of type $(g,gz_i)$ with $g \in F_{\mu}$, such that
$$T^{(1)} = \smash{\coprod^{\mu}_{i = 1}}T^{(i,1)},$$
and let
    $$T^+  =  T \cap G_{\mu}^+ \subseteq T.$$
(iii)\qua For $F = F_\mu$ (resp. $F^+_\mu$) let $G = G_{\mu}$ (resp.
$G^+_{\mu}$).
\end{definition}

We shall only be considering subtrees $T \subseteq G$ containing
the vertex $1 \in G^{(0)}$.

\begin{proposition}\label{MVprop}
Given an $A$--module $P$ let $E = P[F]$ be the induced
$A[F]$--module, regarded as a 0--dimensional $A[F]$--module chain complex.

{\rm (i)}\qua For any subtree $T \subseteq G$ there is defined a
Mayer--Vietoris presentation of $E$
$$\xymatrix{E\langle T \rangle\co
0 \ar[r] & {\smash{\bigoplus\limits_{i = 1}^\mu C^{(i)}[F]}} \ar[r]^-{\di{f}} &
D[F] \ar[r] & E \ar[r] & 0 }$$ with
$$\begin{array}{l}
D = P[T^{(0)}],\quad C^{(i)} = D \cap
z^{-1}_iD =  P[T^{(i,1)}] \subseteq E,\\[1ex]
f^{+,i}\co C^{(i)} \to D;~ xp \mapsto xp,\quad f^{-,i}\co C^{(i)}
\to D;~ xp \mapsto z_ixp.\end{array}$$
{\rm (ii)}\qua The Mayer--Vietoris presentations $E\langle T \rangle$
are such that
$$E\langle T \cap T' \rangle  =  E\langle T \rangle \cap E\langle T' \rangle,
  \quad E\langle T \cup T' \rangle = E\langle T \rangle + E\langle T'
    \rangle \quad(T,T' \subseteq G).$$
If $P$ is f.g.~projective and $T$ is finite then $C^{(i)}$, $D$
are f.g.~projective $A$--modules.

{\rm (iii)}\qua  Given a morphism of induced $A[F]$--modules
$$\phi\co E = P[F]\to E' = P'[F]$$
and a subtree $T \subseteq G$ let $\phi_*T\subseteq G$ be the
smallest subtree such that
$$\phi(P) \subseteq P'[\phi_*T^{(0)}] \subseteq E'.$$
For any subtree $T' \subseteq G$ such that $\phi_*T \subseteq T'$
there is defined a morphism of Mayer--Vietoris presentations
$$\xymatrix{ E\langle T \rangle~:~ 0 \ar[r] &
{{\bigoplus\limits_{i = 1}^\mu C^{(i)}[F]}}\ar[r]^-{\di{f}}
\ar[d]^{\bigoplus \di{g^{(i)}}} & D[F] \ar[r]\ar[d]^{\di{h}} & E
\ar[r]\ar[d]^-{\di{\phi}} & 0 \\
E'\langle T' \rangle~:~ 0 \ar[r] & {\bigoplus\limits_{i = 1}^\mu
C^{\prime(i)} [F]}\ar[r]^-{\di{f'}} & D'[F] \ar[r] & E' \ar[r] &0}$$
with
$$g^{(i)} = \phi\vert\co C^{(i)} \to {C'}^{(i)},\quad h = \phi\vert\co D \to D'.$$
If $P$ is a f.g.~$A$--module and $T \subset G$ is finite, then so
is $\phi_*T \subset G$.
\end{proposition}
\begin{proof} By construction.
\end{proof}

\begin{example}
The Mayer--Vietoris presentation of $E$ associated to the
minimal subtree $T = \{1\}\subset G$ is
$$\xymatrix{0 \ar[r] & 0 \ar[r] & P[F] \ar[r]^-{\rm id.}& E \ar[r]
&0.}$$
\end{example}

\begin{definition} \label{algcan}
The \emph{canonical Mayer--Vietoris presentation} of an
$A[F]$--module chain complex $E$ with each $E_r = P_r[F]$ an
induced $A[F]$--module
$$\xymatrix{E\langle \infty \rangle~:~
0 \ar[r] & \bigoplus\limits^{\mu}_{i = 1} C^{(i)}[F] \ar[r]^-{\di{f}} &
    D[F] \ar[r]& E \ar[r] &0}$$
is the Mayer--Vietoris presentation with $E_r\langle \infty \rangle
 = E_r\langle T \rangle$ the Mayer--Vietoris presentation of $E_r$
associated to the maximal subtree $T = G \subseteq G$, where
$$f^{+,i} = {\rm id.},~f^- = z_i~:~C^{(i)}  = k^!E \to D = k^!E$$
with $k\co A \to A[F]$ the inclusion.
\end{definition}

\begin{remark} (i)\qua The canonical Mayer--Vietoris presentation
can be written in terms of induction and restriction
$$\xymatrix{E\langle \infty \rangle~:~
0 \ar[r] & \lower16pt\hbox{$\bigoplus\limits_{\mu}k_!k^!E$}
\ar[r]^{\di{f}} & k_!k^!E \ar[r]& E \ar[r] &0}$$
with
$$\begin{array}{l}
f\co \bigoplus\limits_{\mu}k_!k^!E \to k_!k^!E;~
x_i \otimes y \mapsto x_iz_i\otimes y-x_i\otimes z_iy\quad (x_i \in A[F],y \in E),\\[1ex]
k_!k^!E \to E;~ x \otimes y \mapsto xy\quad (x \in A[F],y \in E).
\end{array}$$
(ii)\qua The canonical Mayer--Vietoris presentation for $F = F_\mu$ is the
algebraic analogue of the canonical homotopy $F_\mu$--splitting of a
space $W$ with an $F_\mu$--cover $\wwtilde{W}$ in \fullref{canonical}.
\end{remark}

\begin{thm}[Algebraic transversality for chain complexes]\label{MV}
Let $E$ be an $n$--dimensional $A[F]$--module chain complex
$$\xymatrix{E\co E_n \ar[r]^-{d_n} & E_{n-1} \ar[r] &\cdots \ar[r]&
    E_1 \ar[r]^-{d_1} & E_0}$$
with each $E_r = P_r[F]$ induced from an $A$--module $P_r$.

{\rm (i)}\qua For any sequence $T = (T_n,T_{n-1},\ldots,T_0)$ of subtrees
$T_r \subseteq G$ such that
$$(d_r)_*(T_r) \subseteq T_{r-1}\quad(r = n,n-1,\ldots,1)\eqno{(*)}$$
there is defined a Mayer--Vietoris presentation
$$\xymatrix@C+15pt{E\langle T \rangle\co
0 \ar[r] & {{\bigoplus\limits_{i = 1}^\mu C^{(i)}[F]}}
\ar[r]^-{{f^+z{-}f^-}} & D[F] \ar[r] & E \ar[r] & 0
}$$
with
$$E\langle T \rangle_r = E_r \langle T_r \rangle
\quad(0 \leqslant r \leqslant n),\quad
    E\langle T \rangle \subseteq E\langle \infty \rangle.$$
{\rm (ii)}\qua If the $A$--modules $P_r$ are f.g.~projective then for
any finite subtree $T_n \subseteq G$ there exists a sequence
$T = (T_n,T_{n-1},\ldots,T_0)$ of finite subtrees $T_r \subseteq G$
satisfying $(*)$, so that $E\langle T \rangle$ is a finite
Mayer--Vietoris presentation of $E$. Thus
$$E\langle \infty \rangle = \bigcup\limits_T E \langle T \rangle $$
with the union taken over all such sequences $T$. In particular,
$E$ admits a finite Mayer--Vietoris presentation.
\end{thm}
\begin{proof} By repeated applications of \fullref{MVprop},
with the sequences $T = (T_n,T_{n-1},\allowbreak\ldots,T_0)$ the
chain complex analogues of the sequences used to construct the homotopy
$F_\mu$--splittings of $CW$ complexes in the proof of \fullref{combtrans}.
\end{proof}

This completes the proof of \fullref{thm1} of the Introduction.

\section{Blanchfield and Seifert modules}\label{modules}

\subsection{The Magnus--Fox embedding}

This section obtains some technical results on the Magnus--Fox embedding
which we shall need to characterize Blanchfield $A[F_\mu]$--modules, and
to approximate h.d.~1 $F_\mu$--link modules by f.g.~projective Seifert
$A$--modules.

Let $A\llangle x_1,x_2,\ldots,x_{\mu}\rrangle$ be the ring of
$A$--coefficient formal power series in non-commuting indeterminates
$x_1,x_2,\ldots,x_\mu$. The \emph{Magnus--Fox embedding} is defined by
$$i\co A[F_{\mu}] \to \widehat{A[F_\mu]} = 
A\llangle x_1,x_2,\ldots,x_{\mu}\rrangle ;~z_j \mapsto 1+x_j.$$
See the paper of Ara and Dicks \cite{AD} for a recent
account of the Magnus--Fox embedding, including the relationship with
noncommutative Cohn localization.

The augmentations $\epsilon(z_j) = 1$, $\widehat\epsilon(x_j) = 0$ give rise to a
commutative triangle of rings
$$
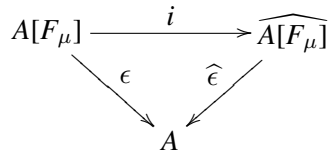
\xymatrix{ A[F_{\mu}] \ar[rr]^-{\di{i}} \ar[dr]^-{\di{\epsilon}}
&&  \widehat{A[F_\mu]}  \ar[dl]_-{\di{\widehat{\epsilon}}}\\
&A &}$$

\begin{proposition}\label{Magnus-Foxon}
{\rm (i)}\qua For projective $\widehat{A[F_\mu]}$--modules
$\what{K},\what{L}$ the augmentation map
$$\what{\epsilon}\co
\Hom_{\scriptsize\widehat{A[F_\mu]}}\bigl(\what{K},\what{L}\bigr)
\to\Hom_A\bigl(A\otimes_{\scriptsize\widehat{A[F_\mu]}}\what{K},
A\otimes_{\scriptsize\widehat{A[F_\mu]}}\what{L}\bigr);~
\what{f} \mapsto 1 \otimes \what{f}$$
is surjective.

{\rm (ii)}\qua  A morphism $\what{f}\co \what{K} \to \what{L}$ of
projective $\widehat{A[F_\mu]}$--modules is an isomorphism if and only
if the $A$--module morphism
$$1\otimes\what{f}\co A\otimes_{\scriptsize\widehat{A[F_\mu]}}\what{K}
\to A\otimes_{\scriptsize\widehat{A[F_\mu]}}\what{L}$$
is an isomorphism.

{\rm (iii)}\qua A morphism $f\co K \to L$ of projective $A[F_\mu]$--modules
induces an $\widehat{A[F_\mu]}$--module isomorphism
$$1\otimes f\co \widehat{A[F_\mu]}\otimes_{A[F_\mu]}K \to
\widehat{A[F_\mu]}\otimes_{A[F_\mu]}L$$
if and only if the $A$--module morphism
$$1\otimes f\co A\otimes_{A[F_\mu]}K \to A\otimes_{A[F_\mu]}L$$
is an isomorphism.
\end{proposition}
\begin{proof}
(i)\qua By additivity this reduces to the special case
$\what{K} = \what{L} = \widehat{A[F_\mu]}$, which is just the
fact that $\widehat{\epsilon}\co \widehat{A[F_\mu]} \to A$ is surjective.

(ii)\qua It suffices to prove that if $1\otimes\what{f}$ is an $A$--module isomorphism
then $\what{f}$ is an $\widehat{A[F_\mu]}$--module isomorphism.

Consider first the special case when $\what{K},\what{L}$
are free $\widehat{A[F_\mu]}$--modules, say $\smash{\widehat{A[F_\mu]}}^k$,
$\smash{\widehat{A[F_\mu]}}^\ell$ for some sets $k,\ell$. The augmentation map
$$\widehat{\epsilon}\co 
\Hom_{\widehat{A[F_\mu]}}\bigl(\smash{\widehat{A[F_\mu]}}^k,
\smash{\widehat{A[F_\mu]}}^{\ell}\bigr)
\to\Hom_A\bigl(A^k,A^{\ell}\bigr);~
\what{f} \mapsto 1 \otimes \what{f}$$
has a canonical splitting. If $1\otimes \what{f}$ is an isomorphism then
all the entries in the matrix of the $\widehat{A[F_\mu]}$--module morphism
$$g = 1-(1\otimes \what{f})^{-1}\what{f}\co \smash{\widehat{A[F_\mu]}}^k
\to \smash{\widehat{A[F_\mu]}}^k$$
have constant term 0, so that $1-g = (1\otimes \what{f})^{-1}\what{f}$
is an $\widehat{A[F_\mu]}$--module isomorphism with inverse
$$(1-g)^{-1} = 1+g+g^2+g^3+g^4+\cdots\co \smash{\widehat{A[F_\mu]}}^k
\to \widehat{A[F_\mu]}^k,$$
and $\what{f} = (1\otimes \what{f})(1-g)$ is an isomorphism.

For the general projective case apply (i) to lift $(1\otimes\what{f})^{-1}$
to an $\widehat{A[F_\mu]}$--module morphism $\what{e}\co \what{L}
\to \what{K}$. Choose a projective $\widehat{A[F_\mu]}$--module
$\what{J}$ such that $\what{J} \oplus \what{K} \oplus
\what{L}$ is a free $\widehat{A[F_\mu]}$--module, and apply the
special case to the $\widehat{A[F_\mu]}$--module morphism
$$1 \oplus \begin{pmatrix} 0 & \smash{\what{e}} \\ \smash{\what{f}} & 0 \end{pmatrix}\co 
\what{J} \oplus \what{K} \oplus
\what{L} \to \what{J} \oplus \what{K} \oplus
\what{L}.$$
(iii)\qua This is a special case of (ii).
\end{proof}

For $j = 1,2,\ldots,\mu$ let $y_j$ be a formal square root of $z_j$, so
that $(y_j)^2 = z_j$.  Let $F_\mu(y)$ be the free group generated by
$y_1,y_2,\ldots,y_\mu$, so that $F_\mu \subset F_\mu(y)$ is the free
subgroup generated by $z_1,z_2,\ldots,z_\mu$.  We can identify
$G_\mu^{(1,j)}$ with the subset $F_\mu y_j \subset F_\mu(y)$: the edge
$(g,gz_j) \in G_\mu^{(1,j)}$ ($g \in F_\mu$) is identified with the
element $gy_j^{-1} \in F_\mu(y)$.

\begin{lemma}\label{decompose_tree}
If $T\subset G_\mu$ is a finite subtree then
$$A[T^{(0)}] =  A[\{1\}] \oplus \Bigl(\bigoplus\limits_{j = 1}^\mu
A[T^{(1,j)}](y_j^{-1}-y_j)\Bigr) \subset A[F_\mu].\eqno{(*)}$$
\end{lemma}
\begin{proof} If $w\in T^{(1,j)}$ then certainly
$w(y_j^{-1}-y_j)\in A[T^{(0)}]$.
Let us check linear independence of the generators on the right hand
side of~$(*)$. Assuming the contrary, let
$$a_1+
\sum\limits_{gy_j^{-1}\in U} a_g g(y_j^{-1}-y_j) = 0 \in A[F_\mu]
$$
be a non-trivial relation with $U\subset T^{(1)}$ non-empty and
minimal.  We reach a contradiction by
observing that if $g(y_j)^{-1}\in U$ is a word of maximal length (in reduced
form) then $a_g = 0$.

We must also show that every $v\in T^{(0)}$ is an element of the
right-hand side of~$(*)$. Indeed there is a (unique)
path in the tree from $1$ to $v$ defined by a sequence of edges
$w_1,w_2,\ldots,w_n \in T^{(1)}$ and we have
$$v = 1+\sum_{i = 1}^nw_i(y_{j(i)}^{-1}-y_{j(i)})\eta_i \in A[F_\mu]$$
if the signs $\eta_i\in\{\pm1\}$ are chosen appropriately and $j(i)$
is such that $w_i\in T^{(1,j(i))}$.
\end{proof}

\begin{proposition} \label{Magnus-Foxtw}
For any finite subset $S \subset F_\mu$ the inclusion
$i\vert\co A[S] \to \widehat{A[F_\mu]}$ is a split $A$--module injection.
\end{proposition}
\begin{proof}
Since every finite $S$ is contained in the vertex set of
some finite tree we may assume that $S = T^{(0)}$ for some finite subtree
$T\subset G_\mu$. We proceed by induction on $|T^{(0)}|$.

If the tree $T$ has only one vertex then $T^{(0)} = \{1\}$
with $i(1) = 1 \in \widehat{A[F_\mu]}$ and
$$\widehat{A[F_\mu]} = A[\{1\}]\oplus \bigoplus\limits^{\mu}_{i = 1}
\widehat{A[F_{\mu}]}x_i~
 = A\oplus \bigoplus\limits^{\mu}_{i = 1}\widehat{A[F_\mu]}(1-z_i^{\eta})
\eqno{(**)}$$
for any $\eta \in \{\pm 1\}$,
and $i|\co A[\{1\}] \to \widehat{A[F_\mu]}$ is a split injection.

Suppose now that $|T^{(0)}|\geqslant 2$. Let $v_0\in T^{(0)}$ be a leaf,
ie a vertex to which only one edge is incident. Let $T\backslash \{v_0\}$
denote the tree obtained by removing the vertex $v_0$ and the incident
edge. By the inductive hypothesis, $i|\co A[T^{(0)}\backslash \{v_0\}] \to
\widehat{A[F_\mu]}$ is a split injection; we denote the image by $P$.

Since $v_0$ is incident to precisely one edge then
$v_0 = w_0y_k^\eta$ for unique $\eta\in\{\pm1\}$,
$k\in\{1,\ldots,\mu\}$ and $w_0\in T^{(1,k)}$.
Now for every $j$ we have
$T^{(1,j)}y_j^{-\eta} \subset T^{(0)}\backslash \{v_0\}$. Thus
$$\begin{array}{ll}
T^{(1,j)}(y_j^{-1}-y_j)& = T^{(1,j)}y_j^{-\eta}(1-y_j^{2\eta})\eta\\[1ex]
& = T^{(1,j)}y_j^{-\eta}(1-z_j^\eta)\eta \subset
(T^{(0)}\backslash \{v_0\})(1-z_j^\eta)\eta.
\end{array}$$
It follows from~$(*)$ that $i(A[T^{(0)}])$ is a direct summand of
$$Ai(v_0)\oplus \bigoplus\limits_{j = 1}^\mu P(1-z_j^\eta)$$
and hence, by the following \fullref{decompose_widehatA[F]}, a direct
summand of $\widehat{A[F_\mu]}$.
\end{proof}

\begin{lemma}\label{decompose_widehatA[F]}
Suppose $P$ is an $A$--module which is a direct summand of $\widehat{A[F_\mu]}$.
If $\theta\in \widehat{A[F_\mu]}$ is an element such that
$\widehat{\epsilon}(\theta) = 1\in A$ and $\eta = 1$ or $-1$ then
$$A\theta\oplus\left(\bigoplus\limits_{j = 1}^\mu P(1-z_j^\eta) \right) \subset
\widehat{A[F_\mu]}$$
is again a direct summand.
\end{lemma}
\begin{proof}
We may write $\widehat{A[F_\mu]} = P\oplus Q$ for some $A$--module $Q$. Let
$\eta = 1$ or $-1$. Now it follows easily from $(**)$ that
$$\begin{array}{ll}
\widehat{A[F_\mu]}& = A\theta\oplus
            \left(\bigoplus\limits_{j = 1}^\mu\widehat{A[F_\mu]}(1-z_j^\eta)\right)\\[3ex]
            & = A\theta \oplus \left(\bigoplus\limits_{j = 1}^\mu
            P(1-z_j^\eta)\right) \oplus
            \left(\bigoplus\limits_{j = 1}^\mu Q(1-z_j^\eta)\right)
\end{array}$$
which completes the proof.\end{proof}

\subsection{Blanchfield modules}\label{bla}

\begin{definition} \label{Fdefinition}
(i)\qua A \emph{Blanchfield $A[F_\mu]$--module} $M$ is an $A[F_\mu]$--module
such that
$$\Tor_*^{A[F_\mu]}(A,M) = 0.$$
(ii)\qua (Sheiham \cite{Sh2}) An \emph{$F_\mu$--link module} $M$ is an $A[F_\mu]$--module
which has a 1--dimensional induced $A[F_\mu]$--module resolution
$$\xymatrix{
0 \ar[r]&P[F_\mu] \ar[r]^-{\di{d}}& P[F_\mu] \ar[r] & M \ar[r] & 0}$$
with $P$ an $A$--module and $d$ an $A[F_\mu]$--module morphism such that
the augmentation $A$--module morphism $\epsilon(d)\co P \to P$ is an isomorphism.
\end{definition}

As before, let $k\co A \to A[F_\mu]$ be the inclusion.

\begin{proposition} \label{B = F} The following conditions on
an $A[F_\mu]$--module $M$ are equivalent:
\begin{itemize}
\item[{\rm (i)}] $M$ is a Blanchfield module,
\item[{\rm (ii)}] $M$ is an $F_\mu$--link module,
\item[{\rm (iii)}] the $A$--module morphism
$$\gamma_M\co \bigoplus \limits_\mu k^!M \to k^!M;~
(m_1,m_2,\ldots,m_\mu)\mapsto \sum\limits^{\mu}_{i = 1}(z_i-1)m_i$$
is an isomorphism.
\end{itemize}
\end{proposition}
\begin{proof}
The canonical Mayer--Vietoris presentation (\fullref{algcan}) of any
$A[F_\mu]$--module $M$ is defined by
$$\xymatrix{0 \ar[r] & \lower16pt\hbox{$\bigoplus\limits_{\mu}k_!k^!M$}
\ar[r]^{\di{d}} & k_!k^!M \ar[r]& M \ar[r] &0}$$
with
$$\begin{array}{l}
d\co \bigoplus_{\mu}k_!k^!M \to k_!k^!M;~
x_i \otimes y \mapsto x_iz_i\otimes y-x_i\otimes z_iy
\quad (x_i \in A[F_\mu],y \in M),\\[1ex]
k_!k^!M = k^!M[F_\mu] \to M;~ x \otimes y \mapsto xy
\quad (x \in A[F_\mu],y \in M),
\end{array}$$
such that $d$ has augmentation $A$--module morphism
$$\epsilon(d) = -\gamma_M\co \bigoplus\limits_{\mu}k^!M \to k^!M.$$
Regarded as a right $A[F_\mu]$--module $A$ has a 1--dimensional f.g.~free resolution
$$\xymatrix@C+25pt{
0 \ar[r]& \bigoplus\limits_{i = 1}^{\mu} A[F_\mu] \ar[r]^-{\di{\oplus
(z_i-1)}}& A[F_\mu] \ar[r]^-{\di{\epsilon}} & A \ar[r] &0,}$$
so that for any $A[F_\mu]$--module $M$
$$\Tor^{A[F_\mu]}_n(A,M) = \begin{cases}
A\otimes_{A[F_{\mu}]}M = \coker(\gamma_M)&{\rm if}~n = 0,\\
\ker(\gamma_M)&{\rm if}~n = 1,\\
0&{\rm if}~n \geqslant 2.
\end{cases}$$
The equivalences (i) $\Longleftrightarrow$ (ii)
$\Longleftrightarrow$ (iii) are now clear.
\end{proof}

\begin{definition} \label{Bdefinition}
(i)\qua Let $\Bla_\infty(A)$ be the category of Blanchfield
$A[F_\mu]$--modules, and let $\Bla(A) \subset \Bla_\infty(A)$ be the
full subcategory of the h.d.~1 Blanchfield $A[F_\mu]$--modules.
(In view of \fullref{B = F} $\Bla_\infty(A)$ is the same as the $F_\mu$--link
module category $\Flk_{\infty}(A)$ of Sheiham \cite{Sh2}).

(ii)\qua Let $\Flk(A) \subset \Bla(A)$ be the full
subcategory of the h.d.~1 Blanchfield $A[F_\mu]$--modules $M$
such that there exists a 1--dimensional induced $A[F_\mu]$--module resolution
$$\xymatrix{
0 \ar[r]&P[F_\mu] \ar[r]^-{\di{d}}& P[F_\mu] \ar[r] & M \ar[r] & 0}$$
with $P$ a f.g.~projective $A$--module.
\end{definition}

\begin{example}
(i)\qua For a principal ideal domain $A$
$$K_0(A[F_\mu]) = K_0(A) = \Z$$
(see Bass \cite{B1}) so that
$$\Bla(A) = \Flk(A).$$
(ii)\qua A finitely presented Blanchfield $\Z[F_\mu]$--module is a
`type $L$' $\Z[F_\mu]$--module in the sense of Sato \cite{Sa}.

(iii)\qua Given a $\mu$--component boundary link $\ell\co \bigsqcup_\mu
S^n \subset S^{n+2}$ let $c\co W \to W_0$ be a $\Z$--homology equivalence
from the exterior $W$ to the exterior $W_0$ of the trivial
$\mu$--component boundary link $\ell_0\co \bigsqcup_\mu S^n \subset
S^{n+2}$, with $F_\mu$--equivariant lift $\smash{\wtilde{c}\co \wwtilde{W}
\to \wwtilde{W}_0}$ to the $F_\mu$--covers.
The homology groups $\smash{\dot H_*\bigl(\wwtilde{W}\bigr) = H_{*+1}\bigl(
\wtilde{c}\co \wwtilde{W} \to \wwtilde{W}_0\bigr)}$
are Blanchfield $\Z[F_\mu]$--modules
of homological dimension $\leqslant 2$.
Each $\dot H_r\bigl(\wwtilde{W}\bigr)$ has a $\Z$--con\-trac\-ti\-ble
f.g.~free $\Z[F_\mu]$--module resolution of the type
$$0 \to \Z[F_\mu]^{a_r} \to \Z[F_\mu]^{b_r} \to \Z[F_\mu]^{c_r}
\to \dot H_r\bigl(\wwtilde{W}\bigr) \to 0\quad(0 \leqslant r \leqslant n+1)$$
with $a_r-b_r+c_r = 0$, and
$\dot H_r(\wwtilde{W})/\Z\hbox{\rm -torsion}$ is an h.d.~1 $F_\mu$--link
module (Levine \cite[3.5]{L2} for $\mu = 1$, Sato \cite[3.1]{Sa} and Duval
\cite[4.1]{Du} for $\mu \geqslant 2$).
See \fullref{chain} below for the construction of an
$(n{+}1)$--dimensional chain complex $C$ in $\Sei(\Z)$ such that
the covering $B(C)$ is an $(n{+}1)$--dimensional chain complex
in $\Flk(\Z)$ with $H_*(B(C)) = \dot H_*(\wwtilde{W})$.
\end{example}

The following \fullref{char} characterizes Blanchfield $A[F_\mu]$--modules
in terms of $A[F_\mu]$--modules $K$ such that
$$\Tor_1^{A[F_\mu]}(A,K) = 0.$$
If $K$ is a flat $A[F_\mu]$--module then ${\rm
Tor}_1^{A[F_\mu]}(B,K) = 0$ for any right $A[F_\mu]$--module $B$, and in
particular $B = A$.
If $K = P[F_\mu]$ is induced from an $A$--module $P$ then
$$\Tor_1^{A[F_\mu]}(A,P[F_\mu]) = \Tor_1^A(A,P) = 0.$$

\begin{proposition} \label{char}
{\rm (i)}\qua If $M$ is an $A[F_\mu]$--module with a resolution
$$\xymatrix{0 \ar[r] & K \ar[r]^-{\di{d}} &L\ar[r] & M\ar[r] & 0}$$
such that
$$\Tor^{A[F_{\mu}]}_1(A,K) = \Tor^{A[F_{\mu}]}_1(A,L) = 0$$
(e.g. the canonical Mayer--Vietoris presentation of \fullref{algcan}) then $M$
is Blanchfield if and only if the $A$--module morphism
$1\otimes d\co A\otimes_{A[F_\mu]}K \to A\otimes_{A[F_\mu]}L$
is an isomorphism.

{\rm (ii)}\qua A morphism $d\co K \to L$ of projective $A[F_\mu]$--modules is
injective and $M = \coker(d)$ is a Blanchfield $A[F_\mu]$--module
if and only if the
$A$--module morphism $1\otimes d\co A\otimes_{A[F_\mu]}K \to
A\otimes_{A[F_\mu]}L$ is an isomorphism.
\end{proposition}
\begin{proof}
(i)\qua It follows from \fullref{B = F} and
the commutative diagram with exact rows and columns
$$\xymatrix@C+10pt{ & 0 \ar[d] & 0 \ar[d] & & \\
\Tor_1^{A[F_\mu]}(A,K) = 0 \ar[r] & \bigoplus\limits_{i = 1}^{\mu}
K  \ar[r]^-{\di{\gamma_K}} \ar[d]^-{\di{\oplus d}}& K \ar[r]\ar[d]^-{\di{d}} &
A\otimes_{A[F_\mu]}K \ar[r]\ar[d]^-{\di{1\otimes d}} & 0  \\
\Tor_1^{A[F_\mu]}(A,L) = 0 \ar[r] & \bigoplus\limits_{i = 1}^{\mu} L \ar[r]^-{\di{\gamma_L}} \ar[d]
& L \ar[r]\ar[d] & A\otimes_{A[F_\mu]}L \ar[r] & 0 \\
&  \bigoplus\limits_{i = 1}^{\mu} M \ar[r]^-{\di{\gamma_M}} \ar[d]& M \ar[d]& &\\
 & 0  & 0 & & }$$
that $M$ is Blanchfield if and only if $1\otimes d$ is an isomorphism.

(ii)\qua If $d$ is injective and $M$ is Blanchfield
then $1\otimes d$ is an isomorphism by (ii), since projective $A[F_\mu]$--modules
are flat. Conversely, if $1\otimes d\co A\otimes_{A[F_\mu]}K \to
A \otimes_{A[F_\mu]}L$ is an isomorphism then
$1\otimes d\co \smash{\widehat{A[F_\mu]}}\otimes_{A[F_\mu]}K \to
\widehat{A[F_\mu]} \otimes_{A[F_\mu]}L$ is an isomorphism
by \fullref{Magnus-Foxon} (iii), and it follows from the injectivity of
$K \to \smash{\widehat{A[F_\mu]}}\otimes_{A[F_\mu]}K$,
$L \to \smash{\widehat{A[F_\mu]}}\otimes_{A[F_\mu]}L$
and the commutative diagram
$$\xymatrix{K \ar[r]^-{\di{d}} \ar[d] & L \ar[d]\\
\widehat{A[F_\mu]}\otimes_{A[F_\mu]}K \ar[r]^-{\di{1\otimes d}}&
\widehat{A[F_\mu]}\otimes_{A[F_\mu]}L}$$
that $d\co K \to L$ is injective.
\end{proof}

The \emph{idempotent completion} $\PP(\E)$ of an additive category
$\E$ is the additive category with objects pairs $(M,p = p^2\co M \to M)$
defined by projections $p$ of objects $M$ in $\E$, and morphisms
$f\co (M,p) \to (N,q)$ defined by morphisms $f\co M \to N$ in $\E$ such that
$qfp = f\co M \to N$. As usual, $\E$ is \emph{idempotent complete} if
the functor $\E \to \PP(\E);M \mapsto (M,1)$ is an equivalence,
or equivalently if for every idempotent $p = p^2\co M \to M$ in $\E$
there exists a direct sum decomposition $M = P\oplus Q$ with
$$p = \begin{pmatrix} 1 & 0 \\ 0 & 0 \end{pmatrix}~:~
M = P\oplus Q \to M = P\oplus Q.$$
For any exact category $\E$ there exists a full embedding
$\E \subset \A$ in an abelian category $\A$ (Gabriel--Quillen),
and the idempotent completion
$\PP(\E)$ is equivalent to the full exact subcategory of $\A$
with objects ${\rm im}(p)$ for objects $(M,p)$ in $\PP(\E)$.
For $\E = \Flk(A) \subset \A = \Bla_{\infty}(A)$ we have that
$\PP(\Flk(A)) \subset \Bla_{\infty}(A)$. In fact, we have:

\begin{proposition} \label{idem}
{\rm (i)}\qua The exact categories $\Prim(A)$, $\Sei(A)$, $\Bla(A)$ are
idempotent complete.

{\rm (ii)}\qua The idempotent completion of $\Flk(A)$ is equivalent to $\Bla(A)$
$$\PP(\Flk(A))~\approx~\Bla(A).$$
\end{proposition}
\begin{proof} (i)\qua The exact categories $\Prim(A)$, $\Sei(A)$, $\Bla(A)$ are
closed under direct summands.

(ii)\qua For any f.g.~projective
$A[F_\mu]$--modules $K,L$ the augmentation map
$$\epsilon\co \Hom_{A[F_\mu]}(K,L) \to \Hom_A(A\otimes_{A[F_\mu]}K,
A\otimes_{A[F_\mu]}L);~d \mapsto 1\otimes d$$
is surjective, by the following argument:
choose f.g.~projective $A[F_\mu]$--modules $K',L'$ such that
$$K\oplus K' = A[F_\mu]^k,~L\oplus L' = A[F_\mu]^{\ell}$$
for some $k,\ell \geqslant 0$, and note that the augmentation map
\begin{multline*}
\epsilon\co \Hom_{A[F_\mu]}\bigl(K\oplus K',L\oplus L'\bigr) = 
\Hom_{A[F_\mu]}\bigl(A[F_\mu]^k,A[F_\mu]^{\ell}\bigr)\\
\longrightarrow \Hom_A\bigl(A\otimes_{A[F_\mu]}\bigl(K\oplus K'\bigr),
A\otimes_{A[F_\mu]}\bigl(L\oplus L'\bigr)\bigr) = \Hom_A\bigl(A^k,A^{\ell}\bigr)
\end{multline*}
is surjective. Given an h.d.~1 Blanchfield $A[F_\mu]$--module $M$ with
a f.g.~projective $A[F_\mu]$--module resolution
$$\xymatrix{0 \ar[r] & K \ar[r]^-{\di{d}} & L \ar[r]& M \ar[r]&0}$$
we know from \fullref{char} (i) that
$1\otimes d\co A\otimes_{A[F_\mu]}K \to A\otimes_{A[F_\mu]}L$ is an
$A$--module isomorphism. By \fullref{Magnus-Foxon} (i)
it is possible to lift $(1\otimes d)^{-1}$ to an $A[F_\mu]$--module
morphism $e\co L \to K$, so that by \fullref{char} (i) $e$ is an injection with
$$N = \coker(e)$$
an h.d.~1 Blanchfield $A[F_\mu]$--module. Let $J$ be a f.g.~projective
$A[F_\mu]$--module such that $J\oplus K \oplus L$ is f.g.~free,
say $A[F_\mu]^m$. The $A[F_\mu]$--module morphism
$$f = 1\oplus \begin{pmatrix}
0 & e \\ d & 0 \end{pmatrix}~:~
J \oplus K \oplus L = A[F_\mu]^m \to J \oplus K \oplus L = A[F_\mu]^m$$
is such that $1\otimes f\co A^m \to A^m$ is an isomorphism, so that
$\coker(f) = M\oplus N$ is an h.d.~1 $F_\mu$--link module.
The functor
$$\Flk(A) \to \Bla(A);~M \mapsto M$$
is a full embedding such that every object in $\Bla(A)$ is a direct summand
of an object in $\Flk(A)$, so that $\Bla(A)$ is (equivalent to) the
idempotent completion $\PP(\Flk(A))$.
\end{proof}

\subsection{Seifert modules}

Let $Q_\mu$ be the complete quiver which has $\mu$ vertices and
$\mu^2$ arrows, one arrow between each ordered pair of vertices.
The path ring is given by
$$Q_{\mu} = 
\Z[e]*\Z\bigl[\pi_1,\pi_2,\ldots,\pi_{\mu}\,\vert\,\pi_i\pi_j = \delta_{ij}\pi_i,
    {\textstyle\sum^{\mu}_{i = 1}}\pi_i = 1\bigr]$$
where $\pi_i e \pi_j$ corresponds to the unique path of length 1
from the $i$th vertex to the $j$th vertex. An $A$--module $P$
together with a ring morphism $\rho\co Q_\mu \to \End_A(P)$ is
essentially the same as a triple $(P,e,\{\pi_i\})$ with $e\co P \to
P$ an endomorphism, and $\{\pi_i\co P \to P\}$ a complete system of
$\mu$ idempotents. (Such representations of $Q_\mu$ were first
considered by Farber \cite{Fa3} for particular $A$.)

\begin{definition} \label{Sdefinition}
(i)\qua A \emph{Seifert $A$--module} $(P,e,\{\pi_i\})$ is an
$A$--module $P$ together with an endomorphism $e\co P \to P$,
and a system $\{\pi_i\co P \to P\}$ of idempotents expressing $P$
as a $\mu$--fold direct sum, with
$$\pi_i\co P = P_1 \oplus P_2 \oplus \cdots \oplus P_{\mu} \to P;~
(x_1,x_2,\ldots,x_\mu) \mapsto (0,\ldots,0,x_i,0,\ldots,0).$$
(ii)\qua A \emph{morphism} of Seifert $A$--modules
$$g\co (P,e,\{\pi_i\}) \to (P',e',\{\pi'_i\})$$
is an $A$--module morphism such that
$$ge = e'g,~g\pi_i = \pi'_ig\co P \to P'.$$
The conditions $g\pi_i = \pi'_ig$ are equivalent to $g$ preserving
the direct sum decompositions, so that
$$g = \begin{pmatrix} g_1 & 0 & \ldots & 0 \\
0 & g_2 & \ldots & 0 \\
0 & 0 & \ldots & 0 \\
\vdots & \vdots & \ddots & 0 \\
0 & 0 & \ldots & g_{\mu} \end{pmatrix}\co P = P_1 \oplus P_2 \oplus
\cdots \oplus P_{\mu} \to P' = P'_1 \oplus P'_2 \oplus \cdots \oplus
P'_{\mu}$$
with $g_i\co P_i \to P'_i$.

(iii)\qua The \emph{Seifert $A$--module category} $\Sei_{\infty}(A)$
has objects Seifert $A$--modules and morphisms as in (ii). Let $\Sei(A)
\subseteq \Sei_{\infty}(A)$ be the full subcategory of the Seifert
$A$--modules $(P,e,\{\pi_i\})$ with $P$ f.g.~projective.
\end{definition}

\subsection{The covering functor $B$}

Seifert modules determine $F_\mu$--link modules by:

\begin{definition} (i)\qua The \emph{covering} of a Seifert $A$--module
$(P,e,\{\pi_i\})$ is the $F_\mu$--link module
$$B(P,e,\{\pi_i\}) = \coker(
1-e+ez\co P[F_\mu] \to P[F_\mu])$$
with Mayer--Vietoris presentation
$$\xymatrix{0 \ar[r] &\bigoplus\limits_{i = 1}^\mu P_i[F_{\mu}]
\ar[r]^-{\di{d}} & P[F_{\mu}] \ar[r] & B(P,e,\{\pi_i\})\ar[r] & 0 ,}$$
where $d = 1-e+ez$.\\
(ii)\qua The \emph{covering} of a Seifert $A$--module
morphism $g\co (P,e,\{\pi_i\}) \to (P',e',\{\pi'_i\})$ is the
$F_\mu$--link module morphism
$$B(g)\co B(P,e,\{\pi_i\}) \to B(P',e',\{\pi'_i\});~x \mapsto g(x)$$
resolved by
$$\xymatrix{0 \ar[r] &P[F_\mu] \ar[r]^-{\di{d}}
\ar[d]^-{\di{g}} & P[F_\mu] \ar[d]^-{\di{g}} \ar[r] &
B(P,e,\{\pi_i\}) \ar[d]^-{\di{B(g)}} \ar[r] &0\\
0 \ar[r] & P'[F_\mu] \ar[r]^-{\di{d'}} & P'[F_\mu] \ar[r] &
B(P',e',\{\pi'_i\})\ar[r] & 0}$$
with $d = 1-e+ez$, $d' = 1-e'+e'z$.
\end{definition}

\begin{example} \label{chain}
Let $\ell\co \bigsqcup_\mu S^n \subset S^{n+2}$ be a
$\mu$--component boundary link with exterior $W$, so that there exists a
$\Z$--homology equivalence $c\co W \to W_0$ to the exterior $W_0$ of the trivial
$\mu$--component boundary link $\ell_0\co \bigsqcup_\mu S^n \subset
S^{n+2}$. The $(n{+}2)$--dimensional f.g.~free $\Z[F_\mu]$--module
chain complex
$$\dot C(\wwtilde{W}) = {\cal C}(\wtilde{c}\co C(\wwtilde{W}) \to
C(\wwtilde{W}_0))_{*+1}$$
is $\Z$--contractible.
For any $\mu$--component Seifert surface $V = V_1\sqcup V_2 \sqcup \ldots
\sqcup V_\mu \subset S^{n+2}$ for $\ell$ there exists a degree 1 map $V
\to V_0$ to the $\mu$--component Seifert surface
$V_0 = \bigsqcup_{\mu}D^{n+1} \subset S^{n+2}$ for $\ell_0$.
Let
$$\dot C(V_i) = {\mathcal C}(C(V_i) \to C(D^{n+1}))_{*+1},\quad
\dot C(V) = {\sum\limits^\mu_{i = 1}}\dot C(V_i).$$
The map $V \to S^{n+2}\backslash V$ pushing $V$ off itself in
the positive normal direction combines with chain level
Alexander duality to induce a $\Z$--module chain map
$$e\co \dot C(V) \to C\bigl(S^{n+2}\backslash V,
    {\textstyle\bigsqcup_\mu} \{{\rm pt.}\}\bigr)~\simeq~ \dot C(V)^{n+1-*},$$
so that there is defined an $(n{+}1)$--dimensional chain complex
$(\dot C(V),e,\{\pi_i\})$ in $\Sei(\Z)$.
The covering $B(\dot C(V),e,\{\pi_i\})$ is an
$(n{+}1)$--dimensional chain complex in $\Flk(\Z)$,
with the projection
\begin{multline*}
{\mathcal C}(1-e+ez\co \dot C(V)[F_\mu] \to \dot C(V)[F_\mu]) = 
\dot C(\wwtilde{W})\\
\longrightarrow B(\dot C(V),e,\{\pi_i\}) = 
\coker(1-e+ez\co \dot C(V)[F_\mu] \to \dot C(V)[F_\mu])
\end{multline*}
a homology equivalence.
\end{example}

The covering construction defines a functor of exact categories
$$B_{\infty}\co \Sei_{\infty}(A) \to \Bla_{\infty}(A);~
(P,e,\{\pi_i\}) \mapsto B(P,e,\{\pi_i\})$$
which restricts to a functor $B\co \Sei(A) \to \Flk(A)$.

\begin{definition} A morphism $f$ in $\Sei_{\infty}(A)$ is a
\emph{$B$--isomorphism} if $B(f)$ is an isomorphism in
$\Bla_{\infty}(A)$. Let $\Xi_{\infty}$ denote the set of
$B$--isomorphisms in $\Sei_{\infty}(A)$, and let $\Xi$ denote the
set of $B$--isomorphisms in $\Sei(A)$.
\end{definition}

\subsection{Blanchfield/Seifert algebraic transversality}

We shall now use the algebraic transversality of \fullref{algebraic
transversality} to establish that every h.d.~1 $F_\mu$--link module $M$ is
isomorphic to the covering $B(P,e,\{\pi_i\})$ of a f.g.~projective
Seifert $A$--module $(P,e,\{\pi_i\})$, uniquely up to morphisms in
$\Xi$.

We refer to Sheiham \cite{Sh2} for the proof that
$B_\infty\co \Sei_\infty(A) \to \Flk_\infty(A)$ induces an
equivalence of exact categories
$\wbar{B}_\infty\co \Xi_\infty^{-1}\Sei_\infty(A)
    \approx \Flk_\infty(A)$.
Algebraic transversality will be used to prove
that the universal localization
$\Sei(A) \to \Xi^{-1}\Sei(A)$ has a calculus of
fractions, and that the covering functor $B\co \Sei(A) \to
\Flk(A)$ induces an equivalence of exact categories
$\wbar{B}\co \Xi^{-1}\Sei(A) \approx \Flk(A)$.

Given an $F_\mu$--link module $M$ let $U(M) = (M,e_M,\{\pi_i\})$ be
the Seifert $A$--module defined in \cite{Sh2} -- the definition is
recalled in the Introduction of this paper, along with the fact
proved in \cite{Sh2} that $B_{\infty}$ is a left adjoint of
$$U_{\infty}\co \Bla_{\infty}(A) \to \Sei_{\infty}(A);~
    M \mapsto U(M).$$
The natural isomorphism of the adjointness
\begin{eqnarray*}
\Hom_{\Bla_{\infty}(A)}(B(Q,f,\{\rho_i\}),M)
&\xymatrix{\ar[r]^-{\di{\cong}}&} &
\Hom_{\Sei_{\infty}(A)}((Q,f,\{\rho_i\}),U(M));\\
g &\longmapsto& \adj(g) = U(g) h
\end{eqnarray*}
is defined for any Seifert $A$--module $(Q,f,\{\rho_i\})$,
with
$$h\co Q \subset Q[F_\mu] \to UB(Q,f,\{\rho_i\})$$
the restriction of the canonical surjection $Q[F_\mu] \to B(Q,f,\{\rho_i\})$.
If $M$ is h.d.~1 and $(Q,f,\{\rho_i\})$ is f.g.~projective
the natural isomorphism can be written as
$$\Hom_{\Flk(A)}(B(Q,f,\{\rho_i\}),M)\cong~
\Hom_{\Sei_{\infty}(A)}((Q,f,\{\rho_i\}),U(M))$$
but note that in general $U(M)$ is not a f.g.~projective Seifert $A$--module.

The following result establishes that for an h.d.~1 $F_\mu$--link
module $M$ the Seifert $A$--module $U(M)$
is the direct limit of a directed system of f.g.~projective Seifert
$A$--modules $(P,e,\{\pi_i\})$ and morphisms in $\Xi$,
with isomorphisms $B(P,e,\{\pi_i\})\cong M$.

\begin{thm}[Blanchfield/Seifert algebraic transversality]
\label{Btrans}
Let $M$ be an h.d.~1 $F_\mu$--link module,
with a 1--dimensional induced f.g.~projective $A[F_\mu]$--module resolution
$$\xymatrix{0 \ar[r] & P[F_\mu] \ar[r]^-{\di{d}} &
P[F_\mu]\ar[r] & M \ar[r] &0}$$
such that $\epsilon(d)\co P \to P$ is an $A$--module isomorphism.

%\begin{itemize}
%\item[\rm (i)]
{\rm(i)}\qua
Let $I_\infty$ be the set of ordered pairs
$T = (T_0,T_1)$ of subtrees $T_0,T_1 \subseteq G_\mu$ such that
$d_*(T_1) \subseteq T_0$.
The set $I_\infty$ is partially ordered by inclusion, with maximal element
$$T_{\max} = \bigcup\limits_{T \in I_\infty}T = (G_\mu,G_\mu) \in I_{\infty}.$$
There is defined a directed system of Seifert $A$--modules
$(P\langle T \rangle,e\langle T\rangle,
\{\pi_i\langle T\rangle\})$ and morphisms in $\Xi_{\infty}$
$$\phi\langle T,T' \rangle\co 
(P\langle T\rangle,e\langle T\rangle,\{\pi_i\langle T\rangle\})
\longrightarrow
  (P\langle T'\rangle,e\langle T'\rangle,\{\pi_i\langle T'\rangle\})
\quad(T \subseteq T' \in I_\infty)$$
with direct limit
$$\varinjlim\limits_{T \in I_\infty}
(P\langle T\rangle,e\langle T\rangle,\{\pi_i\langle T\rangle\}) = 
(P\langle T_{\max}\rangle,e\langle T_{\max}\rangle,\{\pi_i\langle T_{\max}\rangle\}) = 
U(M).$$
For any $T  = (T_0,T_1)\in I_\infty$ the morphism
$\phi\langle T,T_{\max} \rangle\co 
(P\langle T\rangle,e\langle T\rangle,\{\pi_i\langle T\rangle\})\to U(M)$
is the adjoint $\phi\langle T,T_{\max} \rangle  = \adj(\phi\langle T \rangle)$
of an isomorphism in $\Flk(A)$
$$\phi\langle T \rangle\co 
B(P\langle T\rangle,e\langle T\rangle,\{\pi_i\langle T\rangle\})
\xymatrix{\ar[r]^-{\di{\cong}}&} M$$
such that for any $T \subseteq T' \in I_\infty$ there is defined a commutative
triangle of isomorphisms in $\Flk_\infty(A)$
$$\xymatrix@R+20pt{
B(P\langle T\rangle,e\langle T\rangle,\{\pi_i\langle T\rangle\})
\ar[dr]^-{\di{\phi\langle T \rangle}}_{\di{\cong}}
\ar[rr]^-{\di{B(\phi\langle T,T' \rangle)}}_{\di{\cong}}
& &B(P\langle T'\rangle,e\langle T'\rangle,\{\pi_i\langle T'\rangle\})
\ar[dl]_-{\di{\phi\langle T' \rangle}}^-{\di{\cong}}\\
&M&}$$
In particular, $\phi\langle T,T_{\max} \rangle \in \Xi_{\infty}$.

%\item[\rm (ii)]
{\rm(ii)}\qua
Let $I\subset I_\infty$ be the subset of the ordered pairs
$T = (T_0,T_1)$ of finite subtrees $T_0,T_1 \subset G_\mu$ such that
$d_*T_1 \subseteq T_0$. For $T \in I$
$(P\langle T\rangle,e\langle T\rangle,\{\pi_i\langle T\rangle\})$ is
a f.g.~projective Seifert $A$--module, and
$$\varinjlim\limits_{T \in I}
(P\langle T\rangle,e\langle T\rangle,\{\pi_i\langle T\rangle\}) = 
U(M)$$
with $\phi\langle T,T' \rangle \in \Xi$ $(T \subseteq T' \in I)$.

%\item[\rm (iii)]
{\rm(iii)}\qua
For any f.g.~projective Seifert $A$--module $(Q,f,\{\rho_i\})$ every
morphism
$$g\co B(Q,f,\{\rho_i\}) \to M$$
in $\Flk(A)$ factors as
$$g\co B(Q,f,\{\rho_i\}) \xymatrix@C+20pt{\ar[r]^-{\di{B(g\langle T \rangle)}}&}
B(P\langle T \rangle,e\langle T \rangle,\{\pi_i\langle T \rangle\})
\xymatrix@C+10pt{\ar[r]_-{\di{\cong}}^-{\di{\phi\langle T \rangle}}&} M$$
for some $T \in I$, with
$g\langle T \rangle\co (Q,f,\{\rho_i\}) \to
(P\langle T\rangle,e\langle T\rangle,\{\pi_i\langle T\rangle\})$
a morphism in $\Sei(A)$.
%\end{itemize}
\end{thm}
\begin{proof} (i)\qua The induced f.g.~projective
$A[F_\mu]$--module chain complex
$$\xymatrix{E\co E_1 = P[F_\mu] \ar[r]^-{{d}} & E_0 = P[F_\mu]}$$
is such that $H_0(E) = M$, $H_1(E) = 0$. By \fullref{MV} for any subtree
$T_1 \subseteq G_\mu$ there exists a subtree $d_*(T_1) \subseteq G_\mu$
such that for any subtree $T_0 \subseteq G_\mu$ with $d_*(T_1) \subseteq T_0$
$E$ admits a Mayer--Vietoris presentation
$$\xymatrix@C+20pt{
E_1\langle T_1 \rangle\co 0 \ar[r] & \bigoplus\limits_{i = 1}^\mu
C_1^{(i)} [F_{\mu}]  \ar[r]^-{{f_1^+z-f_1^-}} \ar[d]^-{{d_C}}
& D_1[F_{\mu}] \ar[r]\ar[d]^-{{d_D}} & E_1 \ar[r]\ar[d]^-{{d}} & 0  \\
E_0\langle T_0 \rangle\co 0 \ar[r] & \bigoplus\limits_{i = 1}^\mu
C_0^{(i)} [F_{\mu}] \ar[r]^-{{f_0^+z-f_0^-}} & D_0[F_{\mu}] \ar[r] & E_0
\ar[r] & 0 }$$
\begin{align*}
&\text{with}\quad C^{(i)}_j = P\bigl[T_j^{(i,1)}\bigr],\quad
D_j = P\bigl[T_j^{(0)}\bigr]~\subseteq~ E_j = P[F_\mu]\quad(j = 0,1),\\
&\text{and}\quad d_C = {\textstyle\bigoplus_{i = 1}^\mu} d\vert\co 
{\textstyle\bigoplus_{i = 1}^\mu} C^{(i)}_1 \to
{\textstyle\bigoplus_{i = 1}^\mu} C^{(i)}_0,\quad d_D = d\vert\co D_1 \to D_0.
\end{align*}
The $A$--modules defined by
\begin{align*}
P_i\langle T \rangle& = 
\coker(d\vert\co C^{(i)}_1 \to  C^{(i)}_0),\\
P\langle T \rangle& = \coker(d_C) = 
{\textstyle\bigoplus_{i = 1}^\mu} P_i\langle T \rangle,\\
Q\langle T \rangle& = \coker(d_D)
\end{align*}
fit into a commutative diagram of $A[F_\mu]$--modules with exact rows
and columns
$$\xymatrix@C+10pt@R-5pt{ & 0 \ar[d] & 0 \ar[d] & 0 \ar[d] & \\
0 \ar[r] & \bigoplus\limits_{i = 1}^\mu
C_1^{(i)} [F_{\mu}]  \ar[r]^-{{f_1^+z-f_1^-}} \ar[d]^-{{d_C}}
& D_1[F_{\mu}] \ar[r]\ar[d]^-{{d_D}} & P[F_\mu] \ar[r]\ar[d]^-{{d}} & 0  \\
0 \ar[r] & \bigoplus\limits_{i = 1}^\mu
C_0^{(i)} [F_{\mu}] \ar[r]^-{{f_0^+z-f_0^-}} \ar[d]
& D_0[F_{\mu}] \ar[r]\ar[d] & P[F_\mu] \ar[r] \ar[d]& 0 \\
0 \ar[r] & P\langle T \rangle [F_\mu] \ar[r]^-{{f^+z-f^-}} \ar[d]& Q\langle T \rangle
[F_\mu] \ar[r] \ar[d]& M \ar[r] \ar[d]& 0\\
 & 0  & 0 & 0 & }$$
with
$f^+,f^-\co P\langle T \rangle \to Q\langle T \rangle$
the $A$--module morphisms induced by
$$f^+_0,f^-_0\co {\textstyle\bigoplus^\mu_{i = 1}}
C^{(i)}_0 \to D_0.$$
It follows from $\Tor^{A[F_\mu]}_1(A,M) = 0$ that
$f^+-f^-\co P\langle T \rangle \to Q\langle T \rangle$
is an $A$--module isomorphism. The Seifert $A$--module
$(P\langle T \rangle,e\langle T \rangle,\{\pi_i\langle T \rangle\})$
defined by
$$
e\langle T \rangle = (f^+-f^-)^{-1}f^+\co P\langle T \rangle \to
P\langle T \rangle,\quad
\pi_i\langle T \rangle\co P \langle T \rangle\to
P_i\langle T \rangle \to P\langle T \rangle
$$
is such that $P\langle T \rangle[F_\mu]\cong Q\langle T \rangle[F_\mu] \to M$
induces the isomorphism of Blanchfield $A[F_\mu]$--modules
$\phi\langle T \rangle\co 
B(P\langle T \rangle,e\langle T \rangle,\{\pi_i\langle T \rangle\})
\cong M$ adjoint to the natural map
$(P\langle T \rangle,e\langle T \rangle,\{\pi_i\langle T \rangle\})
\to U(M)$. (In particular,
$(P\langle T_{\max} \rangle,e\langle T_{\max} \rangle,\{\pi_i\langle T_{\max} \rangle\})
\allowbreak  = U(M)$
and $\phi\langle T_{\max} \rangle\co BU(M)\cong M$ is the natural
isomorphism $\psi_M$ defined in \cite[5.10]{Sh2}.)
For $T \subseteq T' \in I$ the $B$--isomorphism $\phi \langle T,T' \rangle$
is induced by the inclusion $T \subseteq T'$.

(ii)\qua The augmentation of the $A[F_\mu]$--module morphism $d\co P[F_\mu] \to P[F_\mu]$
is an $A$--module isomorphism $\epsilon(d)\co P \to P$, so that the
induced $\smash{\widehat{A[F_\mu]}}$--module morphism $\widehat{d}\co 
\smash{\widehat{P[F_\mu]}} \to \smash{\widehat{P[F_\mu]}}$ is an isomorphism, by
\fullref{Magnus-Foxon}. For any $T = (T_0,T_1) \in I$ the inclusion
$P\bigl[\smash{T_1^{(0)}}\bigr] \to \smash{\widehat{P[F_\mu]}}$ is a split $A$--module injection
by \fullref{Magnus-Foxtw}. Let $s\co \smash{\widehat{P[F_\mu]}} \to
P\bigl[\smash{T_1^{(0)}}\bigr]$
be a splitting $A$--module surjection. The anticlockwise composition
of the morphisms (inverting $\widehat{d}$) in the diagram
$$\xymatrix{
D_1 = P\bigl[T_1^{(0)}\bigr] \ar[d]_-{\di{d_D = d|}}\ar[r] & P[F_\mu]
\ar[r]\ar[d]^-{\di{d}} &
\widehat{P[F_\mu]}\ar@/_2pc/@{>}[ll]_-{\di{s}}
\ar[d]^-{\di{\widehat{d}}}_-{\di{\cong}} \\
D_0 = P\bigl[T_0^{(0)}\bigr]  \ar[r] & P[F_\mu] \ar[r] & \widehat{P[F_\mu]}
}$$
defines an $A$--module surjection $P\bigl[\smash{T_0^{(0)}}\bigr] \to
P\bigl[\smash{T_1^{(0)}}\bigr]$
splitting $d\vert\co P\bigl[\smash{T_1^{(0)}}\bigr] \to
P\bigl[\smash{T_0^{(0)}}\bigr]$. Thus
$d\vert$ is a split injection of f.g.~projective $A$--modules and
$P\langle T \rangle = \coker(d\vert)$ is a f.g.~projective
$A$--module.

(iii)\qua The morphism $g\co B(Q,f,\{\rho_i\}) \to M$ in $\Flk(A)$
has a canonical resolution
$$\xymatrix@C+10pt@R+10pt{
0 \ar[r] & Q[F_{\mu}] \ar[d]^-{\di{\adj(g)}}\ar[r]^-{\di{1-f+fz}}
& Q[F_{\mu}] \ar[d]^-{\di{\adj(g)}}\ar[r] &
B(Q,f,\{\rho_i\}) \ar[d]^-{\di{g}}\ar[r] \ar[d]& 0 \\
0 \ar[r] & P\langle T_{\max} \rangle  [F_\mu]
\ar[r]& P\langle T_{\max} \rangle [F_\mu]\ar[r] & M \ar[r]& 0}$$
with
$$\adj(g)\co (Q,f,\{\rho_i\}) \to
(P\langle T_{\max} \rangle,e\langle T_{\max} \rangle,
\{\pi_i\langle T_{\max} \rangle\}) = U(M)$$
the adjoint morphism in $\Sei_{\infty}(A)$.
Since $Q$ is f.g.~projective there exists $T \in I$ such that
$${\rm im}(g\co B(Q,f,\{\rho_i\}) \to M) \subseteq {\rm im}(
B(P\langle T \rangle,e\langle T \rangle,\{\pi_i\langle T \rangle\})
\to M)$$
with a lift of $g$ to an $A$--module morphism
$g\langle T \rangle\co Q \to P\langle T \rangle$
which preserves the direct sum structures.
The diagram of $A$--modules and morphisms
$$\xymatrix@C+10pt{
Q \ar[rr]^-{\di{f}} \ar[dd]^-{\di{g\langle T \rangle}}
\ar@/_4pc/[dddd]^{\di{g}}
&&
Q \ar[dd]^-{\di{g\langle T \rangle}}\ar@/^4pc/[dddd]^{\di{g}}
\\
&(*)& \\
P\langle T \rangle
\ar[rr]^-{\di{e\langle T \rangle}}
\ar[dd]^-{\di{\phi\langle T \rangle}} &&
P\langle T \rangle
\ar[dd]^-{\di{\phi\langle T \rangle}} \\
&&\\
U(M)\ar[rr]^-{\di{e}}&&U(M)}$$
commutes except possibly in $(*)$, and $(*)$ commutes if and only if
$$g\langle T \rangle\co (Q,f,\{\rho_i\}) \to (P\langle T \rangle,
e\langle T \rangle,\{\pi_i\langle T \rangle\})$$
is a morphism of Seifert $A$--modules.  Since $Q$ is f.g.~projective and
the composite
$$\xymatrix@C+60pt{
Q \ar[r]^-{{g\langle T \rangle f-e\langle T \rangle g\langle T \rangle}} &
} P\langle T \rangle \xymatrix@C+15pt{\ar[r]^-{{\phi\langle T \rangle}} &}
U(M) = \varinjlim\limits_{T'\in I} P\langle T' \rangle$$
is 0 there exists $T' \in I$ such that $ T \subseteq T'$ and the composite
$$\xymatrix@C+60pt{g\langle T' \rangle f-e\langle T' \rangle g\langle T' \rangle\co 
Q \ar[r]^-{{g\langle T \rangle f-e\langle T \rangle g\langle T \rangle}} &
} P\langle T \rangle \xymatrix{\ar[r]&}  P\langle T' \rangle$$
is 0, so that
$$g\langle T' \rangle\co (Q,f,\{\rho_i\}) \to (P\langle T' \rangle,
e\langle T' \rangle,\{\pi_i\langle T' \rangle\})$$
is a morphism of Seifert $A$--modules as required (except that $T'$ has
to be called $T$).
\end{proof}

\begin{definition}
Let $M = B(P,e,\{\pi_i\})$ for a f.g.~projective Seifert $A$--module
$(P,e,\{\pi_i\})$.

(i)\qua For any $T \in I_{\infty}$ let
$$s\langle T \rangle\co (P,e,\{\pi_i\}) \to
(P\langle T \rangle ,e\langle T \rangle,\{\pi_i\langle T \rangle\})$$
be the $B$--isomorphism determined by the inclusion
$P = P[\{1\}] \subseteq P\bigl[\smash{T^{(0)}_0}\bigr]$.

(ii)\qua For $T = T_{\max}\in I_{\infty}$ write
$$s_M = s\langle T_{\max} \rangle\co 
(P,e,\{\pi_i\}) \to
(P\langle T_{\max} \rangle ,e\langle T_{\max} \rangle,\{\pi_i\langle T_{\max} \rangle\}) = 
U(M).$$
This is the $B$--isomorphism adjoint of $1\co M \to M$, such that
$$\xymatrix@C+10pt{s_M\co (P,e,\{\pi_i\}) \ar[r]^-{{s\langle T \rangle}}&
(P\langle T \rangle ,e,\{\pi_i\})
\ar[r]^-{{\phi \langle T \rangle}} & U(M)}$$
for any $T \in I_{\infty}$.
\end{definition}

Putting everything together:

\begin{thm} \label{calculus}
{\rm (i)}\qua Every h.d.~1 $F_\mu$--link module $M$ is isomorphic
to the covering $B(P,e,\{\pi_i\})$ of a f.g.~projective Seifert $A$--module
$(P,e,\{\pi_i\})$.

{\rm (ii)}\qua For any f.g.~projective Seifert $A$--modules $(P,e,\{\pi_i\})$,
$(Q,f,\{\rho_i\})$ every morphism
$g\co B(Q,f,\{\rho_i\})\to B(P,e,\{\pi_i\})$
in $\Flk(A)$ is of the form $g = B(s)^{-1}B(t)$
for some morphisms
$$s\co (P,e,\{\pi_i\}) \to (P',e',\{\pi'_i\}),~
t\co (Q,f,\{\rho_i\}) \to (P',e',\{\pi'_i\})$$
in $\Sei(A)$ with $s\in \Xi$.

{\rm (iii)}\qua If $u\co (Q,f,\{\rho_i\})\to (P,e,\{\pi_i\})$
is a morphism of f.g.~projective Seifert $A$--modules such that
$B(u) = 0$ there exists an element
$v\co (P,e,\{\pi_i\})\to (P',e',\{\pi'_i\})$ in $\Xi$
such that $vu = 0$.

{\rm (iv)}\qua The localization $\Xi^{-1}\Sei(A)$ has a left calculus of fractions, and the
covering construction defines an equivalence of exact categories
$$\wbar{B}\co 
\Xi^{-1}\Sei(A)
\xymatrix{\ar[r]^-{\di{\approx}}&} \Flk(A);~(P,e,\{\pi_i\}) \mapsto B(P,e,\{\pi_i\}).$$
\end{thm}
\begin{proof} (i)\qua By \fullref{Btrans} (i)--(ii) $M$ is
isomorphic to $B(P\langle T \rangle,e\langle T \rangle,\{\pi_i\langle T \rangle\})$
for any $T \in I$, e.g. for the minimal element
$T_{\min} = (d_*\{1\},\{1\}) \in I$.

(ii)\qua By \fullref{Btrans} (iii) the adjoint of $g$ factors in
$\Sei_{\infty}(A)$ as
$$\xymatrix@R+10pt@C-5pt{
(Q,f,\{\rho_i\}) \ar[rr]^-{\di{\adj(g)}}
\ar[dr]_-{\di{g\langle T \rangle}}
& &UB(P,e,\{\pi_i\})\\
&(P\langle T \rangle,e\langle T \rangle,\{\pi_i\langle T \rangle\})
\ar[ur]_-{\di{\adj(\phi\langle T \rangle)}}&}$$
for some $T\in I$. The morphisms in $\Sei(A)$ defined by
$$\begin{array}{l}
s = s\langle T \rangle\co (P,e,\{\pi_i\}) \to
(P',e',\{\pi'_i\}) = 
(P\langle T \rangle,e\langle T \rangle,\{\pi_i\langle T \rangle\})\\[1ex]
t = g\langle T \rangle\co (Q,f,\{\rho_i\}) \to
(P',e',\{\pi'_i\}) = 
(P\langle T \rangle,e\langle T \rangle,\{\pi_i\langle T \rangle\})
\end{array}$$
are such that $s$ is a $B$--isomorphism (ie $s \in \Xi$)
and $g = B(s)^{-1}B(t)$.

(iii)\qua Let $M = B(P,e,\{\pi_i\})$. We have a commutative diagram in $\Sei_{\infty}(A)$
$$\xymatrix@R+10pt@C-5pt{
(Q,f,\{\rho_i\}) \ar[rr]^-{\di{\adj(B(u)) = 0}}
\ar[dr]_-{\di{u}} & &U(M)\\
&(P,e,\{\pi_i\})
\ar[ur]_-{\di{\theta_M}}&}$$
Since $Q$ is f.g.~projective there exists $T \in I$ such that
$$v = s\langle T \rangle\co (P,e,\{\pi_i\}) \to (P\langle T \rangle,e\langle T \rangle,\{\pi_i\langle T \rangle\})$$
is a $B$--isomorphism in $\Sei(A)$ (ie $v \in \Xi$) with $vu = 0$.

(iv)\qua Immediate from (i)--(iii).
\end{proof}

This completes the proof of \fullref{thm2} of the Introduction.

\section{Primitive Seifert modules} \label{kernel}

This section is devoted to the kernel of the covering functor
$B\co \Sei(A) \to \Flk(A)$.
Following the terminology of Sheiham \cite{Sh2}:

\begin{definition}
(i)\qua A Seifert $A$--module $(P,e,\{\pi_i\})$ is {\it
primitive} if
$$B(P,e,\{\pi_i\}) = 0$$
or equivalently $1-e+ez\co P[F_\mu] \to P[F_\mu]$ is an $A[F_\mu]$--module
isomorphism.

(ii)\qua Let $\Prim(A) \subset \Sei(A)$ be the full subcategory with
objects the primitive f.g.~projective Seifert $A$--modules.
\end{definition}

We shall now obtain an intrinsic characterization of
the objects in $\Prim(A)$, generalizing the results for $\mu = 1$
recalled below.

\begin{definition}[{{L\"uck and Ranicki \cite[Section~5]{LR}}}]
A \emph{near-projection} $(P,e)$ is an $A$--module $P$
together with an endomorphism $e \in \End_A(P)$ such that
$e(1-e) \in \End_A(P)$ is nilpotent.
\end{definition}

\begin{proposition}
[Bass, Heller and Swan \cite{BHS}, L\"uck and Ranicki \cite{LR}]\quad
\label{split}

{\rm (i)}\qua A linear morphism of induced f.g.~projective $A[z]$--modules
$$f_0+f_1z\co P[z] \to Q[z]$$
is an isomorphism if and only if $f_0+f_1\co P \to Q$ is an isomorphism and
$$e = (f_0+f_1)^{-1}f_1\co P \to P$$
is nilpotent.

{\rm (ii)}\qua A linear morphism of induced f.g.~projective
$A[z,z^{-1}]$--modules
$$f_0+f_1z\co P[z,z^{-1}] \to Q[z,z^{-1}]$$
is an isomorphism if and only if $f_0+f_1\co P \to Q$ is an isomorphism and
$$e = (f_0+f_1)^{-1}f_1\co P \to P$$
is a near-projection.

{\rm (iii)}\qua Suppose that $(P,e)$ is a near-projection, or equivalently
that
$$1-e+ze\co P[z,z^{-1}] \to P[z,z^{-1}]$$
is an $A[z,z^{-1}]$--module automorphism. If $N \geqslant 0$ is so large
that $(e(1-e))^N = 0$ then
$$e^N+(1-e)^N\co P \to P$$
is an $A$--module automorphism, and the endomorphism
$$e_{\omega} = (e^N+(1-e)^N)^{-1}e^N\co P \to P$$
is a projection, with $e_{\omega}(1-e_{\omega}) = 0$.
The submodules of $P$
\begin{align*}
P^+& = (1-e_{\omega})(P) = (1-e)^N(P) = \{x \in P\,\vert\, (1-e+ez)^{-1}e(x)
\in P[z]\},\\[1ex]
P^-& = e_{\omega}(P) = e^N(P) = \{x \in P\,\vert\,
(1-e+ez)^{-1}(1-e)(x)\in z^{-1}P[z^{-1}]\}
\end{align*}
are such that
$$(P,e) = (P^+,e^+) \oplus (P^-,e^-)$$
with $e^+\co P^+ \to P^+$ and
$1-e^-\co P^- \to P^-$ nilpotent.
\end{proposition}

\begin{definition} {\rm
A f.g.~projective
Seifert $A$--module $(P,e,\{\pi_i\})$ is \emph{strongly nilpotent} if the
$A[F^+_{\mu}]$--module endomorphism
$$ez = \sum_{i = 1}^\mu e\pi_iz_i\co P[F_\mu^+]\to P[F_\mu^+]$$
is nilpotent, ie $(ez)^N = 0$ for some $N \geqslant 1$.}
\end{definition}

The condition for strong nilpotence is equivalent to the
$A[F_{\mu}]$--module endomorphism
$$ez = \sum_{i = 1}^\mu e\pi_iz_i\co P[F_\mu]\to P[F_\mu]$$
being nilpotent.

Expressed as a representation of the complete quiver $Q_{\mu}$, a
Seifert module $(P, \rho\co Q_{\mu}\to\End_AP)$ is strongly nilpotent
if and only if there exists  $N \geqslant 1$ such that $\rho(p) = 0$
for every path $p\in Q_{\mu}$ of length $\geqslant N$.

\begin{proposition} The following conditions on a f.g.~projective
Seifert $A$--module $(P,e,\{\pi_i\})$
are equivalent:
\begin{itemize}
\item[\rm (i)] $(P,e,\{\pi_i\})$ is strongly nilpotent,
\item[\rm (ii)] the $A[F_{\mu}^+]$--module endomorphism
$$1-ez\co P[F^+_{\mu}] \to P[F^+_{\mu}]$$
is an automorphism,
\item[\rm (iii)] the $A[F_{\mu}^+]$--module endomorphism
$$1-e+ez\co P[F^+_{\mu}] \to P[F^+_{\mu}]$$
is an automorphism.
\end{itemize}
\end{proposition}
\begin{proof}
(i) $\Longrightarrow$ (ii)\qua If $(ez)^N = 0$ then $1-ez$ has inverse
$$\begin{array}{l}
(1-ez)^{-1} = 1+ez + (ez)^2 +\cdots + (ez)^{N-1}\\[1ex]
\hskip100pt
\in \Hom_{A[F^+_{\mu}]}(P[F^+_{\mu}],P[F^+_{\mu}]) = 
\Hom_A(P,P)[F^+_{\mu}].
\end{array}$$
(ii) $\Longrightarrow$ (i)\qua The inverse of $1-ez$ is of the form
$$(1-ez)^{-1} = \sum_{\scriptsize
\begin{array}{c}
1 \leqslant i_1,i_2,\ldots,i_k \leqslant \mu\\
n_1,n_2,\ldots,n_k \geqslant 0\\
n_1+n_2+\cdots+n_k < N
\end{array}} f_{i_1 i_2 \ldots i_k} z_{i_1}^{n_1}z_{i_2}^{n_2}
\ldots z_{i_k}^{n_k}\co P[F^+_{\mu}] \to P[F^+_{\mu}]$$
for some $N \geqslant 1$. We have the identity
$$\begin{array}{l}
(1-ez)^{-1}-(1+ez+(ez)^2+\cdots+(ez)^{N-1}) = (1-ez)^{-1}(ez)^N\\[2ex]
\hskip100pt
\in \Hom_{A[F^+_{\mu}]}(P[F^+_{\mu}],P[F^+_{\mu}]) = \Hom_A(P,P)[F^+_{\mu}]
\end{array}$$
in which the left hand side is a sum of monomials in
$z_{i_1}z_{i_2}^{n_2}\ldots z_{i_k}^{n_k}$ of degree
$n_1+n_2+\cdots+n_k < N$ and the right hand side is a sum of monomials
of degree $\geqslant N$. Both sides of the identity are thus 0,
$$(ez)^N = 0\co P[F^+_{\mu}] \to P[F^+_{\mu}]$$
and $(P,e,\{\pi_i\})$ is strongly nilpotent.

(ii) $\Longleftrightarrow$ (iii)\qua Immediate from the identity
$$1-e+ez = 1-e(1-z)\co P[F^+_{\mu}] \to P[F^+_{\mu}]$$
and the change of variables $z_i \mapsto 1-z_i$.
\end{proof}

\begin{definition}
A $\mu$--component Seifert $A$--module $(P,e,\{\pi_i\})$ is a
\emph{near-projection} if it can be expressed as
$$(P,e,\{\pi_i\}) = \biggl(P^+\oplus P^-,\begin{pmatrix}
e^{++} & e^{+-} \\
e^{-+} & e^{--}
\end{pmatrix},\{\pi_i^+\}\oplus \{\pi_i^-\}\biggr)$$
and the $2\mu$--component Seifert $A$--module
$$(P',e',\pi') = \biggl( P^+\oplus P^-\ ,\
\begin{pmatrix}
e^{++} & -e^{+-} \\
e^{-+} & 1-e^{--}
\end{pmatrix}\ ,\ \{\pi_i^+\}\oplus\{\pi_i^-\}\biggr)$$
is strongly nilpotent.
\end{definition}

\begin{lemma} \label{near}
For a near-projection $(P,e,\{\pi_i\})$
the pairs $(P,e)$, $(P,e')$ are near-projections.
\end{lemma}
\begin{proof} We have a decomposition $P = P^+\oplus P^-$ with respect to
which $e'$ is strongly nilpotent. Now
\begin{align*}
e(1-e) & = \begin{pmatrix}
e^{++} & e^{+-} \\
e^{-+} & e^{--}
\end{pmatrix}\begin{pmatrix}
1-e^{++} & -e^{+-} \\
-e^{-+} & 1-e^{--}
\end{pmatrix} \\[1ex]
& = \begin{pmatrix}
e^{++}-(e^{++})^2-e^{+-}e^{-+} & -e^{++}e^{+-}+e^{+-}(1-e^{--}) \\
e^{-+}-e^{-+}e^{++}-e^{--}e^{-+} & -e^{-+}e^{+-}+e^{--}(1-e^{--})
\end{pmatrix} \\[1ex]
& = \begin{pmatrix}
e^{++}-(e^{++})^2-e^{+-}e^{-+} & -e^{++}e^{+-}+e^{+-}(1-e^{--}) \\
e^{-+} -e ^{-+}e^{++}-e^{--}e^{-+} &
-e^{-+}e^{+-}+(1-e^{--}) - (1-e^{--})^2
\end{pmatrix}.
\end{align*}
The matrix
$$e' = \left(\begin{array}{c|c}
e^{++} & -e^{+-} \\\hline
e^{-+} & 1-e^{--}
\end{array}\right)$$
denotes a strongly nilpotent representation of the complete quiver $Q_{2\mu}$
on $2 \mu$ vertices. In the following illustration $\mu = 1$:
$$
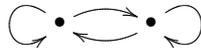
\xymatrix{
\bullet \ar@/^/[r] \ar@(ul,dl)[] & \bullet \ar@/^/[l] \ar@(ur,dr)[]
}$$
Now each entry in the $2\mu\times 2\mu$ matrix $e(1-e)$ above
is (the image of) a linear combination of paths of length at least one in the
quiver. Hence each entry of $(e(1-e))^N$ is the image of a sum of
paths of length at least $N$. It follows that $(e(1-e))^N = 0$ for some 
$N\geqslant 1$.

The pair $(P,e')$ is a near-projection since $e'\co P \to P$
is nilpotent.
\end{proof}

For $\mu = 1$ there is no difference between a near-projection
$(P,e,\{\pi_i\})$ and a near-projection $(P,e)$. For $\mu
\geqslant 2$ a near-projection $(P,e,\{\pi_i\})$ has $(P,e)$ a
near-projection (\fullref{near}) but the splitting
$(P,e) = (P^+,e^+) \oplus (P^-,e^-)$ given by \fullref{split} does not in general extend to a direct sum
decomposition of Seifert $A$--modules
$$(P,e,\{\pi_i\}) = (P^+,e^+,\{\pi^+_i\}) \oplus (P^-,e^-,\{\pi_i^-\}).$$
This is illustrated by the following example.

\begin{example}
Let $A$ be a field, and consider the $2$--component Seifert $A$--module
$(P,e,\{\pi_1,\pi_2\})$ given by
$$P = A^4,\quad e = \left(\begin{array}{cc|cc}
0 & 0 & 0 & 0 \\
0 & 1 & 1 & 0 \\ \hline
0 & 0 & 0 & 0 \\
1 & 0 & 0 & 1
\end{array}\right),~
\pi_1 = \begin{pmatrix} 1 & 0 & 0 & 0 \\
0 & 1 & 0 & 0 \\
0 & 0 & 0 & 0 \\
0 & 0 & 0 & 0 \end{pmatrix},~
\pi_2 = \begin{pmatrix} 0 & 0 & 0 & 0 \\
0 & 0 & 0 & 0 \\
0 & 0 & 1 & 0 \\
0 & 0 & 0 & 1 \end{pmatrix}.$$
In this case $e\co P \to P$ is a projection, with $e(1-e) = 0$.
This f.g.~projective Seifert $A$--module has just one submodule
$$(\wwbar{P},\bar{e},\{\bar{\pi}_1,\bar{\pi}_2\})
\subseteq (P,e,\{\pi_1,\pi_2\})$$
namely
$$\wwbar{P} = e(P) = \{(0,x,0,y)\in P\,\vert\,(x,y) \in A^2\}.$$
It is not possible to decompose $(P,e,\{\pi_1,\pi_2\})$ as a direct
sum, since $(\wwbar{P},\bar{e},\{\bar{\pi}_1,\bar{\pi}_2\})$ is not
a summand. Neither $e$ nor $1-e$ is nilpotent but
{\small$$\begin{array}{ll}
1-e+ez& = \begin{pmatrix}
1     & 0   & 0     & 0 \\
0     & z_1 & z_2-1 & 0 \\
0     & 0   & 1     & 0 \\
z_1-1 & 0   & 0     & z_2
\end{pmatrix}
 = \begin{pmatrix}
1     & 0   & 0     & 0 \\
0     & 1   & z_2-1 & 0 \\
0     & 0   & 1     & 0 \\
z_1-1 & 0   & 0     & 1
\end{pmatrix}
\begin{pmatrix}
1     & 0   & 0     & 0 \\
0     & z_1 & 0     & 0 \\
0     & 0   & 1     & 0 \\
0     & 0   & 0     & z_2
\end{pmatrix}
\end{array}$$}%
and
$$\begin{pmatrix}
0     & 0   & 0     & 0 \\
0     & 0   & z_2-1 & 0 \\
0     & 0   & 0     & 0 \\
z_1-1 & 0   & 0     & 0
\end{pmatrix}^2 = 0$$
so $1-e+ez$ is invertible. Moreover, $(P,e,\{\pi_1,\pi_2\})$ is a near-projection, with
$$P_1^+ = A \oplus 0 \oplus 0 \oplus 0,~
P_1^- = 0 \oplus A \oplus 0 \oplus 0,~
P_2^+ = 0 \oplus 0 \oplus A \oplus 0,~
P_2^- = 0 \oplus 0 \oplus 0 \oplus A$$
such that
$$e' = \begin{pmatrix}
0 & 0 & 0 & 0 \\
0 & 0 & 1 & 0 \\
0 & 0 & 0 & 0 \\
1 & 0 & 0 & 0\end{pmatrix}\co
P = P_1^+ \oplus P_1^- \oplus P_2^+ \oplus P_2^- \to
P = P_1^+ \oplus P_1^- \oplus P_2^+ \oplus P_2^-$$
is strongly nilpotent.
\end{example}

The main result of this section is:
\begin{thm}\label{characterize_primitives}
A f.g.~projective
Seifert $A$--module $(P,e,\{\pi_i\})$ is primitive if and only if it
is a near-projection.
\end{thm}
\begin{proof}
Suppose that $(P,e,\{\pi_i\})$ is a near-projection, with
$e' = \begin{pmatrix}
e^{++} & -e^{+-} \\
e^{-+} & 1-e^{--}
\end{pmatrix}$ strongly nilpotent. We have
{\small\begin{align*}
1{-}e{+}ez & =  1{-}e(1{-}z)\\[1ex]
& = \begin{pmatrix}
1{-}e^{++}(1{-}z) & -e^{+-}(1{-}z) \\
-e^{-+}(1{-}z)  & 1{-}e^{--}(1{-}z)
\end{pmatrix} \\[1ex]
& = \begin{pmatrix}
1{-}e^{++}(1{-}z) & e^{+-}(1{-}z^{-1}) \\
-e^{-+}(1{-}z)  & 1{-}(1{-}e^{--})(1{-}z^{-1})
\end{pmatrix}
\begin{pmatrix}
1 & 0 \\
0 & z
\end{pmatrix} \\[1ex]
& = \left(\begin{pmatrix}
1 & 0 \\
0 & 1
\end{pmatrix}{-}\begin{pmatrix}
e^{++} & -e^{+-} \\
e^{-+}& 1-e^{--}
\end{pmatrix}
\begin{pmatrix}
1{-}z & 0 \\
0 & 1{-}z^{-1}
\end{pmatrix}\right)
\begin{pmatrix}
1 & 0 \\
0 & z
\end{pmatrix}
\co(P^+ {\oplus} P^-)[F_{\mu}] \\
&\hskip242pt\longrightarrow (P^+ {\oplus} P^-)[F_{\mu}].
\end{align*}}
It follows from the strong nilpotence of $e'$ that
$e'((1-z)\oplus (1-z^{-1}))$ is nilpotent, and hence that
$$1-e(1-z) = (1-e'((1-z)\oplus (1-z^{-1})))(1\oplus z)\co 
(P^+ \oplus P^-)[F_{\mu}] \to (P^+ \oplus P^-)[F_{\mu}]$$
is an isomorphism, so that $B(P,e,\{\pi_i\}) = 0$ and $(P,e,\{\pi_i\})$ is 
primitive.

Conversely, suppose that $(P,e,\{\pi_i\})$ is a primitive
f.g.~projective Seifert $A$--module,
ie such that the $A[F_{\mu}]$--module morphism
$$1-e+ez\co P[F_{\mu}] \to P[F_{\mu}]$$
is an isomorphism. We shall use a variant $\wwbar{G}_\mu$ of
the Cayley tree $G_{\mu}$ (\fullref{cayley}) to prove that
$1-e+ez\co P[F_{\mu}] \to P[F_\mu]$ is a near-projection. Define
$$\wwbar{G}_{\mu}^{(0)} = F_\mu,
  \quad\wwbar{G}_{\mu}^{(1)} = \bigl\{(w,z_iw)\,\vert\,
    w \in F_{\mu},i \in \{1,2,\ldots,\mu\}\bigr\}$$
so that there is defined a right $F_\mu$--action
$$\wwbar{G}_{\mu} \times F_\mu \to \wwbar{G}_{\mu};~(w,g)
    \mapsto wg.$$
For each $i = 1,2,\ldots,\mu$ partition $F_{\mu}$ as
$$F_{\mu} = F^{+,i}_{\mu} \sqcup F^{-,i}_{\mu} \sqcup \{1\}$$
with $F_{\mu}^{+,i}$ (resp.  $F_{\mu}^{-,i}$) consisting of the
reduced words in $z_1,z_2,\ldots,z_{\mu}$ which start (resp.  do
not start) with $z_i$. Removing the edge $(w,z_iw)$ disconnects
$\wwbar{G}_{\mu}$, and the complement is a disjoint union of
trees
$$\wwbar{G}_{\mu}-\{(w,z_iw)\} = \wwbar{G}_\mu^+(w,z_iw) \sqcup \wwbar{G}_\mu^-(w,z_iw)$$
with
$$\wwbar{G}_\mu^+(w,z_iw)^{(0)} = F^{+,i}w,~\wwbar{G}_\mu^-(w,z_iw)^{(0)} = 
    (F^{-,i}\cup \{1\})w.$$
 In the diagram
\medskip

{\footnotesize
\Draw
\LineAt(-40,60,-40,20)
\LineAt(-40,20,40,20)
\LineAt(40,20,40,60)
\LineAt(80,60,80,20)
\LineAt(80,20,120,20)
\LineAt(120,-20,80,-20)
\LineAt(80,-20,80,-60)
\LineAt(80,-20,80,-60)
\LineAt(40,-60,40,-20)
\LineAt(40,-20,-40,-20)
\LineAt(-40,-20,-40,-60)
\LineAt(-80,-60,-80,-20)
\LineAt(-80,-20,-120,-20)
\LineAt(-120,20,-80,20)
\LineAt(-80,20,-80,60)
\MoveTo(33,3)
\Text(--$
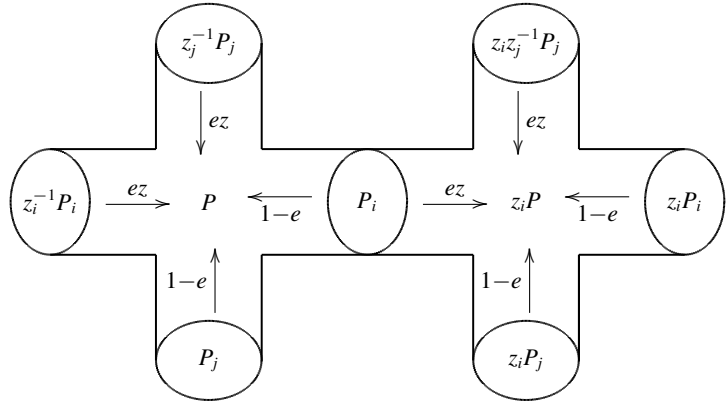
\xymatrix{\ar[r]^-{\di{ez}}&}$--)
\MoveTo(-87,3)
\Text(--$\xymatrix{\ar[r]^-{\di{ez}}&}$--)
\MoveTo(-33,-3)
\Text(--$\xymatrix{&\ar[l]^-{\di{1{-}e}}}$--)
\MoveTo(87,-3)
\Text(--$\xymatrix{&\ar[l]^-{\di{1{-}e}}}$--)
\MoveTo(72,30)
\Text(--$\xymatrix{\ar[d]^-{\di{ez}}&\\&}$--)
\MoveTo(67,-30)
\Text(--$\xymatrix{&\\ \ar[u]^-{\di{1{-}e}}&}$--)
\MoveTo(-48,30)
\Text(--$\xymatrix{\ar[d]^-{\di{ez}}&\\&}$--)
\MoveTo(-52,-30)
\Text(--$\xymatrix{&\\ \ar[u]^-{\di{1{-}e}}&}$--)
\MoveTo(0,0)
\Text(--$P_i$--)
\DrawOval(15,20)
\MoveTo(120,0)
\Text(--$z_iP_i$--)
\DrawOval(15,20)
\MoveTo(-120,0)
\Text(--$z^{-1}_iP_i$--)
\DrawOval(15,20)
\MoveTo(-60,60)
\Text(--$z^{-1}_jP_j$--)
\DrawOval(20,15)
\MoveTo(60,60)
\Text(--$z_iz^{-1}_jP_j$--)
\DrawOval(20,15)
\MoveTo(-60,-60)
\Text(--$P_j$--)
\DrawOval(20,15)
\MoveTo(60,-60)
\Text(--$z_iP_j$--)
\DrawOval(20,15)
\MoveTo(60,0)
\Text(--$z_iP$--)
\MoveTo(-60,0)
\Text(--$P$--)
\EndDraw
}

\medskip

we are placing the components of the range (resp.
domain) $P[F_{\mu}]$ at the vertices (resp. edges) of
$\wwbar{G}_{\mu}$, with the $A$--module $wP$ at $w \in
\wwbar{G}^{(0)}_{\mu}$, and the $A$--module $wP_i$ at
$(w,z_iw)\in \wwbar{G}^{(1)}_\mu$. An element
$$x \in P[F_{\mu}]  =  \!\!\!\!\sum\limits_{(w,z_iw)\in\wwbar{G}^{(1)}_{\mu}}
\!\!wP_i$$
    is sent to
$$(1-e)(x) + ez(x) \in P[F_{\mu}]  =  \!\!\sum\limits_{w \in
\wwbar{G}^{(0)}_{\mu}}\!\!wP,$$
as indicated by the arrows in the diagram.  For $i = 1,2,\ldots,\mu$
define the $A$--modules
\begin{align*}
P_i^+ & =  \biggl\{x \in P_i\,\vert\,
(1-e+ez)^{-1}ez(x) \in \!\!\sum\limits_{w \in F^{+,i}_{\mu}}\!\!wP \biggr\},\\
P_i^- & =  \biggl\{x \in P_i\,\vert\, (1-e+ez)^{-1}(1-e)(x)
\in\sum\limits_{j \neq i}P_j \oplus
\!\!\sum\limits_{w \in F^{-,i}_{\mu}}\!\!wP \biggr\}.
\end{align*}
 An element $x^+ \in P_i$ belongs to $P_i^+$ if and only
if there exist elements $y^+(w) \in P$ ($w \in F^{+,i}_{\mu}$)
such that
$$
ez(x^+) = (1-e+ez)\biggl(\sum\limits_{w \in F^{+,i}_{\mu}}wy^+(w)\biggr) \in
    \sum\limits_{w \in F_{\mu}^{+,i}}wP.\eqno{(*)}$$
There is one component $y^+(w)$ for each edge in
$\wwbar{G}^+_{\mu}(1,z_i)^{(1)}$, and one equation for each
vertex in $\wwbar{G}^+_{\mu}(1,z_i)^{(0)}$. Similarly, an
element $x^- \in P_i$ belongs to $P_i^-$ if and only if there
exist elements $y_j \in P_j$ ($j \neq i$) and $y^-(w) \in P$ ($w
\in F^{-,i}_{\mu}$) such that
$$((1-e)(x^-),0) = 
(1-e+ez)\biggl(\sum\limits_{j\neq i}y_j+\sum\limits_{w \in
F^{-,i}_{\mu}}wy^-(w)\biggr) \in P\oplus\sum\limits_{w\in
F^{-,i}_{\mu}}wP.\eqno{(**)}$$ There is one component $y_j$ ($j
\neq i$) or $y^-(w)$ for each edge
$\wwbar{G}^-_{\mu}(1,z_i)^{(1)}$, and one equation for each
vertex in $\wwbar{G}^-_{\mu}(1,z_i)^{(0)}$. For
$i = 1,2,\ldots,\mu$ partition
$$F^{+,i}_{\mu} = F^{++,i}_{\mu} \sqcup F^{-+,i}_{\mu},\quad
F^{-,i}_{\mu} = F^{+-,i}_{\mu} \sqcup F^{--,i}_{\mu}$$ with
$F^{\alpha +,i}_{\mu}$ consisting of the words
$w = \smash{z_{i_0}^{\epsilon_0}}\ldots\smash{z_{i_k}^{\epsilon_k}} \in F_{\mu}$
with $(i_0,\epsilon_0) = (i,+)$, $\epsilon_k = \alpha$, and $F^{\alpha
-,i}_{\mu}$ consisting of the words $w = \smash{z_{i_0}^{\epsilon_0}}\ldots
\smash{z_{i_k}^{\epsilon_k}} \in F_{\mu}$ with $(i_0,\epsilon_0)\neq
(i,+)$, $\epsilon_k = \alpha$. For any $x^+ \in P^+_i$ and $w \in
F^{\alpha+,i}_{\mu}$ we have that $y^+(w) \in \smash{P^\alpha_j}$, as
given by all the terms in $(*)$ involving
$\wwbar{G}^+(w,z_iw)$. Similarly, for any $x^- \in P^-_i$ and
$w \in F^{\alpha-,i}_{\mu}$ we have that $y^-(w) \in P^\alpha_j$,
as given by all the terms in $(**)$ involving
$\wwbar{G}^-(w,z_iw)$.

Regarded as an $A$--module isomorphism
$1{-}e{+}ez\co P[F_{\mu}] \to P[F_{\mu}]$ can be expressed as
\begin{multline*}
1{-}e{+}ez  =  \\
\begin{pmatrix} ez\vert & (1{-}e{+}ez)\vert & 0 \\
(1{-}e)\vert & 0 & (1{-}e{+}ez)\vert \end{pmatrix}\co
P_i \oplus  \biggl( \sum\limits_{w \in F^{+,i}_{\mu}}\!\!wP \biggr)\oplus
\biggl( \sum\limits_{j \neq i} P_j \oplus\!\! \sum\limits_{w \in
F^{-,i}_{\mu}}\!\!wP \biggr) \\
\longrightarrow \biggl(\sum\limits_{w \in F^{+,i}_{\mu}}\!\!wP \biggr)\oplus
\biggl(P\oplus\!\! \sum\limits_{w \in F^{-,i}_{\mu}}\!\!wP \biggr)
\end{multline*}
so that there is induced an $A$--module isomorphism
\begin{multline*}
\left[\begin{matrix} ez\vert \\ (1-e)\vert \end{matrix} \right]\co 
P_i \to \biggl(\coker\biggl((1-e+ez)\vert\co \!\!\sum\limits_{w \in
F^{+,i}_{\mu}}\!\!wP
\to \!\!\sum\limits_{w \in F^{+,i}_{\mu}}\!\!wP\biggr)\biggr)\\
\oplus \biggl(\coker\biggl((1-e+ez)\vert\co 
\sum\limits_{j \neq i} P_j \oplus \!\!\sum\limits_{w \in
F^{-,i}_{\mu}}\!\!wP
\to P\oplus \!\!\sum\limits_{w \in F^{-,i}_{\mu}}\!\!wP\biggr)\biggr)
\end{multline*}
and
$$P_i = P_i^+ \oplus P_i^-,$$
with
\begin{align*}
(1-e+ez)^{-1}ez(P^+_i) &\subseteq \sum\limits^{\mu}_{j = 1}
\sum\limits_{w \in F^{++,i}_{\mu}}\!\!wP^+_j \oplus
\sum\limits^{\mu}_{j = 1}
\sum\limits_{w \in F^{-+,i}_{\mu}}\!\!wP^-_j ,\\
(1-e+ez)^{-1}(1-e)(P^-_i) &\subseteq \sum\limits^{\mu}_{j = 1}
\sum\limits_{w \in F^{+-,i}_{\mu}}\!\!wP^+_j \oplus
\sum\limits^{\mu}_{j = 1} \sum\limits_{w \in F^{--,i}_{\mu}}\!\!wP^-_j.
\end{align*}
For $\alpha,\beta \in \{\pm\}$ let
$$e_{ji}^{\beta \alpha}\co P^{\alpha}_i \to P^{\beta}_j$$
be the $A$--module morphisms such that
$$e = \begin{pmatrix} e_{ji}^{++} & e_{ji}^{+-} \\ e_{ji}^{-+} &
e_{ji}^{--} \end{pmatrix}\co P = \sum\limits^{\mu}_{i = 1}(P_i^+ \oplus
P_i^-) \to P = \sum\limits^{\mu}_{j = 1}(P_j^+ \oplus P_j^-).$$
Let
$$\nu_{ji}^{\beta \alpha}(w)\co P^{\alpha}_i \to P^{\beta}_j$$
be the $A$--module morphisms such that
\begin{align*}
-(1{-}e{+}ez)^{-1}ez\vert & = 
\begin{pmatrix} \sum\limits^{\mu}_{j = 1}
  \sum\limits_{w \in F^{++,i}_{\mu}}w\nu^{++}_{ji}(w)\\
  \sum\limits^{\mu}_{j = 1}
  \sum\limits_{w \in F^{-+,i}_{\mu}}w\nu^{-+}_{ji}(w)\end{pmatrix}
\co P^+_i\\[-2ex]
&\hskip70pt\longrightarrow \sum\limits^{\mu}_{j = 1}
\sum\limits_{w \in F^{++,i}_{\mu}}wP^+_j \oplus
\sum\limits^{\mu}_{j = 1}
\sum\limits_{w \in F^{-+,i}_{\mu}}wP^-_j ,\\
-(1{-}e{+}ez)^{-1}(1-e)\vert & = 
\begin{pmatrix} \sum\limits^{\mu}_{j = 1}
\sum\limits_{w \in F^{+-,i}_{\mu}}w\nu^{+-}_{ji}(w)\\
\sum\limits^{\mu}_{j = 1}
\sum\limits_{w \in F^{--,i}_{\mu}}w\nu^{--}_{ji}(w)\end{pmatrix}\co
P^-_i\\[-2ex]
&\hskip70pt\longrightarrow \sum\limits^{\mu}_{j = 1}
\sum\limits_{w \in F^{+-,i}_{\mu}}wP^+_j \oplus
\sum\limits^{\mu}_{j = 1}
\sum\limits_{w \in F^{--,i}_{\mu}}wP^-_j.
\end{align*}
Composing with $1-e+ez$ gives
\begin{equation}
\begin{aligned}
-e_{ji}^{++}z_i & =  \sum\limits^{\mu}_{k = 1} \biggl(\sum\limits_{w \in
F^{++,i}_{\mu}}\!\!\!\!w(\delta_{jk}-e_{jk}^{++})\nu^{++}_{ki}(w)
+\!\!\!\! \sum\limits_{w \in
F^{-+,i}_{\mu}}\!\!\!\!wz_k(e_{jk}^{+-})\nu^{-+}_{ki}(w)\biggr)\co\\[-2ex]
&\hskip220pt P_i^+ \longrightarrow \!\!\!\!\sum\limits_{w \in
F^{++,i}_{\mu}}\!\!\!\!P_j^+,\\
-e_{ji}^{-+}z_i & =  \sum\limits^{\mu}_{k = 1}
\biggl(
\sum\limits_{w \in F^{++,i}_{\mu}}\!\!\!\!we_{jk}^{-+}\nu^{++}_{ki}(w)
+\!\!\!\!\sum\limits_{w \in
F^{-+,i}_{\mu}}\!\!\!\!wz_ke_{jk}^{--}\nu^{-+}_{ki}(w)\biggr)\co\\[-2ex]
&\hskip220pt P_i^+ \longrightarrow \!\!\!\!\sum\limits_{w \in
F^{-+,i}_{\mu}}\!\!\!\!wP_j^-,\\
-(-e_{ji}^{+-}) & =  \sum\limits^{\mu}_{k = 1} \biggl(\sum\limits_{w \in
F^{+-,i}_{\mu}}\!\!\!\!w(\delta_{jk}{-}e_{jk}^{++})\nu^{+-}_{ki}(w)
+ \!\!\!\!\sum\limits_{w \in
F^{-+,i}_{\mu}}\!\!\!\!wz_ke_{jk}^{+-}\nu^{--}_{ki}(w)\biggr)\co\\[-2ex]
&\hskip220pt P_i^- \longrightarrow
\!\!\!\!\sum\limits_{w \in F^{+-,i}_{\mu}}\!\!\!\!wP_j^+,\\
-(\delta_{ji}{-}e_{ji}^{--})& =  \sum\limits^{\mu}_{k = 1}
\biggl(\sum\limits_{w \in F^{+-,i}_{\mu}}\!\!\!\!we_{jk}^{-+}\nu^{+-}_{ki}(w)
+\!\!\!\!\sum\limits_{w \in F^{--,i}_{\mu}}\!\!\!\!wz_ke_{jk}^{--}\nu^{--}_{ki}(w)
\biggr)\co\\[-2ex]
&\hskip220pt P_i^- \longrightarrow \!\!\!\!\sum\limits_{w \in
F^{--,i}_{\mu}}\!\!\!\!wP_j^-.
\end{aligned}
\tag{$*$}
\label{eqstar}
\end{equation}
Comparing the coefficients of $z_i$ and 1 gives
\begin{align*}
-\begin{pmatrix} e_{ji}^{++} \\ e_{ji}^{-+} \end{pmatrix} & = 
\sum\limits^{\mu}_{k = 1}
\left(\begin{pmatrix} \delta_{jk}-e_{jk}^{++}  \\ -e_{jk}^{-+} \end{pmatrix}
\nu^{++}_{ki}(z_i)+\begin{pmatrix} e_{jk}^{+-}  \\ e_{jk}^{--} \end{pmatrix}
\nu^{-+}_{ki}(z_iz_k^{-1})\right)\co \\
&\hspace*{225pt}
P^+_i \to P^+_j \oplus P^-_j,\\
-\begin{pmatrix} -e_{ji}^{+-} \\ \delta_{ji}-e_{ji}^{--} \end{pmatrix} & = 
\sum\limits^{\mu}_{k = 1}\left(
\begin{pmatrix} \delta_{jk}-e_{jk}^{++}\\ -e_{jk}^{-+} \end{pmatrix}
\nu^{+-}_{ki}(1)+
\begin{pmatrix} e_{jk}^{+-}\\ e_{jk}^{--} \end{pmatrix}
\nu^{--}_{ki}(z_k^{-1})\right)\co \\
&\hspace*{225pt}
P^-_i \to P^+_j \oplus P^-_j.
\end{align*}
Writing
$$\begin{pmatrix} \nu^{++} & \nu^{+-} \\ \nu^{-+} & \nu^{--} \end{pmatrix} = 
\begin{pmatrix}
\nu_{ki}^{++}(z_i) & \nu_{ki}^{+-}(1) \\
\nu_{ki}^{-+}(z_iz_k^{-1}) & \nu_{ki}^{--}(z_k^{-1}) \end{pmatrix}\co
P^+ \oplus P^- \to P^+ \oplus P^-,$$
we thus have
$$-\begin{pmatrix} e^{++} & -e^{+-} \\ e^{-+} & 1-e^{--} \end{pmatrix} = 
\begin{pmatrix} 1-e^{++} & e^{+-} \\ -e^{-+} & e^{--} \end{pmatrix}
\begin{pmatrix} \nu^{++} & \nu^{+-} \\ \nu^{-+} & \nu^{--} \end{pmatrix}.$$
Let $Q^{+,-}_{\mu}$ be the quiver with $2\mu$ vertices
$(i,\pm)_{1 \leqslant i \leqslant \mu}$ and one edge
$(i_0,\epsilon_0) \to (i_1,\epsilon_1)$ for each pair of vertices
with $(i_0,\epsilon_0) \neq (i_1,-\epsilon_1)$. (The path ring is given by
\begin{multline*}
Q^{+,-}_{\mu}  =  \\
\Z[s]*\Z\Bigl[\pi^+_1,,\ldots,\pi^+_{\mu},
\pi^-_1,\ldots,\pi^-_{\mu}\,\vert\,\pi^{\alpha}_i\pi^{\beta}_j
 =  \delta_{\alpha\beta}\delta_{ij}\pi^{\alpha}_i,
\sum^{\mu}_{i = 1}(\pi^+_i+\pi^-_i) = 1\Bigr] \\
/\bigl\{\pi^{\alpha}_is\pi^{-\alpha}_i\bigr\}
\end{multline*}
where $\pi^{\alpha}_i s \pi^{\beta}_j$ ($(i,\alpha) \neq (j,-\beta)$)
corresponds to the unique path of length 1 from $(i,\alpha)$ to
$(j,\beta)$.) In the illustration $\mu = 2$:
$$
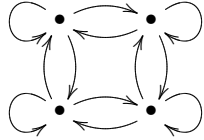
\xymatrix{
\bullet \ar@/^/[r] \ar@(ul,dl)[] \ar@/^/[d]
& \bullet \ar@/^/[l] \ar@(ur,dr)[] \ar@/^/[d] \\
\bullet \ar@/^/[r] \ar@(ul,dl)[] \ar@/^/[u]
& \bullet \ar@/^/[l] \ar@(ur,dr)[] \ar@/^/[u]}$$
Regard a word $w = z_{i_0}^{\epsilon_0}z_{i_1}^{\epsilon_1} \ldots
z_{i_k}^{\epsilon_k} \in F_{\mu}$ as a path of length $\vert w \vert = k$ in
$Q^{+,-}_{\mu}$
$$(i_0,\epsilon_0) \to (i_1,\epsilon_1) \to\cdots \to (i_k,\epsilon_k)$$
and for $k \geqslant 1$ define an $A$--module morphism
$\nu(w)\co P^{\epsilon_0}_{i_0} \to P^{\epsilon_1}_{i_1}$ as follows.
Define
$$[w] = \bigl[z_{i_0}^{\epsilon_0}z_{i_1}^{\epsilon_1}\bigr]
\bigl[z_{i_1}^{\epsilon_0}z_{i_2}^{\epsilon_1}\bigr] \ldots
\bigl[z_{i_{k-1}}^{\epsilon_{k-1}}z_{i_k}^{\epsilon_k}\bigr]\in F_{\mu}$$
with
$$\bigl[z_{i_0}^{\epsilon_0}z_{i_1}^{\epsilon_1}\bigr] = 
\begin{cases}
z_{i_0}^{\epsilon_0}&\hbox{\rm if $(\epsilon_0,\epsilon_1) = (+,+)$}\\[1ex]
z_{i_0}^{\epsilon_0}z_{i_1}^{\epsilon_1}&
\hbox{\rm if $(\epsilon_0,\epsilon_1) = (+,-)$}\\[1ex]
z_{i_1}^{\epsilon_1}&\hbox{\rm if $(\epsilon_0,\epsilon_1) = (-,-)$}\\[1ex]
1&\hbox{\rm if $(\epsilon_0,\epsilon_1) = (-,+)$}.
\end{cases}$$
For $k = 1$ set
$$\nu\bigl(z_{i_0}^{\epsilon_0}z_{i_1}^{\epsilon_1}\bigr) = 
\nu^{\epsilon_1\epsilon_0}_{i_1i_0}\bigl(\bigl[z_{i_1}^{\epsilon_1}z_{i_0}^{\epsilon_0}\bigr]\bigr)$$
and for $k \geqslant 2$ set
$$\nu\bigl(z_{i_0}^{\epsilon_0}z_{i_1}^{\epsilon_1}\ldots
z_{i_k}^{\epsilon_k}\bigr) = 
\nu\bigl(z_{i_{k-1}}^{\epsilon_{k-1}}z_{i_k}^{\epsilon_k}\bigr)\ldots
\nu\bigl(z_{i_1}^{\epsilon_1}z_{i_2}^{\epsilon_2}\bigr)
\nu\bigl(z_{i_0}^{\epsilon_0}z_{i_1}^{\epsilon_1}\bigr).$$
The identities
$$\nu(w) = \nu^{\epsilon_k\epsilon_0}_{i_ki_0}([w])\co 
P^{\epsilon_0}_{i_0} \to P^{\epsilon_k}_{i_k}$$
may be verified by induction on $k$,
since both sides satisfy the equations \eqref{eqstar} and so
\begin{align*}
-(1-e+ez)^{-1}ez\vert & =  \sum\limits_{w\in F^{+,i}_{\mu}}\!\!w\nu(w)\co 
P^+_i \to \!\!\sum\limits_{w \in F^{+,i}_{\mu}}wP,\\
-(1-e+ez)^{-1}(1-e)\vert & =  \sum\limits_{j \neq i}\nu^{+-}_{ji} +
\!\!\sum\limits_{w \in F^{-,i}_{\mu}}\!\! w\nu(w)\co
P^-_i \to \sum\limits_{j \neq i}P_j \oplus \!\!\sum\limits_{w \in F^{-,i}_{\mu}}wP.
\end{align*}
For $\alpha,\beta \in \{\pm\}$ let
$F^{\beta\alpha}_{\mu}$ be the set of paths
$$(i_0,\epsilon_0) \to (i_1,\epsilon_1) \to\cdots \to (i_k,\epsilon_k)$$
in $Q^{+,-}_{\mu}$ with $\epsilon_0 = \alpha$, $\epsilon_k = \beta$.
The $A[F_{\mu}]$--module endomorphism
$$\nu' = \begin{pmatrix} z & 0 \\ 0 & 1 \end{pmatrix}
\begin{pmatrix} \nu^{++} & \nu^{+-} \\ \nu^{-+} & \nu^{--} \end{pmatrix}
\begin{pmatrix} 1 & 0 \\ 0 & z^{-1} \end{pmatrix}\co
(P^+ \oplus P^-)[F_{\mu}] \to (P^+ \oplus P^-)[F_{\mu}]  $$
is such that for any $N \geqslant 1$
\begin{multline*}
(\nu')^N = 
\begin{pmatrix} \sum\limits_{w\in F^{++}_{\mu},\vert w \vert  = N}w\nu^{++}(w) &
\sum\limits_{w\in F^{+-}_{\mu},\vert w \vert  = N}w\nu^{+-}(w)\\
\sum\limits_{w\in F^{-+}_{\mu},\vert w \vert  = N}w\nu^{-+}(w) &
\sum\limits_{w\in F^{--}_{\mu},\vert w \vert  = N}w\nu^{--}(w)
\end{pmatrix}\co\\
(P^+ \oplus P^-)[F_{\mu}] \to (P^+ \oplus P^-)[F_{\mu}].
\end{multline*}
If $N \geqslant 1$ is so large that
$$(1-e+ez)^{-1} = \sum\limits_{w \in F_{\mu},\vert w \vert<N}a_ww\co
P[F_{\mu}] \to P[F_{\mu}]\quad (a_w \in \Hom_A(P_{i_0},P_{i_k})),$$
then for any word $w \in F_{\mu}$ of length $\vert w \vert = k > N$
$$\nu(w) = 0\co P^{\epsilon_0}_{i_0} \to P^{\epsilon_k}_{i_k}.$$
The $2\mu$--component Seifert module
$$(P',\nu',\pi') = \biggl(P^+ \oplus P^-,
\begin{pmatrix} \nu^{++} & \nu^{+-} \\ \nu^{-+} & \nu^{--} \end{pmatrix},
\{\pi_i^+ \oplus \pi_i^-\}\biggr)$$
is strongly nilpotent, with $(\nu' z')^N = 0$, regarding
$F_{2\mu}$ as free group on $2\mu$ generators
$z'_1,z'_2,\ldots,z'_{2\mu}$ and letting
$$z' = \begin{pmatrix}
z'_1 & 0 & \ldots & 0 \\
0 & z'_2 & \ldots & 0 \\
\vdots & \vdots & \ddots & \vdots \\
0 & 0 & \ldots & z'_{2\mu}
\end{pmatrix}\co P'[F_{2\mu}] \to P'[F_{2\mu}].$$
Define the $2\mu$--component Seifert module
$$(P',e',\pi') = \biggl(P^+\oplus P^-,\begin{pmatrix}
e^{++} & -e^{+-} \\
e^{-+} & 1-e^{--}
\end{pmatrix},\{\pi^+_i\oplus \pi^-_i\}\biggr).$$
Applying the augmentation $\epsilon\co z_i \mapsto 1$ to the
$A[F_{\mu}]$--module morphisms
\begin{align*}
-\sum\limits_{w \in F^{+,i}_{\mu}}w\nu(w)&\co
P^+_i[F_\mu] \xymatrix{\ar[r]^{ez}&} P[F_{\mu}]
\xymatrix@C+20pt{\ar[r]^{(1-e+ez)^{-1}}&} P[F_{\mu}],\\[1ex]
-\bigg(\sum\limits_{j \neq i}\nu^{+-}_{ji}+\sum\limits_{w \in
F^{-,i}_{\mu}}w\nu(w)\bigg)&\co
P^-_i[F_\mu] \xymatrix{\ar[r]^{1-e}&} P[F_{\mu}]
\xymatrix@C+20pt{\ar[r]^{(1-e+ez)^{-1}}&} P[F_{\mu}]
\end{align*}
shows that the components of $e'$ are given by linear combinations
of paths of length $\geqslant 1$
\begin{align*}
e^{++} & =  -\sum\limits_{w \in F^{++}_{\mu}}\nu(w)\co P^+ \to P^+,\\
e^{-+} & =  -\sum\limits_{w \in F^{+-}_{\mu}}\nu(w)\co P^+ \to P^-,\\
-e^{+-} & =  -\sum\limits_{w \in F^{-+}_{\mu}}\nu(w)\co P^- \to P^+,\\
1-e^{--} & =  -\sum\limits_{w \in F^{--}_{\mu}}\nu(w)\co P^- \to P^-.
\end{align*}
The $A[F_{2\mu}]$--module endomorphism $e'z'\co P'[F_{2\mu}] \to P'[F_{2\mu}]$
is nilpotent, with
$$(e'z')^N = 0,$$
so that $(P,e,\pi)$ is strongly nilpotent.
\end{proof}

This completes the proof of \fullref{thm3} of the Introduction.

\section{Algebraic $K$--theory} \label{ktheory}

We shall obtain our results on the algebraic $K$--theory of $A[F_\mu]$
and Blanchfield and Seifert modules using the Waldhausen \cite{Wald3}
algebraic $K$--theory of categories with cofibrations and weak
equivalences, and the noncommutative localization algebraic $K$--theory
exact sequence of Neeman and Ranicki \cite{NR1,NR2}.

\subsection{The algebraic $K$--theory of exact categories}

The higher algebraic $K$--groups $K_n(\E)$ of an exact category $\E$ are
defined by Quillen \cite{Q} to be the homotopy groups of
a connective spectrum $K(\E)$
$$\pi_n(K(\E)) = K_n(\E)\quad (n \geqslant 0)$$
with $K_0(\E)$ the Grothendieck class group.  The idempotent completion
$\E \to \PP(\E)$ induces an injection $K_0(\E) \to K_0(\PP(\E))$ and
isomorphisms $K_n(\E) \to K_n(\PP(\E))$ for $n \geqslant 1$, by the
cofinality theorem of Grayson \cite{Gr2}.  The lower $K$--groups
$K_n(\E)$ ($n \leqslant -1$) are defined by Schlichting \cite{Schl}
(following on from the definitions of Karoubi and Pedersen--Weibel for
the lower $K$--groups of filtered additive categories) to be the lower
homotopy groups of a nonconnective spectrum $K\PP(\E)$ such that
$$\pi_n(K\PP(\E)) = K_n(\PP(\E))\quad (n \in \Z),$$
with $K_n(\E) = K_n(\PP(\E))$ for $n \neq 0$.

The algebraic $K$--groups of a ring $R$ are the algebraic
$K$--groups of the idempotent complete exact category $\E = \Proj(R)$
of f.g.~projective $R$--modules
$$K_n(R) = K_n(\Proj(R))\quad (n \in \Z),$$
as defined for $-\infty < n \leqslant 1$ in Bass \cite{B2}, and for $2
\leqslant n < \infty$ in Quillen \cite{Q}. The nonconnective spectrum
defined by $K(R) = K\PP(\Proj(R))$ has homotopy groups $\pi_*(K(R)) = K_*(R)$.

A \emph{Waldhausen category} $(\C,w)$ is a small category $\C$ with cofibrations
together with a subcategory $w \subset \C$ of weak equivalences satisfying
the axioms of \cite{Wald3}.  As usual, there is defined a connective
algebraic $K$--theory spectrum
$$K(\C,w) = \Omega \vert w S_{\bullet} \C\vert $$
with homotopy groups the algebraic $K$--theory groups
$$K_n(\C,w) = \pi_n(K(\C,w))\quad (n \geqslant 0).$$
A functor $F\co (\C,w) \to (\C',w')$ of Waldhausen categories
induces a long exact sequence of algebraic $K$--groups
$$\xymatrix@C-10pt
{\cdots \ar[r]&K_{n+1}(F) \ar[r] & K_n(\C,w)
\ar[r]^-{\di{F}} & K_n(\C',w') \ar[r] & \cdots \ar[r]&K_0(F) \ar[r] &0}$$
with $K_n(F) = \pi_n(F\co K(\C,w) \to K(\C',w'))$ ($n \geqslant 0$).

As in Thomason and Trobaugh \cite[1.9]{TT} we shall only be considering
Waldhausen categories $(\C,w)$ which are `complicial biWaldhausen', so
that in particular $\C$ is a full subcategory of the category of chain
complexes in an abelian category ${\mathcal A}$, the cofibrations are
chain maps which are split injections in each degree, $w$ contains the
quasi-isomorphisms ( =  the chain maps inducing isomorphisms in homology),
and which in addition are closed under the formation
of canonical homotopy pushouts and pullbacks.

The \emph{homotopy} (or \emph{derived}) \emph{category} \cite[page~269]{TT}
of a Waldhausen category $(\C,w)$ is the category of fractions
$$D(\C,w) = w^{-1}\C,$$
which is a triangulated category under the above hypotheses.
The idempotent completion $\PP D(\C,w)$ is then
also triangulated (Balmer and Schlichting \cite{BS}),
and the class groups $K_0(D(\C,w))$, $K_0(\PP D(\C,w))$ are defined,
with $K_0(D(\C,w)) = K_0(\C,w)$.
Schlichting \cite{Schl} defined the lower $K$--groups $K_n(\PP D(\C,w))$
for $n \leqslant -1$ for Waldhausen categories as above,
and constructed a nonconnective spectrum $K\PP(\C,w)$
with homotopy groups
$$\pi_n(K\PP(\C,w)) = K\PP_n(\C,w) = \begin{cases}
K_n(\C,w)&\hbox{\rm for $n \geqslant 1$}\\
K_0(\PP D(\C,w))&\hbox{\rm for $n = 0$}\\
K\PP_n(\C,w)&\hbox{\rm for $n\leqslant -1$.}
\end{cases}$$
A functor $F\co (\C,w) \to (\C',w')$ of Waldhausen categories
induces a long exact sequence of algebraic $K$--groups
$$\xymatrix@C-10pt{\cdots \ar[r]&K\PP_{n+1}(F) \ar[r] & K\PP_n(\C,w)
\ar[r]^-{\di{F}} & K\PP_n(\C',w') \ar[r] & K\PP_n(F)\ar[r] &\cdots,}$$
with $K\PP_n(F) = \pi_n(F\co K\PP(\C,w) \to K\PP(\C',w'))$ $(n \in \Z)$.

Given an exact category $\E$ let $C^b(\E)$ be the
category of bounded chain complexes in $\E$ and chain maps.
An object $C$ in $C^b(\E)$ is \emph{acyclic} (in the sense of Keller \cite[Chapter 11]{K})
if each differential $d\co C_r \to C_{r-1}$ factors as
$C_r \to Z_r \to C_{r-1}$ with
$$0 \to Z_{r+1} \to C_r \to Z_r \to 0$$
exact. A morphism $f\co C \to D$ in $C^b(\E)$ is a \emph{quasi-isomorphism}
if the mapping cone $\C(f)$ is chain equivalent to an acyclic complex.
If $\E$ is fully embedded in an abelian category $\A$ with the embedding
closed under extensions and the idempotent completion $\PP(\E)$ is closed
under taking kernels of surjections then a quasi-isomorphism is the same
as a chain map inducing isomorphisms in homology in the ambient abelian
category $\A$ \cite[Appendix A]{TT}.

Let $(\C^b(\E),w_{\E})$ be the Waldhausen category with cofibrations
the chain maps which are degreewise split injections, and
$w_\E\subset\C^b(\E)$ the subcategory of quasi-isomorphisms.  The
derived category
$$D^b(\E) = D(\C^b(\E),w_{\E})$$
is the category of bounded chain complexes
in $\E$ and fractions of chain homotopy classes of chain maps,
with denominators quasi-isomorphisms. As usual, let $K^b(\E)$ be the
category of bounded chain complexes in $\E$ and chain homotopy classes
of chain maps, and let $wK_{\E} \subset K^b(\E)$ be the subcategory
of quasi-isomorphisms: the localization
$$D^b(\E) = (wK_{\E})^{-1}K^b(\E)$$
has both a left and a right calculus of fractions.  The derived
category $D^b(\E)$ is a triangulated category \cite[1.9.6]{TT}.
Balmer and Schlichting \cite[2.12]{BS} prove that the idempotent
completion of the derived category is the derived category of the
idempotent completion
$$\PP D^b(\E) = D^b(\PP(\E))$$
and the algebraic $K$--groups are such that
$$\begin{cases}
K_n(\C^b(\E),w_\E) = K_n(\E)&
\hbox{\rm for $n \geqslant 0$ (Gillet \cite{Gi})}\\[1ex]
K\PP_n(\C^b(\E),w_\E) = K_n(\PP(\E))
&\hbox{\rm for $n \in \Z$ (Schlichting \cite{Schl})}.
\end{cases}$$
By \cite[1.9.2]{TT} the Waldhausen category defined in the same way but
with cofibrations the chain maps which are degreewise admissible
monomorphisms has the same algebraic $K$--theory.

\begin{definition} \label{w}
Let $F\co \E \to \D$ be a functor of exact categories.

{\rm (i)}\qua The \emph{algebraic $K$--groups} $K\PP_*(\E,\D)$ are
$$K\PP_n(\E,\D) = K\PP_n(C^b(\E,\D),w_{(\E,\D)})\quad (n \in \Z)$$
with $(\C^b(\E,\D),w_{(\E,\D)}) \subset (\C^b(\E),w_{\E})$ the
Waldhausen subcategory with $\C^b(\E,\D)\subset \C^b(\E)$ the full
subcategory with objects the bounded chain complexes $C$ in $\E$
which are chain equivalent in $\D$ to acyclic complexes, and
$$w_{(\E,\D)} = w_{\E}\cap \C^b(\E,\D) \subset \C^b(\E,\D)$$
the subcategory of the quasi-isomorphisms.

{\rm (ii)}\qua The \emph{algebraic $\Gamma K$--groups} of $F$ are
$$\Gamma K_n(F) = K\PP_n(C^b(\E),w_{\D})\quad (n \in \Z)$$
with $w_{\D} \subset \C^b(\E)$ the
subcategory with morphisms the chain maps in $\E$ which become
quasi-isomorphisms in $\D$, or equivalently such that
the mapping cones are in $\C^b(\E,\D)$.
\end{definition}

The groups $\Gamma K_*(F)$ are the algebraic
$K$--theory analogues of the algebraic $L$--theory groups
$\Gamma_*(F)$ of Cappell and Shaneson \cite{CS1}.

\begin{thm} \label{long}
Let $F\co \E \to \D$ be a functor of exact categories.

{\rm (i)}\qua The algebraic $K$--groups fit into a commutative braid
of exact sequences
$$
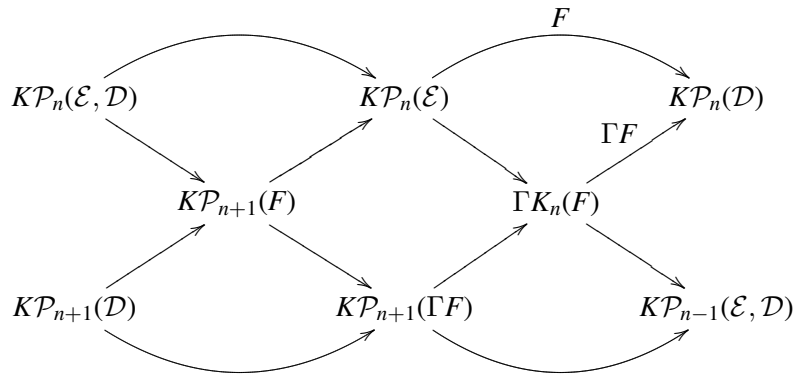
\xymatrix@C-15pt{
K\PP_n(\E,\D)\ar[dr] \ar@/^2pc/[rr]^{} && K\PP_n(\E)  \ar[dr]
\ar@/^2pc/[rr]^-{\di{F}} &&K\PP_n(\D)\\&
K\PP_{n+1}(F)\ar[ur] \ar[dr] &&
\Gamma K_n(F) \ar[ur]^-{\di{\Gamma F}}\ar[dr]&&\\
K\PP_{n+1}(\D)  \ar[ur] \ar@/_2pc/[rr]_-{}&&K\PP_{n+1}(\Gamma F)
    \ar[ur]\ar@/_2pc/[rr]_{}&&K\PP_{n-1}(\E,\D)}$$
with $\Gamma F\co (C^b(\E),w_{\D}) \to (C^b(\D),w_{\D})$ induced by $F$.

{\rm (ii)}\qua If $\Gamma F\co \PP D(\C^b(\E),w_{\D})  \to \PP
D(\C^b(\D),w_{\D})$ is an equivalence of categories then
$$K\PP_*(\Gamma F) = 0,~K\PP_{*+1}(F)\cong~K\PP_*(\E,\D),~
\Gamma K_*(F)\cong~K\PP_*(\D)$$
and the braid of {\rm (i)} collapses to the exact sequence
$$\xymatrix{ \cdots \ar[r]& K\PP_{n+1}(\D) \ar[r]&
K\PP_n(\E,\D) \ar[r] & K\PP_n(\E)
\ar[r]^-{\di{F}} & K\PP_n(\D) \ar[r] &\cdots}$$
{\rm (iii)}\qua The hypothesis of {\rm (ii)} is satisfied if
$F\co \E\to\D = \Sigma^{-1}\E$ is the canonical functor to a
category of fractions and $\D$ has a calculus of left fractions.
\end{thm}
\begin{proof} (i)\qua The cases $n \geqslant 0$ are a direct application of
the version of the localization theorem of \cite[1.6.4]{Wald3} stated
in Theorem 2.3 and Lemma 2.5 of Neeman and Ranicki \cite{NR2}, with
$$\begin{array}{l}
{\mathcal R}^c = D(\C^b(\E,\D),w_{(\E,\D)})
\subset {\mathcal S}^c = D(\C^b(\E),w_{\E}),\quad
{\mathcal S}^c/{\mathcal R}^c \approx D(\C^b(\E),w_{\D}),\\
{\bf R}  =  (\C^b(\E,\D),w_{(\E,\D)}),\quad
{\bf S}  =  (\C^b(\E),w_\E),\quad {\bf T}  =  {\bf S}_{\bf R}  =  (\C^b(\E),w_\D)
\end{array}$$
giving a fibration sequence of connective spectra
$$K(\C^b(\E,\D),w_{(\E,\D)}) \to  K(\C^b(\E),w_\E) \to K(\C^b(\E),w_{\D}).$$
The cases $n <0$ follow from Theorems 2.4, 3.7 of \cite{NR2} and
Schlichting \cite[Theorems 1,9]{Schl}, which give a fibration sequence of
nonconnective spectra
$$K\PP(\C^b(\E,\D),w_{(\E,\D)}) \to K\PP(\C^b(\E),w_\E) \to
K\PP(\C^b(\E),w_{\D}).$$
(ii)\qua This is a direct application of the Approximation Theorem
of Waldhausen \cite[Theorem~1.6.7]{Wald3}: if $F\co (\C,w) \to (\C',w')$ is a functor which induces an equivalence
of the homotopy categories $F\co D(\C,w) \to D(\C',w')$ then
$F\co K(\C,w) \to K(\C',w')$ is a homotopy equivalence inducing isomorphisms
$F\co K_*(\C,w)\cong K_*(\C',w')$. Similarly, if $F\co \PP D(\C,w) \to \PP
D(\C',w')$ is an equivalence there are induced isomorphisms
$F\co K\PP_*(\C,w)\cong K\PP_*(\C',w')$ (\cite{Schl}).

(iii)\qua Every object $D$ in $\C^b(\D)$ is chain equivalent to $F(E)$ for an
object $E$ in $\C^b(\E)$, and the functors $F\co C^b(\E) \to C^b(\D)$,
$F\co D(\C^b(\E),w_{\D}) \to D(\C^b(\D),w_{\D})$ are localizations.
\end{proof}

\begin{definition}
(i)\qua Write the algebraic $K$--groups of the exact categories
$\Prim(A)$, $\Sei(A)$, $\Bla(A)$, $\Flk(A)$ as
$$\begin{array}{l}
\PRIM_*(A) = K_*(\Prim(A)),~\SEI_*(A) = K_*(\Sei(A)),\\[1ex]
\BLA_*(A) = K_*(\Bla(A)),~\FLK_*(A) = K_*(\Flk(A)).
\end{array}$$
(ii)\qua Write the algebraic $K$--groups of the idempotent completion
of the homotopy category of
$(\C^b(\Sei(A),\Bla(A)),w_{(\Sei(A),\Bla(A))})$ as
$$(\SEI,\BLA)_*(A) = K\PP_*(\C^b(\Sei(A),\Bla(A)),w_{(\Sei(A),\Bla(A))}).$$
\end{definition}

\begin{proposition} \label{blink}
The covering functor $B\co \Sei(A) \to \Bla(A)$
induces morphisms $B\co \SEI_*(A) \to \BLA_*(A)$ which fit into
a long exact sequence
$$\xymatrix@C-5pt{ \cdots \ar[r]& (\SEI,\BLA)_n(A) \ar[r]&
\SEI_n(A) \ar[r]^-{\di{B}} & \BLA_n(A)
\ar[r] & (\SEI,\BLA)_{n-1}(A) \ar[r] &\cdots}$$
with
$${\rm im}(B\co \SEI_0(A) \to \BLA_0(A)) = \FLK_0(A) \subseteq \BLA_0(A).$$
\end{proposition}
\begin{proof} Apply \fullref{long} (iii) with
$$F\co \E = \Sei(A) \to \D = \Xi^{-1}\Sei(A)~\approx~\Flk(A),$$
noting that $\Sei(A)$ is idempotent complete (\fullref{idem} (i)), that
$\Xi^{-1}\Sei(A)\approx \Flk(A)$ has a
left calculus of fractions by \fullref{calculus}, and
that $\Bla(A)\approx \PP(\Flk(A))$ (\fullref{idem}(ii)).
\end{proof}

In the next section it will be shown that the functor
$$\Prim(A) \to \C^b(\Sei(A),\Bla(A));~(P,e,\{\pi_i\}) \mapsto
(\cdots \to 0 \to (P,e,\{\pi_i\}))$$
induces isomorphisms of algebraic $K$--groups $\PRIM_*(A)\cong
(\SEI,\BLA)_*(A)$.

\subsection{The algebraic $K$--theory of noncommutative localizations}

Given a ring $R$ let $\Mod(R)$ be the abelian category of $R$--modules,
so that $\Proj(R) \subset \Mod(R)$ is an exact subcategory.
Write the Waldhausen category of $\Proj(R)$ as
$$(\C^b(R),w_R) = (\C^b(\Proj(R)),w_{\Proj(R)}).$$
An object in $\C^b(R)$ is a bounded chain complex $C$ of f.g.~projective
$R$--modules; $C$ is acyclic if and only if $H_*(C) = 0$.
A morphism $f\co C \to D$ in $\C^b(R)$ is a chain map;
$f$ is in $w_R$ if and only if $f_*\co H_*(C) \to H_*(D)$ is an isomorphism.
The algebraic $K$--groups of $R$ are given by
$$K_*(R) = K_*(\Proj(R)) =  K\PP_*(\C^b(R),w_R).$$
\indent
A ring morphism $\F\co R \to S$ induces a functor of abelian categories
$$\F = S\otimes_R-\co \Mod(R) \to \Mod(S);~P \mapsto S\otimes_RP$$
which restricts to an exact functor
$F\co \Proj(R) \to \Proj(S)$. There is also induced
a functor of Waldhausen categories
$$\F\co (\C^b(R),w_R) \to (\C^b(S),w_S);~C\mapsto S\otimes_RC.$$
The relative homotopy groups of $\F\co K(R) \to K(S)$
are the relative $K$--groups $K_*(\F)$ in the long exact sequence
$$\xymatrix{\cdots \ar[r]&  K_n(R) \ar[r]^{\F }& K_n(S) \ar[r] &
    K_n(\F) \ar[r] & K_{n-1}(R) \ar[r] & \cdots.}$$

Let $R$ be a ring, and let $\Sigma$ be a set of morphisms
of f.g.~projective $R$--modules. A ring morphism $R\to T$ is
\emph{$\Sigma$--inverting}
if each $(s\co P \to Q) \in \Sigma$ induces a
$T$--module isomorphism $1\otimes s\co T\otimes_RP \to T\otimes_RQ$.
By Cohn \cite{Co} there exists a
\emph{universal $\Sigma$--inverting localization} ring morphism
$$\F\co R \to S = \Sigma^{-1}R$$
such that any $\Sigma$--inverting ring morphism $R \to T$ has a unique
factorization
$$\xymatrix{R\ar[r]^-{\F }& S \ar[r] &T.}$$
The category of fractions $\Sigma^{-1}\Proj(R)$ is equivalent to the
full subcategory
$$\Proj_R(S) \subseteq \Proj(S)$$
with objects isomorphic to the f.g.~projective $S$--modules
$\Sigma^{-1}P = S\otimes_RP$ induced from f.g.~projective $R$--modules
$P$, and $\Proj(S) = \PP(\Proj_R(S))$ is the idempotent completion.

\begin{definition}
(i)\qua For any ring morphism $\F\co R \to S$
write the Waldhausen categories defined in \fullref{w} as
\begin{align*}
(\C^b(\Proj(R),\Proj(S)),w_{(\Proj(R),\Proj(S))}) & =  (\C^b(R,S),w_{(R,S)}),\\
(\C^b(\Proj(R)),w_S) & =  (\C^b(R),w_S)
\end{align*}
with corresponding nonconnective algebraic $K$--theory spectra
$$K\PP(\C^b(R,S),w_{(R,S)}) = K(R,S),~K\PP(\C^b(R),w_S) = \Gamma K(\F)$$
and algebraic $K$--groups $K_*(R,S)$, $\Gamma K_*(\F)$.
An object in $\C^b(R,S)$ is a bounded chain complex $C$ of f.g.~projective
$R$--modules such that $H_*(S\otimes_RC) = 0$.
A morphism $f\co C \to D$ in $\C^b(R,S)$ is a chain map;
$f$ is in $w_{(R,S)}$ if and only if $f_*\co H_*(C) \to H_*(D)$ is an
isomorphism. A morphism $f\co C \to D$ in $\C^b(R)$ is
in $w_S$ if and only if $1\otimes f\co H_*(S\otimes_RC) \to H_*(S\otimes_RD)$ is
an isomorphism.

(ii)\qua For an injective universal localization $\F\co R \to S = \Sigma^{-1}R$
let $H(R,\Sigma)$ be the exact category of \emph{h.d.~1
$\Sigma$--torsion $R$--modules}, ie the cokernels of injective
morphisms $s\co P \to Q$ of f.g.~projective $R$--modules which induce
an $S$--module isomorphism $1\otimes s\co S\otimes_RP \to S\otimes_RQ$
(eg if $s \in \Sigma$).

\rm (iii)\qua (Neeman and Ranicki \cite{NR1,NR2})
A universal localization $\F\co R \to S = \Sigma^{-1}R$ is \emph{stably flat} if
$$\Tor^R_i(S,S) = 0 \quad(i \geqslant 1).$$
\end{definition}

In particular, a universal localization $\F\co R \to S$ is stably flat if
$S$ has flat dimension $\leqslant 1$ as an $R$--module, ie  if there
exists a 1--dimensional
flat $R$--module resolution
$$0 \to F_1 \to F_0 \to S \to 0.$$

\begin{proposition} \label{Gamma}
{\rm (i)}\qua For any ring morphism $\F\co R \to S$ the functor
$${\Gamma \mathcal F}\co (\C^b(R),w_S) \to (\C^b(S),w_S);~C \mapsto  S\otimes_RC$$
induces morphisms of algebraic $K$--groups
$\Gamma \F\co \Gamma K_*(\F) \to K_*(S)$ which fit into a
commutative braid of exact sequences
$$
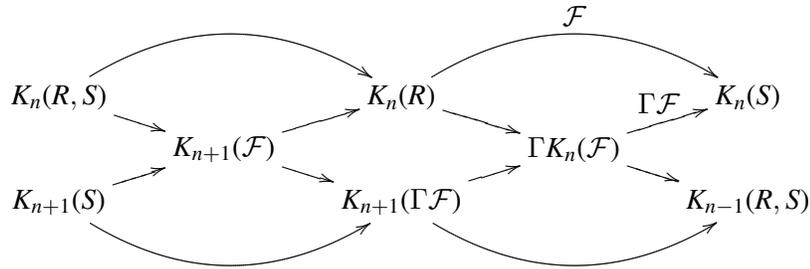
\xymatrix@C-5pt@R-20pt{
K_n(R,S)\ar[dr] \ar@/^2pc/[rr]^{} && K_n(R)  \ar[dr]
\ar@/^2pc/[rr]^-{\di{\F}} &&K_n(S)\\&
K_{n+1}(\F)\ar[ur] \ar[dr] &&
\Gamma K_n(\F) \ar[ur]^-{\di{\Gamma\mathcal F}} \ar[dr]&&\\
K_{n+1}(S)  \ar[ur] \ar@/_2pc/[rr]_-{}&&K_{n+1}(\Gamma\F)
    \ar[ur]\ar@/_2pc/[rr]_{}&&K_{n-1}(R,S)}$$
{\rm (ii)}\qua For any universal localization $\F\co R \to S = \Sigma^{-1}R$
$$\begin{array}{l}
\Gamma K_n(\F) = K_n(S),~K_n(\F) = K_{n-1}(R,S),~
K_n(\Gamma\F) = 0\quad (n \leqslant 1).
\end{array}$$
{\rm (iii)}\qua For a stably flat universal localization
$\F\co R \to S = \Sigma^{-1}R$
$$\Gamma K_*(\F) = K_*(S),~K_{*+1}(\F) = K_*(R,S),~K_*(\Gamma\F) = 0,$$
and there is induced a localization exact sequence in the algebraic $K$--groups
$$\xymatrix@C-7pt{\cdots \ar[r]& K_n(R,S) \ar[r]&
K_n(R) \ar[r]^-{\di{\F}}& K_n(S) \ar[r]& K_{n-1}(R,S) \ar[r]& \cdots.}$$
{\rm (iv)}\qua For an injective universal localization
$\F\co R \to S = \Sigma^{-1}R$ there is defined an equivalence
of homotopy categories
$$D(\C^b(R,S),w_{(R,S)})~\approx~D(C^b(H(R,\Sigma)),w_{H(R,\Sigma)})$$
inducing isomorphisms
$$K_*(R,S)\cong~ K_*(H(R,\Sigma)).$$
{\rm (v)}\qua For an injective stably flat universal localization
$\F\co R \to S = \Sigma^{-1}R$ there is defined
a localization exact sequence in the algebraic $K$--groups
$$\xymatrix@C-7pt{\cdots \ar[r]& K_n(H(R,\Sigma)) \ar[r]&
K_n(R) \ar[r]^-{\di{\F}}& K_n(\Sigma^{-1}R) \ar[r]& K_{n-1}(H(R,\Sigma))
\ar[r]& \cdots}$$
\end{proposition}
\begin{proof} (i)\qua Immediate from \fullref{long} (i) and (ii) applied to
$\F\co \C^b(R) \to \C^b(S)$.

(ii)--(v)\qua See Neeman and Ranicki \cite{NR1,NR2}.
\end{proof}

\subsection{Triangular matrix rings}

We refer to Haghany and Varadarajan \cite{HV} for the
general theory of modules over triangular matrix rings,
and to Schofield \cite{Sc}, Ranicki \cite{RNLAT} and
Sheiham \cite{Sh4} for previous accounts of the universal
localization of triangular matrix rings.

\begin{proposition} \label{triangular}
Let
$$A = \begin{pmatrix} A_1 & B \\ 0 & A_2
\end{pmatrix}$$
be the triangular $2\times 2$ matrix ring
defined by rings $A_1,A_2$ and an $(A_1,A_2)$--bimodule $B$. 

{\rm (i)}\qua An $A$--module $L = (L_1,L_2,\lambda)$ is defined by an $A_1$--module
$L_1$, an $A_2$--module $L_2$ and an $A_1$--module morphism
$\lambda\co B\otimes_{A_2}L_2 \to L_1$. As an additive group $L = L_1 \oplus L_2$,
written
$\bigl(\begin{smallmatrix} L_1 \\ L_2 \end{smallmatrix}\bigr)$,
with
$$\begin{pmatrix} A_1 & B \\ 0 & A_2 \end{pmatrix} \times
\begin{pmatrix} L_1 \\ L_2 \end{pmatrix} \to
\begin{pmatrix} L_1 \\ L_2 \end{pmatrix}\co
\biggl(\begin{pmatrix} a_1 & b \\ 0 & a_2 \end{pmatrix},
\begin{pmatrix} x_1 \\ x_2 \end{pmatrix}\biggr) \to
\begin{pmatrix} a_1x_1+\lambda(b\otimes x_2) \\ a_2x_2 \end{pmatrix} .$$
{\rm (ii)}\qua An $A$--module morphism $(f_1,f_2)\co (L_1,L_2,\lambda) \to
(L'_1,L'_2,\lambda')$ is defined by an $A_1$--module morphism $f_1\co L_1 \to L'_1$,
and an $A_2$--module morphism $f_2\co L_2 \to L'_2$ such that the diagram
$$\xymatrix{
B\otimes_{A_2}L_2 \ar[r]^-{\di{\lambda}}\ar[d]_-{\di{1\otimes f_2}}
& L_1 \ar[d]^-{\di{f_1}}\\
B\otimes_{A_2}L'_2 \ar[r]^-{\di{\lambda'}}
& L'_1}$$
commutes.

{\rm (iii)}\qua An $A$--module $L = (L_1,L_2,\lambda)$ is
f.g.~projective if and only if $\lambda$ is injective,
$\coker(\lambda)$ is a f.g.~projective $A_1$--module, and $L_2$ is
a f.g.~projective $A_2$--module.

{\rm (iv)}\qua The projection
$$\Proj(A) \to \Proj(A_1) \times \Proj(A_2);~
L = (L_1,L_2,\lambda) \mapsto (\coker(\lambda),L_2)$$
induces isomorphisms
$$K_*(A)\cong~ K_*(A_1) \oplus K_*(A_2).$$
{\rm (v)}\qua If an $A$--module $L = (L_1,L_2,\lambda)$ is h.d.~1 then
\begin{enumerate}
\item the 1--dimensional $A_1$--module chain complex
$$\xymatrix{K\co \cdots \ar[r] & 0 \ar[r] &
B\otimes_{A_2}L_2 \ar[r]^-{\di{\lambda}}& L_1}$$
is such that there exists a quasi-isomorphism ( =  homology equivalence)
$J \to K$ for a 1--dimensional f.g.~projective $A_1$--module chain complex $J$,
and
\item $L_2$ is an h.d.~1 $A_2$--module.
\end{enumerate}
If $B$ is a flat right $A_2$--module the converse also holds:
an $A$--module $L$ is h.d.~1 if and only if conditions 1.  and 2.
are satisfied.

{\rm (vi)}\qua The columns of $A$ are f.g.~projective $A$--modules
$$S_1 = (A_1,0,0),~S_2 = (B,A_2,1)$$
with
$$S_1\oplus S_2 = A,~{\rm End}(S_1) = A_1,~{\rm End}(S_2) = A_2.$$
The universal localization of $A$ inverting a non-empty subset
$\Sigma \subseteq \Hom_A(S_1,S_2)$ is a morphism of $2 \times 2$ matrix rings
$$A = \begin{pmatrix} A_1 & B \\ 0 & A_2 \end{pmatrix} \to
\Sigma^{-1}A = M_2(C) = \begin{pmatrix} C & C \\ C & C \end{pmatrix}$$
with $C$ the endomorphism ring of the induced f.g.~projective
$\Sigma^{-1}A$--module
$$\Sigma^{-1}S_1\cong \Sigma^{-1}S_2.$$
The composite of the functor
$$\Sigma^{-1}\co \Mod(A) \to \Mod(\Sigma^{-1}A);~P \mapsto
\Sigma^{-1}P = \Sigma^{-1}A\otimes_AP$$
and the Morita equivalence of categories
\begin{eqnarray*}
(C~C)\otimes_{\Sigma^{-1}A}-~:~
\Mod(\Sigma^{-1}A) &\xymatrix{\ar[r]^-{\di{\simeq}}&}& \Mod(C); \\
L = (L_1,L_2,\lambda) &\longmapsto& (C~C)\otimes_{\Sigma^{-1}A}L
\end{eqnarray*}
is the {\rm assembly} functor
\begin{align*}
\Mod(A) &\longrightarrow \Mod(C);\\
L = (L_1,L_2,\lambda) &\longmapsto
(C~C)\otimes_{\Sigma^{-1}A}\Sigma^{-1}L = (C~C)\otimes_AL\\
& = \coker\biggl(\begin{pmatrix} 1 \otimes \lambda \\ \kappa \otimes 1
\end{pmatrix}\co C \otimes_{A_1}B\otimes_{A_2}L_2 \to C \otimes_{A_1}L_1 \oplus
C\otimes_{A_2}L_2 \biggr)
\end{align*}
with
$$\kappa\co C \otimes_{A_1}B \to C;~x\otimes y \mapsto xy$$
the $(C,A_2)$--bimodule morphism defined by multiplication in $C$,
using the $A_1$--module morphism $B \to C$.
The assembly functor $\Proj(A) \to \Proj(C)$ induces the morphisms
$$\Sigma^{-1}\co K_*(A) = K_*(A_1)\oplus K_*(A_2) \to
K_*(\Sigma^{-1}A) = K_*(C).$$
{\rm (vii)}\qua If $B$ and $C$ are flat $A_1$--modules and $C$ is a flat
$A_2$--module then the $A$--module
$\bigl(\begin{smallmatrix}C \\ C \end{smallmatrix}\bigr)$
has a 1--dimensional flat $A$--module resolution
$$0 \to \begin{pmatrix} B \\ 0 \end{pmatrix}\otimes_{A_2}C
\to \begin{pmatrix} A_1 \\ 0 \end{pmatrix}\otimes_{A_1}C \oplus
\begin{pmatrix} B \\ A_2 \end{pmatrix}\otimes_{A_2}C
\to \begin{pmatrix} C \\ C \end{pmatrix} \to 0$$
so that $\Sigma^{-1}A = \bigl(\begin{smallmatrix} C \\ C
\end{smallmatrix}\bigr) \oplus \bigl(\begin{smallmatrix} C \\ C
\end{smallmatrix}\bigr)$ is stably flat.
An h.d.~1 $A$--module $L = (L_1,L_2,\lambda)$ is $\Sigma$--torsion if and only if
the $C$--module morphism
$$\begin{pmatrix} 1 \otimes \lambda \\ \kappa \otimes 1\end{pmatrix}\co 
C \otimes_{A_1}B\otimes_{A_2}L_2 \to C \otimes_{A_1}L_1 \oplus
C\otimes_{A_2}L_2$$
is an isomorphism.
\end{proposition}
\begin{proof}
(i) and (ii)\qua Standard.

(iii)\qua For any $A$--module $L = (L_1,L_2,\lambda)$ there is defined an
exact sequence
$$0 \to (\ker(\lambda),0,0)\to
(B\otimes_{A_2}L_2,L_2,1) \xymatrix{\ar[r]^-{\di{(\lambda,1)}} &}
(L_1,L_2,\lambda) \to (\coker(\lambda),0,0) \to 0.$$
Now $(A_1,0,0) = \bigl(\begin{smallmatrix} A_1 \\ 0 \end{smallmatrix}\bigr)$ and
$(B,A_2,1) = \bigl(\begin{smallmatrix} B \\ A_2 \end{smallmatrix}\bigr)$ are
f.g.~projective $A$--modules, since
$$(A_1,0,0) \oplus (B,A_2,1) = 
\biggl(A_1\oplus B,A_2,\begin{pmatrix} 0 \\ 1\end{pmatrix}\biggr) = A.$$
If $\ker(\lambda) = 0$ and $\coker(\lambda)$ is a f.g.~projective
$A_1$--module then
$(\coker(\lambda),0,0) = (A_1,0,0)\otimes_{A_1} \coker(\lambda)$
is a f.g.~projective $A$--module. If $L_2$ is a f.g.~projective $A_2$--module
then
$$(B\otimes_{A_2}L_2,L_2,1) = \begin{pmatrix} B \\ A_2 \end{pmatrix}\otimes_{A_2}L_2$$
is a f.g.~projective $A$--module. Thus if these two conditions are satisfied
then the exact sequence splits and $L$ is a f.g.~projective $A$--module.

Conversely, suppose that $(L_1,L_2,\lambda)$ is a f.g.~projective
$A$--module, so that there exists an $A$--module $(L'_1,L'_2,\lambda')$
with an $A$--module isomorphism
$$(L_1,L_2,\lambda) \oplus (L'_1,L'_2,\lambda')\cong~
\biggl((A_1)^k\oplus B^k,(A_2)^k,\begin{pmatrix} 0 \\ 1\end{pmatrix}\biggr) = A^k$$
for some $k \geqslant 0$. It follows from $\ker(\lambda \oplus
\lambda') = 0$ that $\ker(\lambda) = 0$, and from
$\coker(\lambda \oplus \lambda')\cong (A_1)^k$ that
$\coker(\lambda)$ is a f.g.~projective $A_1$--module.
Also, $L_2 \oplus L'_2\cong (A_2)^k$, so that $L_2$ is a f.g.~projective
$A_2$--module.

(iv)\qua The result that the inclusion and projection
$$i = A_1 \times A_2 \to A,\quad j\co A \to A_1 \times A_2$$
induce inverse isomorphisms
$$K_*(A_1 \times A_2) = K_*(A_1) \oplus K_*(A_2) \xymatrix{\ar@<0.5ex>[r]^{{i_*}} &
\ar@<0.5ex>[l]^{{j_*}}} K_*(A)$$
was first obtained by Berrick and Keating \cite{BK}.
Here is a proof using Waldhausen $K$--theory. It is immediate from $ji = 1$ that
$$j_*i_* = 1\co K_*(A_1 \times A_2) \to K_*(A) \to K_*(A_1 \times A_2).$$
Every f.g.~projective $A$--module $L = (L_1,L_2,\lambda\co B\otimes_{A_2}L_2 \to L_1)$
fits into a natural short exact sequence of f.g.~projective $A$--modules
$$0 \to (B\otimes_{A_2}L_2,L_2,1) \xymatrix{\ar[r]^-{{(\lambda,1)}}&}
(L_1,L_2,\lambda) \to  (\coker(\lambda),0)\to 0.$$
The functors
$$\begin{array}{l}
F_1\co \Proj(A) \to \Proj(A);~ L \mapsto
A\otimes_{A_1}A_1\otimes_AL = (\coker(\lambda),0),\\[2ex]
F_2\co \Proj(A) \to \Proj(A);~ L \mapsto
A\otimes_{A_2}A_2\otimes_AL = (B\otimes_{A_2}L_2,L_2,1)
\end{array}$$
fit into a cofibration sequence
$$F_2 \to 1_{\Proj(A)} \to F_1,$$
and are such that
$$F_k\co K_*(A) \to K_*(A_k) \to K_*(A)\quad (k = 1,2).$$
Now apply the additivity theorem for Quillen $K$--theory
\cite[Proposition 1.3.2 (4)]{Wald3} to identify
$$i_*j_* = F_1 + F_2 = 1\co K_*(A) \to K_*(A),$$
so that $i_*$, $j_*$ are inverse isomorphisms.

(v)\qua If $L = (L_1,L_2,\lambda)$ is an h.d.~1 $A$--module there
exists a 1--dimensional f.g.~projective $A$--module resolution
$$\xymatrix@C+10pt{
0 \ar[r] & (P_1,P_2,f) \ar[r]^-{{(h_1,h_2)}}&
(Q_1,Q_2,g) \ar[r] & (L_1,L_2,\lambda) \ar[r] & 0,}$$
so that $\coker(f),\coker(g)$ are f.g.~projective $A_1$--modules
and $P_2,Q_2$ are f.g.~projective $A_2$--modules.
The 1--dimensional $A_1$--module chain complex
$$\xymatrix{K\co \cdots \ar[r] & 0 \ar[r] &
B\otimes_{A_2}L_2 \ar[r]^-{{\lambda}}& L_1}$$
and the 1--dimensional f.g.~projective $A_1$--module chain complex
$$J\co J_1 = \coker(f) \xymatrix{\ar[r]^-{{h_1}}&}
J_0 = \coker(g)$$
are related by a homology equivalence $J \to K$. Furthermore,
$L_2 = \coker(h_2)$ is an h.d.~1 $A_2$--module. Thus both conditions
1. and 2. are satisfied.

Conversely, suppose that $B$ is a flat right $A_2$--module and
that $L = (L_1,L_2,\lambda)$ is an $A$--module
such that there exists a homology equivalence $J \to K$ with $J$
a 1--dimensional f.g.~projective $A_1$--module chain complex and that
$L_2$ is an h.d.~1 $A_2$--module with a 1--dimensional
f.g.~projective $A_2$--module resolution
$$0 \to P_2 \to Q_2 \to L_2 \to 0.$$
There is induced a short exact sequence of $A_1$--modules
$$0 \to B\otimes_{A_2}P_2 \to B\otimes_{A_2}Q_2 \to B\otimes_{A_2}L_2 \to 0$$
and it follows from the 1--dimensional f.g.~projective $A$--module resolution of $L$
$$0\to (B\otimes_{A_2}P_2,P_2,1) \oplus (J_1,0,0) \to
(B\otimes_{A_2}Q_2,Q_2,1) \oplus (J_0,0,0) \to L \to 0$$
that $L$ is an h.d.~1 $A$--module.

(vi) and (vii)\qua See \cite[2.2]{RNLAT}.
\end{proof}

We shall actually be working with $(\mu+1) \times (\mu{+}1)$--matrix rings:

\begin{definition}
For any ring $A$ and $\mu \geqslant 1$ define
the triangular $(\mu+1)\times (\mu{+}1)$--matrix ring
$$T_\mu(A) = \begin{pmatrix} A & A \oplus A & A \oplus A & \ldots &
A\oplus A \\
0 & A & 0 & \ldots & 0 \\
0 & 0 & A & \ldots & 0 \\
\vdots & \vdots & \vdots & \ddots & \vdots \\
0 & 0 & 0 &\ldots & A
\end{pmatrix}.$$
\end{definition}

The ring $T_\mu(A)$ is the $A$--coefficient path algebra of the quiver
with vertices $0,1,\ldots,\mu$ and two arrows $s^+_i,s^-_i\co i \to 0$
for $i = 1,2,\ldots,\mu$.  A $T_\mu(A)$--module $L = (L_i,f^+_i,f^-_i)$
consists of $A$--modules $L_0,L_1,\ldots,L_\mu$ and $A$--module morphisms
$f^+_i,f^-_i\co L_i \to L_0$ ($1 \leqslant i \leqslant \mu$).

Let $S_0,S_1,\ldots,S_\mu$ be the $T_\mu(A)$--modules defined
by the columns of $T_\mu(A)$, so that
\begin{align*}
S_0 & =  (A,0,\ldots,0;0,\ldots,0),\\
S_i & =  (A\oplus A,0,\ldots,0,A,0,\ldots,0;0,\ldots,0,{\rm id.},0,\ldots,0)
\quad (1 \leqslant i \leqslant \mu).
\end{align*}
It follows from
$$S_0 \oplus S_1 \oplus \cdots \oplus S_\mu = T_\mu(A)$$
that each $S_i$ is a f.g.~projective $T_\mu(A)$--module.
Let $\sigma = \{s^+_i,s^-_i\}$ be the set of
f.g.~projective $T_\mu(A)$--module morphisms
$$s^+_i = ((1~0),0,\ldots,0),~s^-_i = ((0~1),0,\ldots,0)\co
S_i \to S_0 \quad (1 \leqslant i \leqslant \mu).$$

\begin{proposition}\label{Floc}
{\rm (i)}\qua The universal $\sigma$--inverting localization of $T_\mu(A)$
is given by the inclusion
$$\F\co T_\mu(A) \to \sigma^{-1}T_\mu(A) = M_{\mu+1}(A[F_\mu])$$
with $M_{\mu+1}(A[F_\mu])$ the ring of all $(\mu+1)\times (\mu{+}1)$--matrices
with entries in $A[F_\mu]$.
The universal localization $\F$ is both injective and stably flat.

{\rm (ii)}\qua The composite
$$\xymatrix@C+10pt{
\Mod(T_\mu(A)) \ar[r]^-{\di{\F}} &
\Mod(M_{\mu+1}(A[F_\mu])) \ar[r]^-{\rm Morita}_-{\approx}& \Mod(A[F_\mu])}$$
sends a $T_\mu(A)$--module $L = (L_i,f^+_i,f^-_i)$
to the  {\rm assembly} $A[F_\mu]$--module
$$M = \coker\Bigl((f_1^+z_1{-}f_1^-~\ldots~ f_\mu^+z_\mu{-}f_\mu^-)\co
\bigoplus^\mu_{i = 1}L_i[F_\mu] \to L_0[F_\mu]\Bigr).$$
{\rm (iii)}\qua A $T_\mu(A)$--module $L = (L_i,f^+_i,f^-_i)$ is
f.g.~projective if and only if $L_0,\ldots,L_\mu$ are f.g.~projective
$A$--modules and the $A$--module morphism
$$\begin{pmatrix}
f^+_1 & f^-_1 & f^+_2 & f^-_2 & \ldots & f_{\mu}^+ & f_{\mu}^-
\end{pmatrix}\co \bigoplus^\mu_{i = 1}L_i \oplus L_i \to  L_0$$
is a split injection. The projection
$$\Proj(T_\mu(A)) \to \prod\limits_{\mu+1}\Proj(A);~
(L_i,f^+_i,f^-_i) \mapsto (L_0,L_1,L_2,\ldots,L_\mu)$$
induces isomorphisms in algebraic $K$--theory
$$K_*(T_\mu(A))\cong~\bigoplus\limits_{\mu+1}K_*(A).$$
{\rm (iv)}\qua A $T_\mu(A)$--module $L = (L_i,f^+_i,f^-_i)$ is h.d.~1
$\sigma$--torsion if and only if $L_0,\ldots,L_\mu$ are f.g.~projective
$A$--modules and the $A[F_\mu]$--module morphism
$$\begin{pmatrix} f^+_1z_1-f^-_1 & f^+_2z_2-f^-_2 & \ldots &
f^+_\mu z_\mu -f^-_\mu
\end{pmatrix}\co \bigoplus^\mu_{i = 1}L_i[F_\mu]\to L_0[F_\mu]$$
is an isomorphism. A f.g.~projective Seifert $A$--module $(P,e,\{\pi_i\})$ is
primitive if and only if $(P,P_i,f^+_i,f^-_i)$ is an h.d.~1
$\sigma$--torsion $T_\mu(A)$--module. The functor
$$\Prim(A) \to H(T_\mu(A),\sigma);~(P,e,\{\pi_i\}) \mapsto
(P,P_i,e\pi_i,(e-1)\pi_i)$$
is an equivalence of exact categories, so that
$$\PRIM_*(A) = K_*(H(T_\mu(A),\sigma)).$$
The forgetful functor
\begin{multline*}
\Prim(A) \to\prod\limits_{2\mu}\Proj(A);\\
\biggl(P^+ \oplus P^-,\begin{pmatrix} e^{++} & e^{+-} \\ e^{-+} & e^{--}
\end{pmatrix},\{\pi^+_i\}\oplus \{\pi^-_i\}\biggr)
\mapsto (P^+_1,P^-_1,\ldots,P^+_\mu,P^-_\mu)
\end{multline*}
(defined using \fullref{characterize_primitives}) is split by
$$\prod\limits_{2\mu}\Proj(A) \to \Prim(A);\quad 
(P^+_1,P^-_1,\ldots,P^+_\mu,P^-_\mu) \mapsto (P^+ \oplus P^-,0,\{\pi^+_i\}
\oplus \{\pi^-_i\}).$$
The reduced $K$--groups defined by
$$\widetilde{\PRIM}_*(A) = 
\ker(\PRIM_*(A) \to \bigoplus\limits_{2\mu} K_*(A))$$
are such that
$$K_*(H(T_\mu(A),\sigma)) = 
\PRIM_*(A) = \bigoplus\limits_{2\mu}K_*(A)
\oplus\widetilde{\PRIM}_*(A).$$
\end{proposition}
\begin{proof} The universal localization $\sigma^{-1}T_{\mu}(A)$
is the $(\mu+1)\times (\mu{+}1)$--matrix
ring $M_{\mu+1}(R)$ with $R$ the endomorphism ring
of the induced f.g.~projective $\sigma^{-1}T_\mu(A)$--module $\sigma^{-1}S_0$,
and there is defined an isomorphism
$$A[F_\mu] \to R;~z_i \mapsto s^+_i(s^-_i)^{-1}.$$
The remaining parts are given by \fullref{triangular},
viewing the $(\mu+1)\times (\mu+1)$ matrix ring
$T_\mu(A)$ as a triangular $2 \times 2$ matrix ring
$$T_\mu(A) = \begin{pmatrix} A_1 & B \\ 0 & A_2 \end{pmatrix}$$
with
$$A_1 = A,\quad A_2 = \begin{pmatrix}
A & 0 & \ldots & 0 \\
0 & A & \ldots & 0 \\
\vdots & \vdots & \ddots & \vdots \\
0 & 0 &\ldots & A
\end{pmatrix},\quad B = \begin{pmatrix} A\oplus A & \ldots & A\oplus A\end{pmatrix}$$
such that
$$\Mod(A_2) =  \prod\limits_{\mu}\Mod(A).$$
An $A_2$--module is just a $\mu$--tuple $(L_1,L_2,\ldots,L_{\mu})$
of $A$--modules. By the $2 \times 2$ theory a $T_\mu(A)$--module
$L$ just a $(\mu{+}1)$--tuple $(L_0,L_1,\ldots,L_{\mu})$ of $A$--modules,
together with $A$--module morphisms $f^+_i,f^-_i\co L_i \to L_0$
($1 \leqslant i \leqslant \mu$).
Note that $B$ is a flat right $A_2$--module, and that for an h.d.~1
$\sigma$--torsion $T_\mu(A)$--module $L = (L_i,f^+_i,f^-_i)$
each $L_i$ ($0 \leqslant i \leqslant \mu$) is a f.g.~projective
$A$--module, by the following argument.
The necessary and sufficient conditions of \fullref{triangular} (v) and (vii)
for a $T_\mu(A)$--module $L$ to be h.d.~1 $\sigma$--torsion
are:
\begin{itemize}
\item[(i)] there exists a 1--dimensional f.g.~projective $A$--module chain
complex $J\co J_1 \to J_0$ with a homology equivalence
$$\xymatrix@C+20pt{J_1 \ar[r] \ar[d] & J_0 \ar[d] \\
\bigoplus\limits^\mu_{i = 1}L_i \oplus L_i \ar[r]^-{\di{(f^+_i~f^-_i)}}&
L_0,}$$
\item[(ii)]  each $L_i$ $(1 \leqslant i \leqslant \mu)$
is an h.d.~1 $A$--module,
\item[(iii)] the $A[F_\mu]$--module morphism
$$(f^+_iz_i-f^-_i)\co \bigoplus\limits^\mu_{i = 1}L_i[F_\mu]\to L_0[F_\mu]$$
is an isomorphism.
\end{itemize}
If $L$ satisfies these conditions there is defined a
commutative diagram of $A$--modules
$$\xymatrix@C+10pt{0 \ar[r] & \bigoplus\limits^\mu_{i = 1}L_i
\ar[r]^-{\bigl(\begin{smallmatrix} 1 \\ -1 \end{smallmatrix}\bigr)}
\ar[d]_{\di{\cong}}^-{{(f^+_i-f^-_i)}} &
\bigoplus\limits^\mu_{i = 1}L_i \oplus L_i \ar[d]^-{{(f^+_i~f^-_i) }}
\ar[r]^-{\di{(1~1)}} & \bigoplus\limits^{\mu}_{i = 1}L_i\ar[d] \ar[r]& 0
\\
0 \ar[r] & L_0\ar[r]^-{\di{1}} & L_0 \ar[r]& 0 \ar[r] & 0
}$$
with exact rows and with $(f^+_i-f^-_i)$ an isomorphism.
There are defined $A$--module isomorphisms
$$J_0 \oplus L_0\cong~ J_0 \oplus \bigoplus\limits^{\mu}_{i = 1}L_i\cong~J_1,$$
so that each $L_i$ ($0 \leqslant i \leqslant \mu$) is a f.g.~projective
$A$--module.
\end{proof}

\begin{example} The assembly of
$A[F_\mu]$--modules is an algebraic analogue of the geometric
construction of an $F_\mu$--cover $\wwtilde{W}$ of a space $W$
from a fundamental domain $U \subset \wwtilde{W}$. The subspaces
$$V_i = U \cap z^{-1}_iU,\quad z_iV_i = z_iU \cap U \subset
    U \quad (1 \leqslant i \leqslant \mu)$$
are disjoint, with embeddings
$$f^+_i\co V_i \to U;~x \mapsto x,\quad f^-_i\co V_i \to U;~x \mapsto
    z_ix \quad (1 \leqslant i \leqslant \mu),$$
and  $\wwtilde{W}$ can be
constructed by glueing together $F_\mu$ copies of $U$
\begin{align*}
\wwtilde{W}& = (F_\mu \times U)/\{(g,f^+_i(x)) \sim (gz_i,f^-_i(x))\,\vert\,
g \in F_\mu,x \in V_i,1 \leqslant i \leqslant \mu\}\\[1ex]
& = \bigcup\limits_{g \in F_\mu} gU~{\rm with}~U \cap z_i^{-1}U = V_i.
\end{align*}
Such a situation arises if $W$ is a manifold with
a surjection $\pi_1(W) \to F_\mu$, eg a boundary link exterior.
The surjection is induced by a map
$$c\co W \to BF_\mu = \bigvee\limits_\mu S^1$$
which is transverse regular at $\{1,2,\ldots,\mu\} \subset BF_\mu$.
Cutting $W$ open at the inverse image codimension--1 submanifolds
$V_i = c^{-1}(\{i\}) \subset W$ there is obtained a fundamental
domain $U \subset \wwtilde{W}$ for the pullback
$\wwtilde{W} = c^*EF_\mu$ to $W$ of the universal cover $EF_{\mu}$ of $BF_\mu$.
More generally, suppose that $W$ is a finite $CW$ complex with
an $F_\mu$--cover $\wwtilde{W}$, and that $U \subset \wwtilde{W}$
is a fundamental domain which is a subcomplex.  The embeddings
$f^+_i,f^-_i\co V_i\to U$ induce inclusions of the
cellular f.g.~free $\Z$--module chain complexes $f^+_i,f^-_i\co C(V_i)\to
C(U)$. The f.g.~projective $T_\mu(\Z)$--module chain complex
$C = (C(U),C(V_i),f^+_i,f^-_i)$ has assembly the cellular f.g.~free
$\Z[F_\mu]$--module chain complex of $\wwtilde{W}$
$$\coker\bigg((f_1^+z_1-f_1^- ~\ldots~f_\mu^+z_\mu-f_\mu^-)\co 
\bigoplus^\mu_{i = 1}C(V_i)[F_\mu] \to C(U)[F_\mu]\bigg) = C(\wwtilde{W}),$$
such that
$$C(\wwtilde{W})_r = \coker\bigg((f_1^+ ~\ldots~
f_\mu^+)\co \bigoplus^\mu_{i = 1}C(V_i)_r \to C(U)_r\bigg)[F_\mu].$$
\end{example}

\begin{thm} The algebraic $K$--groups of $A[F_\mu]$ split as
$$K_*(A[F_\mu]) = K_*(A) \oplus \bigoplus\limits_{\mu} K_{*-1}(A) \oplus
\widetilde{\PRIM}_{*-1}(A).$$
\end{thm}
\begin{proof} By \fullref{Floc} the universal localization
$$\F\co A[F_\mu]\to\sigma^{-1}T_\mu(A) = M_{\mu+1}(A)$$
is injective and stably flat.
The noncommutative localization exact sequence of Neeman and Ranicki
\cite{NR1,NR2}
\begin{multline*}
\cdots \longrightarrow K_{n+1}(\sigma^{-1}T_\mu(A)) \longrightarrow
K_n(H(T_\mu(A),\sigma)) \\
\longrightarrow K_n(T_\mu(A)) \longrightarrow
K_n(\sigma^{-1}T_\mu(A)) \longrightarrow \cdots
\end{multline*}
is given by
$$\cdots \to K_{n+1}(A[F_\mu]) \to \PRIM_n(A) \to
\bigoplus\limits_{\mu+1} K_n(A) \to K_n(A[F_\mu]) \to \cdots$$
with $\PRIM_n(A) \to K_n(T_\mu(A)) = \bigoplus\limits_{\mu+1}K_n(A)$
induced by
$$\Prim(A) \to \prod\limits_{\mu+1} \Proj(A);~
(P,e,\{\pi_i\}) \mapsto (P,P_1,P_2,\ldots,P_\mu),$$
so that
\begin{align*}
&\PRIM_n(A)  =  \bigoplus\limits_{2\mu} K_n(A) \oplus
\widetilde{\PRIM}_n(A) \to \bigoplus\limits_{\mu+1}K_n(A);\\[1ex]
&(x^+_1,x^-_1,x^+_2,x^-_2,\ldots,x^+_\mu,x^-_\mu,\wtilde{x})\\
&\hspace{100pt}\longmapsto \Bigl(\sum\limits^\mu_{i = 1}(x^+_i+x^-_i),x^+_1+x^-_1,
x^+_2+x^-_2,\ldots,x^+_\mu+x^-_\mu\Bigr).
\end{align*}
This completes the proof.
\end{proof}

\begin{definition}{\rm
Let $\G\co A[F_\mu]\to\Sigma^{-1}A[F_\mu]$ be the universal
localization inverting the set $\Sigma$ of morphisms of f.g.~projective
$A[F_\mu]$--modules which induce an isomorphism of f.g.~projective
$A$--modules under the augmentation
$\epsilon\co A[F_\mu] \to A;z_i \mapsto 1$.}
\end{definition}

\begin{proposition} \label{Delta}
{\rm (i)}\qua The universal localization
$\G\co A[F_\mu]\to\Sigma^{-1}A[F_\mu]$ is injective.
The h.d.~1 $\Sigma$--torsion $A[F_\mu]$--module category is
$$H(A[F_\mu],\Sigma) = \Bla(A).$$
{\rm (ii)}\qua The composite
$$\G\F\co T_\mu(A) \stackrel{\F}{\longrightarrow}
\sigma^{-1}T_\mu(A) = M_{\mu+1}(A[F_\mu])
\stackrel{\G}{\longrightarrow}
\tau^{-1}T_\mu(A) = M_{\mu+1}(\Sigma^{-1}A[F_\mu])$$
is the universal localization inverting the set $\tau$ of morphisms of
f.g.~projective $T_\mu(A)$--modules which become isomorphisms under the
composite
$$\xymatrix{\epsilon \F\co 
T_\mu(A) \ar[r]^-{\di{\F}} & \sigma^{-1}T_\mu(A) = 
M_{\mu+1}(A[F_\mu]) \ar[r]^-{\di{\epsilon}} &     M_{\mu+1}(A).}$$
{\rm (iii)}\qua A $T_\mu(A)$--module $L = (L_i,f^+_i,f^-_i)$ is h.d.~1
$\tau$--torsion if and only if $L_0,\ldots,L_\mu$ are f.g.~projective
$A$--modules and the $A$--module morphism
$$f = \begin{pmatrix} f^+_1-f^-_1 & f^+_2-f^-_2 & \ldots & f^+_\mu-f^-_\mu
\end{pmatrix}\co L_1 \oplus L_2 \oplus \cdots \oplus L_\mu \to L_0$$
is an isomorphism, if and only if
$$(P,e,\{\pi_i\}) = \Bigl(\bigoplus\limits_{i = 1}^\mu L_i,
f^{-1}(f^+_1~f^+_2~\ldots~f^+_\mu),\{\pi_i\}\Bigr)$$
is a f.g.~projective Seifert $A$--module. The functor
$$\Sei(A) \to H(T_\mu(A),\tau);~(P,e,\{\pi_i\}) \mapsto
(P,P_i,e\pi_i,(e-1)\pi_i)$$
is an equivalence of exact categories. The assembly of
$(L_i,f^+_i,f^-_i)$ is the covering Blanchfield $A[F_\mu]$--module of
$(P,e,\{\pi_i\})$
\begin{multline*}
\coker\Bigl((f_1^+z_1-f_1^-~\ldots~
f_\mu^+z_\mu-f_\mu^-)\co \bigoplus\limits^\mu_{i = 1}L_i[F_\mu] \to
L_0[F_\mu]\Bigr)\\
 = \coker\big(1-e+ze\co P[F_\mu] \to P[F_\mu]\big)
 = B(P,e,\{\pi_i\}),
\end{multline*}
so that up to equivalence
$$\F = B\co H(T_\mu(A),\tau) = \Sei(A) \to
H(M_{\mu+1}(A[F_\mu]),\tau) = \Bla(A).$$
{\rm (iv)}\qua The forgetful functor
$$\Sei(A) \to\prod\limits_{\mu}\Proj(A) ;~
(P,e,\{\pi_i\}) \mapsto (P_1,P_2,\ldots,P_\mu)$$
is split by
$$\prod\limits_{\mu}\Proj(A) \to \Prim(A);~
(P_1,P_2,\ldots,P_\mu) \mapsto
  \Bigl(\bigoplus\limits_{i = 1}^{\mu}P_i,0,\{\pi_i\}\Bigr)$$
The reduced $K$--groups defined by
$$\widetilde{\SEI}_*(A) = 
\ker\Bigl(\SEI_*(A) \to \bigoplus\limits_{\mu} K_*(A)\Bigr)$$
are such that
$$K_*(H(T_\mu(A),\tau)) = \SEI_*(A) = \bigoplus\limits_{\mu}K_*(A)
\oplus\widetilde{\SEI}_*(A).$$
\end{proposition}
\begin{proof} (i)\qua The Magnus--Fox embedding
$A[F_\mu] \to A\llangle x_1,
\ldots,x_\mu \rrangle$ is $\Sigma$--inverting, so that there is
a unique factorization
$$A[F_{\mu}] \to \Sigma^{-1}A[F_\mu] \to A\llangle x_1,x_2,
\ldots,x_\mu \rrangle.$$
The identification $H(A[F_\mu],\Sigma) = \Bla(A)$ is a formality,
as is the identification $\Proj(A[F_\mu]) = \PP(\Proj_A(A[F_\mu]))$
with $\Proj_A(A[F_\mu]) \subseteq \Proj(A[F_\mu])$ the full subcategory
with objects isomorphic to the f.g.~projective $A[F_\mu]$--modules $P[F_\mu]$ induced
from f.g.~projective $A$--modules $P$.

(ii)--(iv)\qua By construction, working as in the proof of
\fullref{Floc} (iv) to show that if $L = (L_i,f^+_i,f^-_i)$ is an h.d.~1
$\tau$--torsion $T_\mu(A)$--module then $L_0,L_1,\ldots,L_\mu$ are
f.g.~projective $A$--modules.
\end{proof}

\begin{thm} \label{final}
{\rm (i)}\qua The algebraic $K$--groups of $\Prim(A)$, $\Sei(A)$ and
$\Bla(A)$ fit into a commutative braid of exact sequences
$$
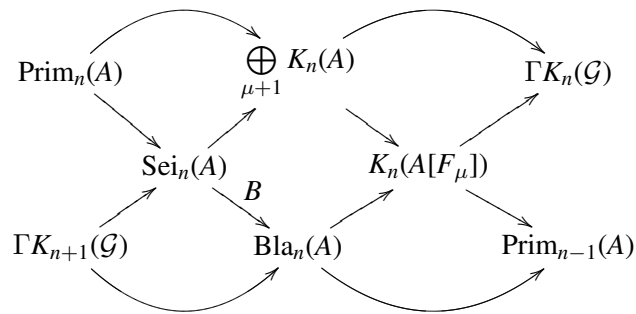
\xymatrix@C-25pt@R-10pt{
\PRIM_n(A)\ar[dr] \ar@/^2pc/[rr]^{} &&
\bigoplus\limits_{\mu+1}K_n(A) \ar[dr] \ar@/^2pc/[rr]
&&\Gamma K_n(\G)  \\&
\SEI_n(A)\ar[ur] \ar[dr]^-{\di{B}} && K_n(A[F_\mu]) \ar[ur] \ar[dr]&&\\
\Gamma K_{n+1}(\G)  \ar[ur]
\ar@/_2pc/[rr]_-{}&&\BLA_n(A)
\ar[ur]\ar@/_2pc/[rr]_{}&&\PRIM_{n-1}(A)}$$
for $n \in \Z$, with $\G\co A[F_\mu] \to \Sigma^{-1}A[F_\mu]$ the
universal localization and
\begin{align*}
K_*(T_\mu(A))& = \bigoplus\limits_{\mu+1}K_*(A),\\[-1ex]
K_*(H(T_\mu(A),\sigma))& = (\SEI,\BLA)_*(A) = \PRIM_*(A) = 
\bigoplus\limits_{2\mu}K_*(A)\oplus \widetilde{\PRIM}_*(A),\\[-1ex]
K_*(H(T_\mu(A),\tau))& = \SEI_*(A) = \bigoplus\limits_{\mu}K_*(A)
\oplus\widetilde{\SEI}_*(A),\\[-1ex]
K_*(H(A[F_\mu],\Sigma))& = \BLA_*(A) = \bigoplus\limits_{\mu}K_{*-1}(A)
\oplus \widetilde{\BLA}_*(A),\\[-1ex]
\Gamma K_*(\G)& = K_*(A) \oplus \widetilde{\SEI}_{*-1}(A)~( = K_*(\Sigma^{-1}A[F_\mu])~
\emph{for}~* \leqslant 1)
\end{align*}
The reduced $K$--groups fit into a long exact sequence
$$\cdots \to \widetilde{\PRIM}_n(A) \to
\widetilde{\SEI}_n(A) \to \widetilde{\BLA}_n(A)
\to \widetilde{\PRIM}_{n-1}(A) \to \cdots.$$
{\rm (ii)}\qua If $\G\co A[F_\mu]\to\Sigma^{-1}A[F_\mu]$ is stably flat then
$$\Gamma K_n(\G) = K_n(\Sigma^{-1}A[F_\mu]) = 
K_n(A) \oplus \widetilde{\SEI}_{n-1}(A)$$
for all $n \in \Z$.
\end{thm}
\begin{proof}
(i)\qua Consider the commutative square of Waldhausen categories
$$\xymatrix{
(\C^b(T_\mu(A)),w_{T_{\mu}(A)}) \ar[r] \ar[d]_-{\di{\F}} &
(\C^b(T_\mu(A)),w_{\tau^{-1}T_{\mu}(A)}) \ar[d] \\
(\C^b(T_\mu(A)),w_{\sigma^{-1}T_{\mu}(A)})\ar[r]^-{\di{\G}} &
(\C^b(\sigma^{-1}T_\mu(A)),w_{\tau^{-1}T_{\mu}(A)}).}$$
Since $\F\co T_\mu(A)\to\sigma^{-1}T_\mu(A) = M_{\mu+1}(A[F_\mu])$ is stably flat
there are defined equivalences
$$(\C^b(T_\mu(A)),w_{\sigma^{-1}T_{\mu}(A)})~\approx~
(\C^b(\sigma^{-1}T_\mu(A)),w_{\sigma^{-1}T_{\mu}(A)})~\approx~
(\C^b(A[F_\mu]),w_{A[F_\mu]})$$
which induce  homotopy equivalences
$$K\PP(\C^b(T_\mu(A)),w_{\sigma^{-1}T_{\mu}(A)})~\simeq~K(\sigma^{-1}T_\mu(A))~
\simeq~K(A[F_\mu]).$$
Also, since $\tau^{-1}T_\mu(A) = M_{\mu+1}(\Sigma^{-1}A[F_\mu])$
the functor
$$(\C^b(T_\mu(A)),w_{\tau^{-1}T_{\mu}(A)}) \to
(\C^b(\sigma^{-1}T_\mu(A)),w_{\tau^{-1}T_{\mu}(A)})$$
induces an equivalence of the homotopy categories
$$D(\C^b(T_\mu(A)),w_{\tau^{-1}T_{\mu}(A)})~\approx~
(\C^b(\sigma^{-1}T_\mu(A),w_{\tau^{-1}T_{\mu}(A)}).$$
The composite of this equivalence and the Morita equivalence
$$D(\C^b(\sigma^{-1}T_\mu(A)),w_{\tau^{-1}T_{\mu}(A)})~\approx~
D(\C^b(A[F_\mu]),w_{\Sigma^{-1}A[F_\mu]})$$
induces a homotopy equivalence
$$\begin{array}{ll}
K\PP(\C^b(\sigma^{-1}T_\mu(A)),w_{\tau^{-1}T_{\mu}(A)})~
\simeq&K\PP(\C^b(A[F_\mu]),w_{\Sigma^{-1}A[F_\mu]})\\[1ex]
& = \Gamma K(\G\co A[F_\mu] \to\Sigma^{-1}A[F_\mu]).
\end{array}$$
Thus Propositions \ref{Gamma}, \ref{Floc} and \ref{Delta}
give a braid of Waldhausen categories
$$\xymatrix@C-85pt{
(\C^b(T_\mu(A),\sigma),w_{(T_{\mu}(A),\sigma)})\ar[dr] \ar@/^2pc/[rr]^{} &&
(\C^b(T_\mu(A)),w_{T_{\mu}(A)})\quad \ar[dr] \ar@/^2pc/[rr]
&&(\C^b(A[F_\mu]),w_{\Sigma^{-1}A[F_\mu]})\\
&(\C^b(T_\mu(A),\tau),w_{(T_{\mu}(A),\tau)})\ar[ur] \ar[dr]&&
(\C^b(A[F_\mu]),w_{A[F_\mu]}) \ar[ur] &&\\
&&(\C^b(A[F_\mu],\Sigma),w_{\C^b(A[F_\mu],\Sigma)})\ar[ur]&&}$$
inducing a commutative braid of exact sequences
$$
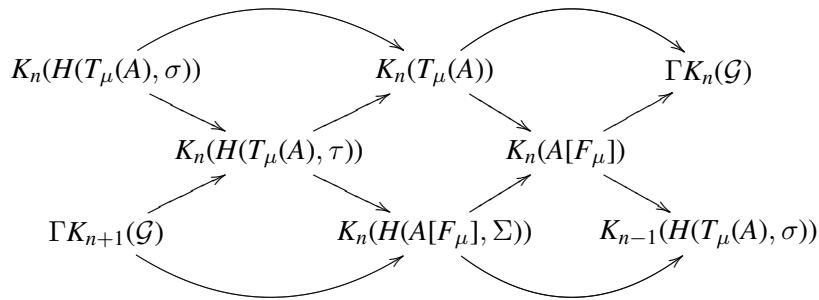
\xymatrix@C-40pt@R-10pt{
K_n(H(T_\mu(A),\sigma))\ar[dr] \ar@/^2pc/[rr]^{} &&
K_n(T_\mu(A)) \ar[dr] \ar@/^2pc/[rr]&&\Gamma K_n(\G)  \\&
K_n(H(T_\mu(A),\tau))\ar[ur] \ar[dr] && K_n(A[F_\mu]) \ar[ur] \ar[dr]&&\\
\Gamma K_{n+1}(\G)  \ar[ur] \ar@/_2pc/[rr]_-{}&&
K_n(H(A[F_\mu],\Sigma))
\ar[ur]\ar@/_2pc/[rr]_{}&&K_{n-1}(H(T_\mu(A),\sigma))}$$
Split off the reduced $K$--groups in
\begin{align*}
\PRIM_*(A)& =  \bigoplus\limits_{2\mu}K_*(A) \oplus \widetilde{\PRIM}_*(A),\\
\SEI_*(A)& =  \bigoplus\limits_{\mu}K_*(A) \oplus \widetilde{\SEI}_*(A)
\end{align*}
from the long exact sequence
$$\cdots \to \PRIM_n(A) \to\SEI_n(A) \to \BLA_n(A)\to \PRIM_{n-1}(A) \to
\cdots$$
to define the reduced $K$--groups in
$$\BLA_*(A) = \bigoplus\limits_{\mu}K_{*-1}(A) \oplus \widetilde{\BLA}_*(A)$$
and to obtain the long exact sequence
$$\cdots \to \widetilde{\PRIM}_n(A) \to
\widetilde{\SEI}_n(A) \to \widetilde{\BLA}_n(A)
\to \widetilde{\PRIM}_{n-1}(A) \to \cdots.$$
(ii)\qua This is a special case of \fullref{Gamma} (ii).
\end{proof}

This completes the proofs of Theorems \ref{thm4} and \ref{thm5} of the
Introduction.

\begin{remark}
Unfortunately, we do not know if the universal localization
$\Sigma^{-1}A[F_\mu]$ is stably flat in general.
See Dicks and Sontag \cite{DS}, Farber and Vogel \cite{FV}
for proofs that $\Sigma^{-1}A[F_\mu]$ is stably flat
when $A$ is a principal ideal domain, and
Ara and Dicks \cite[Theorem 4.4]{AD} when $A$ is a von Neumann regular
ring or a commutative Bezout domain.
\end{remark}

\begin{remark}
\label{det}
Sheiham \cite{Sh3} computed
$$K_1(\Sigma^{-1}A[F_\mu]) = K_1(A) \oplus \epsilon_\Sigma^{-1}(1)/C$$
with $\epsilon_\Sigma\co \Sigma^{-1}A[F_\mu] \to A$ the
factorization of the augmentation map $\epsilon\co A[F_\mu] \to A$ and
$C \subseteq \epsilon_{\Sigma}^{-1}(1)$ the subgroup generated by the
commutators
$$(1+ab)(1+ba)^{-1}\qquad(a,b \in \Sigma^{-1}A[F_\mu],~
\epsilon(ab) = \epsilon(ba) = 0).$$
It follows from the splitting given by \fullref{final} (i)
$$K_1(\Sigma^{-1}A[F_\mu]) = K_1(A) \oplus \widetilde{\SEI}_0(A)$$
that there is defined an isomorphism
$$\widetilde{\SEI}_0(A) \xymatrix{\ar[r]^-{{\cong}}&}
\epsilon_{\Sigma}^{-1}(1)/C;~
(P,e,\{\pi_i\}) \mapsto D(1-e+ez\co P[F_\mu] \to P[F_\mu])$$
with $D$ the generalized Dieudonn\'e noncommutative determinant of
\cite[4.3]{Sh3}.
\end{remark}

\begin{example}
{\rm (i)}\qua The algebraic $K$--groups of $\Z[F_\mu]$ are such that
$$\begin{array}{l}
K_*(\Z[F_\mu]) = K_*(\Z) \oplus \bigoplus\limits_\mu K_{*-1}(\Z),\\[1ex]
K_n(\Z[F_\mu]) = K_n(\Z) = \begin{cases}
\Z&{\rm if}~n = 0\\
0&{\rm if}~n \leqslant -1
\end{cases}
\end{array}$$
by Stallings \cite{St}, Gersten \cite{Ge}, Bass \cite[XII]{B2} and
Waldhausen \cite{Wald2,Wald2a}, so that
\begin{align*}
&\Flk(\Z)  =  \Bla(\Z), \qquad \widetilde{\PRIM}_*(\Z)  =  0,\\
&K_{*+1}(\Sigma^{-1}\Z[F_\mu])/K_{*+1}(\Z)  = 
\widetilde{\SEI}_*(\Z)  =  \FLK_*(\Z)  =  \BLA_*(\Z),\\
&K_*(H(\Z[F_\mu],\Sigma))  =  \bigoplus\limits_{\mu}K_{*-1}(\Z)\oplus
\widetilde{\SEI}_*(\Z),\\
&K_n(H(\Z[F_\mu],\Sigma))  =  \widetilde{\SEI}_n(\Z)\quad(n \leqslant 0).
\end{align*}
{\rm (ii)}\qua Given a $\mu$--component boundary link
$\ell\co \bigsqcup_\mu S^n \subset S^{n+2}$ with exterior $W$
and given a $\mu$--component Seifert surface $V = V_1 \sqcup V_2
\sqcup \ldots \sqcup V_\mu \subset S^{n+2}$ for $\ell$
let $\dot C(W)$, $(\dot
C(V),e,\{\pi_i\})$ be the chain complexes defined in \fullref{chain}. Thus $\dot C(\wwtilde{W})$ is a
$\Sigma^{-1}\Z[F_\mu]$--acyclic $(n{+}2)$--dimensional f.g.~free
$\Z[F_\mu]$--module chain complex, $(\dot C(V),e,\{\pi_i\})$ is an
$(n{+}1)$--dimensional chain complex in $\Sei(\Z)$, and $B(\dot
C(V),e,\{\pi_i\})$ is an $(n{+}1)$--dimensional chain complex in
$\Flk(\Z)$ with a homology equivalence $\dot C(\wwtilde{W}) \to
B(\dot C(V),e,\{\pi_i\})$. The torsion
\begin{align*}
\tau(\ell)& =  \tau(\Sigma^{-1}\dot C(\wwtilde{W}))\\
& =  (\dot C(V),e,\{\pi_i\})
 = \sum\limits^{n+1}_{r = 0}(-)^r(\dot C_r(V),e,\{\pi_i\})  =  [\dot C(\wwtilde{W})]\\
&\hspace*{40pt}
\in K_1(\Sigma^{-1}\Z[F_\mu])/K_1(\Z)  =  K_0(H(\Z[F_\mu],\Sigma))  = 
\widetilde{\rm \SEI}_0(\Z)  =  \BLA_0(\Z)
\end{align*}
is an isotopy invariant of $\ell$, given by
Sheiham \cite{Sh3} to be the generalized Dieudonn\'e determinant
$$\tau(\ell) = 
\sum\limits^{n+1}_{r = 0}(-)^r
D(1-e+ez\co \dot C_r(V)[F_\mu] \to \dot C_r(V)[F_\mu])
\in \widetilde{\SEI}_0(\Z) = \epsilon_{\Sigma}^{-1}(1)/C$$
with $\epsilon_{\Sigma}\co \Sigma^{-1}\Z[F_\mu] \to \Z$ and
$C \subseteq \epsilon_{\Sigma}^{-1}(1)$ as recalled in \fullref{det}.
The $\Z[F_\mu]$--modules
$\dot H_r(\wwtilde{W})/\Z\hbox{\rm -torsion}$
($0 \leqslant r \leqslant n+1$) are h.d.~1 $F_\mu$--link modules, and
\begin{align*}
\tau(\ell)& = \sum\limits^{n+1}_{r = 0}(-)^r
D(1-e+ez\co \dot H_r(V)[F_\mu] \to \dot H_r(V)[F_\mu])\\
& = \sum\limits^{n+1}_{r = 0}(-)^r[\dot H_r(\wwtilde{W})/\Z\hbox{\rm--torsion}] \\
&\in
K_1(\Sigma^{-1}\Z[F_\mu])/K_1(\Z) = K_0(H(\Z[F_\mu],\Sigma))\\
&\hskip75pt
 = \widetilde{\rm \SEI}_0(\Z) = \BLA_0(\Z) = 
(\Sigma^{-1}\Z[F_\mu])^{\bullet}/\{\pm 1\}.
\end{align*}
For $\mu = 1$ this is just the Reidemeister torsion of a knot $\ell\co S^n
\subset S^{n+2}$, which is the alternating product of the Alexander
polynomials
$$\begin{array}{ll}
\tau(\ell)& = \sum\limits^{n+1}_{r = 0}(-)^r
{\rm det}(1-e+ez\co \dot H_r(V)[z,z^{-1}] \to \dot H_r(V)[z,z^{-1}])\\
& = \sum\limits^{n+1}_{r = 0}(-)^r[\dot H_r(\wwtilde{W})/\Z\hbox{\rm -torsion}]
\\
&\hskip20pt\in
K_1(\Sigma^{-1}\Z[z,z^{-1}])/K_1(\Z) = K_0(H(\Z[z,z^{-1}],\Sigma))\\
&\hskip40pt
 = \widetilde{\rm \SEI}_0(\Z) = \BLA_0(\Z) = \widetilde{\ENDD}_0(\Z)
 = (\Sigma^{-1}\Z[z,z^{-1}])^{\bullet}/\{\pm 1\}
\end{array}$$
(Milnor \cite{M}, cf \cite[Example 17.11]{RHK}).

(iii)\qua The isotopy classes of simple $\mu$--component boundary links
$\ell\co \bigsqcup_{\mu}S^{2q-1} \subset S^{2q+1}$ for $q
\geqslant 3$ are in one-one correspondence with the
`$l$--equivalence classes of Seifert matrices' (Liang \cite{Li},
generalizing the case $\mu = 1$ due to Levine \cite{L1}), and also
with the `$R$--equivalence classes of $(-)^q$--symmetric isometry
structures of multiplicity $\mu$' (Farber \cite[4.7]{Fa2}).  For
simple $\ell$ $H_q(\wwtilde{W})$ is an h.d.~1 $F_\mu$--link
module, and the torsion
$$\tau(\ell) = (-)^q[H_q(\wwtilde{W})] \in \widetilde{\SEI}_0(\Z) = 
\FLK_0(\Z) = \BLA_0(\Z)$$
is just the $K$--theory part of these complete isotopy invariants for
$q \geqslant 3$.\end{example}

\bibliographystyle{gtart}
\bibliography{link}

\end{document}